\documentclass[a4paper,12pt]{amsart}
\usepackage{amssymb}%,amsthm,amsfonts,amssymb,euscript,mathrsfs}
\usepackage{enumitem}
\usepackage{appendix}
\usepackage{url}
\usepackage[ansinew]{inputenc}
\usepackage[english]{babel}%%Francais par defaut
\usepackage{color}
\usepackage{graphicx}
\usepackage{tikz}
\usetikzlibrary{calc}
\usepackage[normalem]{ ulem }
\usepackage{soul}

\ifx \undefined \cftil \def \cftil#1{\~#1}\fi

\newcommand{\s}{{\mathfrak s}}

 \DeclareMathAlphabet{\got}{U}{euf}{m}{n}     %% les gothiques
\DeclareMathAlphabet{\mat}{U}{msb}{m}{n}     %% les lettres Bbb=20
\DeclareMathAlphabet{\mathbold}{OML}{cmm}{bx}{it}
\newtheorem{theo}{Theorem}[section]
\newtheorem{pr}{Proposition}[section]
\newtheorem{lem}{Lemma}[section]
\newtheorem{defi}{Definition}[section]

\newtheorem{co}{Corollary}[section]
\newtheorem{rem}{Remark}[section]

\newtheorem{qu}{Question}

\newenvironment{dem}{\vskip 1em{\it Proof} :}%
{\unskip\hfill\null\nobreak\hfill\carre\vskip1em\par}
\newcommand{\carre}{\rule{1ex}{1ex}} 

\makeatletter
\providecommand*{\diff}%
{\@ifnextchar^{\DIfF}{\DIfF^{}}}
\def\DIfF^#1{%
\mathop{\mathrm{\mathstrut d}}%
\nolimits^{#1}\gobblespace}
\def\gobblespace{%
\futurelet\diffarg\opspace}
\def\opspace{%
\let\DiffSpace\!%
\ifx\diffarg(%
\let\DiffSpace\relax
\else
\ifx\diffarg[%
\let\DiffSpace\relax
\else
\ifx\diffarg\{%
\let\DiffSpace\relax
\fi\fi\fi\DiffSpace}

\providecommand*{\eu}%
{\ensuremath{\mathrm{e}}}

\providecommand*{\iu}%
{\ensuremath{\mathrm{i}}}

\newcommand{\im}{\mathop{\mathrm{Im}}}

\newcommand{\Hom}{\mathop{\mathrm{Hom}}}

\newcommand{\vol}{\mathop{\mathrm{vol}}}

\newcommand{\Norm}{\mathop{\mathrm{N}}}
\newcommand{\tr}{\mathop{\mathrm{tr}}}
\newcommand{\rank}{\mathop{\mathrm{rank}}}
\newcommand{\End}{\mathop{\mathrm{End}}}
\newcommand{\codim}{\mathop{\mathrm{codim}}}

\begin{document}
\title[Restriction of the metaplectic representation]{Restriction of
  the metaplectic representation over a $p$-adic field to an
  anisotropic torus}

    \author{Khemais Maktouf} \address{{University of Sousse,
        Laboratoire Physique-Mathématique, Fonctions Spéciales et
        Applications (LR 11 ES 35), ESSTHS 4002 Sousse, Tunisia}}
    \email{{khemais.maktouf@fsm.rnu.tn}}

\author{Pierre Torasso} \address{{Université de Poitiers, CNRS, UMR
    6086 Laboratoire de Mathémati\-ques et Applications, Site du
    Futuroscope, TSA 61125, 86073 Poitiers Cedex 9, France}}
\email{{pierre.torasso@math.univ-poitiers.fr}}

\subjclass[2020]{22E50}

\keywords{
  Metaplectic group, Weil representation, maximal irreducible
  tori, anisotropic tori, admissibility, restriction, symplectic reduction}

\begin{abstract}
 
  \noindent
In this article, we examine the restriction of the metaplectic
representation $\pi$ over a $p$-adic field $k$, $p\neq2$, of zero
characteristic to an anisotropic torus $S$ contained in the symplectic
group. First we give necessary and sufficient conditions on the
momentum map in order that $S$ be admissible, that is $\pi_{\vert S}$
decomposes with finite multiplicities.

Let us say that a torus contained in the symplectic group is
irreducible if its action on the symplectic space is irreducible over
$k$. Then we examine the case when $S$ is a proper subtorus of a maximal
irreducible torus $T$ in the symplectic group and give sufficient
conditions on $T$ in order that $S$ never be admissible. When these
conditions are not satisfied, we give examples of admissible proper
tori of a maximal irreducible torus.

Finally, for any admissible subtorus $S$ of a certain type of maximal
irreducible torus, we compute the multiplicity of the unitary
characters of $S$ appearing into $\pi_{\vert S}$. We also show that
the multiplicity of such a character is equal to the volume of the
symplectic reduction of the inverse image under the momentum map of a
linear form associated to it.

%Let $W$ be a symplectic space over a $p$-adic field of zero
%characteristic, $Mp(W)$ the metaplectic group which is the non-trivial
%double covering of the symplectic group $Sp(W)$, and $S$ an anistropic
%torus over $k$ contained in $Sp(W)$. One knows that $S$ admits a
%lifting in $Mp(V)$. An important question is to study the restriction
%of the so called metaplectic or Weil representation of $Mp(W)$ to
%$S$. In particular, we say that the metaplectic representation is
%admissible with respect to $S$ or simply that $S$ is admissible, if
%the restriction of the metaplectic representation to $S$ decomposes
%with finite mutiplicities.

%The symplectic action of $Mp(W)$ into $V$ induces naturally a momentum
%map $\phi$ from $W$ to the dual $\mathfrak{s}^{*}$ of the Lie algebra
%$\mathfrak{s}$ of $S$.

%We show that the following assertions are equivalent :

%(i) the torus $S$ is admissible

%(ii) $0\notin\phi(W\smallsetminus\{0\})$

%(iii) $\phi$ is a proper map

%(iv) the restriction of $\phi$ to $W\smallsetminus\{0\}$ is proper on its image.
\end{abstract}

\maketitle

\section*{Introduction}
The so-called metaplectic or Weil representation appears naturally in
the theory of representations of Lie or algebraic groups, both in the
real and the $p$-adic case, as a corner stone for the construction of
irreducible unitary representations.

One important question in this theory is to understand branching
rules, that is, how the restriction of an irreducible representation of
a group $G$ to one of its subgroups $H$ decomposes into
irreducible ones. For instance one would like to characterize situations
in which the restriction to $H$ of a representation of $G$ decomposes
with finite multiplicities, and in such a situation give a recipe to
compute these multiplicities.

%In the real case this problem when $G$ and $H$ are reductive groups
%has been addressed by many people who obtained a lot of interesting
%results. We may cite Duflo, Guillemin,
%Meinrenken, Paradan, Sjamar, Sternberg, Vargas, Vergne...

Due to the central place occupied by the Weil representation in the
construction of representations of algebraic groups, it seems to be of
interest to study its restriction to various subgroups of the
metaplectic group. This article is an attempt to do that in the
$p$-adic case for anisotropic tori contained in the metaplectic group.

Let $k$ be a local field of zero characteristic. % and residual
%characteristic $p\neq2$.
Let $\psi$ be a unitary character of the additive group $k$. Let $W$
be a symplectic space over $k$ with symplectic form $\beta$. Let
$Sp(W)$ be the symplectic group of $W$, that is the group of $g\in
GL(W)$ such that $\beta(gx,gy)=\beta(x,y)$, $x,y\in W$.

Let $S\subset Sp(W)$ be an anisotropic torus over $k$ and
$\mathfrak{s}$ be its Lie algebra. The action of $S$ in $W$ is
symplectic. Let $\phi:W\rightarrow \mathfrak{s}^{*}$ be the associated
momentum map which is defined as follows :
\begin{equation*}
  \langle\phi(w),X\rangle=-\frac{1}{2}\beta(X\cdot w,w)\mbox{, }w\in W\mbox{, }
  X\in\mathfrak{s}.
\end{equation*}

Let $Mp(W)$ be the corresponding metaplectic group, the
non-trivial double covering of $Sp(W)$. We have the exact
sequence of groups :
\begin{equation*}
1\rightarrow\{1,\epsilon\}\rightarrow Mp(W)\rightarrow Sp(W)\rightarrow1.
\end{equation*}

Let $\pi$ be the so-called
metaplectic or Weil representation of $Mp(W)$ associated to the
character $\psi$. See for instance \cite{lion-vergne-1980} for the
real case and \cite[chapitre 2]{moeglin-vigneras-waldspurger-1987} for
the $p$-adic case. Then $\pi$ is a genuine representation of $Mp(W)$
which, although not irreducible, is the hilbertian sum of two
irreducible ones and satisfies $\pi(\epsilon)=-Id$.

Let $\tilde{S}$ be the inverse image of $S$ in $Mp(W)$. It turns
out that $\tilde{S}$ is an abelian compact group. Thus, the restriction
of $\pi$ to $\tilde{S}$ decomposes as a direct sum of unitary
characters. We are interested in studying this decomposition.

Let $\hat{\tilde{S}}$ be the group of unitary characters of
$\tilde{S}$. Then there exists $m(\chi)\in\mathbb{N}\cup\{\infty\}$
such that
\begin{equation*}
\pi_{\vert S}=\oplus_{\chi\in\hat{\tilde{S}}}m(\chi)\chi.
\end{equation*}

We say that the torus $S$ is admissible if $\pi_{\vert\tilde{S}}$
decomposes with finite multiplicities, that is $m(\chi)\in\mathbb{N}$,
$\chi\in\hat{\tilde{S}}$.

The question is twofold: First to give necessary and sufficient conditions on
the momentum map for $S$ to be admissible. Second, assuming these
conditions satisfied, to express the multiplicity of a character
of $\tilde{S}$ in the Weil representation in terms of geometric
objects built upon the momentum map.

Theorem \ref{theo3.1}, the first main result of this paper, answers to
the first question. More precisely, it says that the following
conditions are equivalent

\smallskip
(i) $S$ is admissible,

(ii) $\phi^{-1}(0)=\{0\}$,

(iii) $\phi$ is a proper map,

(iv) the restriction of $\phi$ to $W\smallsetminus\{0\}$ is proper on
its image.
\smallskip

Before going ahead in the description of the results contained in this
paper, let us say few words about what is known in the real case. In
this case, the symplectic space $W$ may be identified with
$\mathbb{C}^{r}$ equipped with the canonical hermitian product $h$,
the symplectic form $\beta$ being the real part of $\iu h$. Then the
unitary group $U(r)$ is a maximal compact subgroup of $Sp(W)$. Remark
that $U(1)$ naturally identifies with the group of complex
numbers with module $1$.  Moreover, the inverse image of $U(r)$ in
$Mp(W)$ may be identified with
\begin{equation*}
\tilde{U}(r)=\{(g,z)\in U(r)\times U(1)\vert z^{2}=\det g\}.
\end{equation*}

Then the subgroup $T$ of diagonal matrices in $U(r)$ is a maximal
anisotropic torus in $Sp(W)$. Also the inverse image $\tilde{T}$ of
$T$ in $Mp(W)$ is a non-trivial double covering of $T$. If $S$ is an
anisotropic torus of $Sp(W)$, up to conjugation we may suppose that
$S$ is contained in $T$. Let $\tilde{S}$ be the inverse image of $S$
in $Mp(W)$. The problem is to study the restriction of the Weil
representations to $\tilde{S}$.

The tensor product of the restriction to $\tilde{U}(r)$ of the Weil
representation by the character $(g,z)\mapsto z$ factors to a
representation of $U(r)$. The space of $U(r)$-finite vectors in this
one identifies itself to the space
$\mathbb{C}[z_{1},\ldots,z_{r}]$ of polynomial functions on
$\mathbb{C}^{r}$ on which $U(r)$ acts by left translation on the
variables. Thus, the problem of studying the restriction of the
metaplectic representation to $\tilde{S}$ reduces to study the
restriction to $S$ of the natural representation of $U(r)$ in
$\mathbb{C}[z_{1},\ldots,z_{r}]$.

The kernel of the exponential map from $\mathfrak{s}$ onto $S$ is a
$\mathbb{Z}$-lattice $L$ of rank $\dim S$. Let $L^{*}$ be the dual
lattice, that is the set of linear forms $f\in\mathfrak{s}^{*}$ such
that $f(L)\subset2\pi\mathbb{Z}$. If $f\in L^{*}$, the formula
$\chi_{f}(\exp X)=\eu^{\iu f(X)}$, $X\in\mathfrak{s}$, defines a
unitary character of $S$ and the map $f\mapsto\chi_{f}$ is a bijection
from $L^{*}$ onto the group $\hat{S}$ of unitary characters of $S$. If
$f\in L^{*}$, we denote by
$m(f)=m(\chi_{f})\in\mathbb{N}\cup\{\infty\}$, the multiplicity of
$\chi_{f}$ in the restriction of the above representation of $U(r)$.

As $S$ is contained in $T$, there exists
$\alpha_{1},\ldots,\alpha_{r}\in L^{*}$ such that one has
\begin{equation*}
\exp
X\cdot z=(\eu^{\iu\alpha_{1}(X)}z_{1},\ldots,\eu^{\iu\alpha_{r}(X)}z_{r})\mbox{,
}X\in\mathfrak{s}\mbox{, }z\in\mathbb{C}^{r}.
\end{equation*}

It follows that the momentum map is given by
\begin{equation}\label{eqI1}
\phi(z)=\frac{1}{2}\sum_{1\leq j\leq r}\vert z_{j}\vert^{2}\alpha_{j}.
\end{equation}

Moreover, it is well known that for $f\in L^{*}$, one has
\begin{equation*}
  m(\chi_{f})=
  \vert\{m\in\mathbb{Z}^{r}\vert f=m_{1}\alpha_{1}+\cdots+m_{r}\alpha_{r}\}\vert
\end{equation*}
where, for any set $A$, $\vert A\vert$ denotes its cardinal.

Then, it follows easily from formula \ref{eqI1} and \cite[first
  page of the introduction]{guillemin-1997} that the fact that $S$ is
admissible is equivalent to conditions (ii), (iii) and (iv) of our
Theorem \ref{theo3.1} and also to
\smallskip

(v) the linear forms $\alpha_{1},\ldots,\alpha_{r}$ all lie in a fixed
open half-space of $\mathfrak{s}^{*}$.

\smallskip
Suppose that these conditions are satisfied and let $f\in L^{*}$. We
say that $f$ is in general position if for every subset $I$ of
$\{1,\ldots,r\}$ such that $f\in\sum_{i\in I}\mathbb{R}\alpha_{i}$,
the system $(\alpha_{i})_{i\in I}$ generates $\mathfrak{s}^{*}$.

When $f$ is in general position, $\phi^{-1}(f)$ is a submanifold of
$W$ and its quotient $X_{f}=S\backslash\phi^{-1}(f)$ is an orbifold
which inherits a symplectic structure from $\beta$ (see \cite[bottom
  of page 54]{guillemin-1997}). The symplectic orbifold $X_{f}$ is a
so-called reduced symplectic space or also the symplectic reduction of
$\phi^{-1}(f)$. Then as explained in \cite[top of page
  55]{guillemin-1997} it follows from \cite{meinrenken-1996} and
\cite{vergne-1994}  that the multiplicity $m(\chi_{f})$
is the Kawasaki-Riemann-Roch number of $X_{f}$, which is obtained by
integration over $X_{f}$ of a particular differential form. 

Recall that a torus contained in a reductive group $G$ is said to be
elliptic if it is anisotropic modulo the center of $G$. In our case,
$T$ is up to conjugation the unique elliptic maximal torus of
$Sp(W)$. %Then the first part of condition (v) above reads as follows :
%$S$ is contained in an elliptic maximal torus.

Now we suppose that $k$ is a non-archimedean field of residual characteristic
$p\neq2$.  In this case, the landscape is richer than in the real
case.

Let us say that a torus $S$ contained in $Sp(W)$ is irreducible if $W$ is
irreducible as a $S$-module. In the real case this may happen only if
$\dim W=2$ and $S=U(1)$. In the $p$-adic case there always exist
irreducible tori contained in $Sp(W)$. See paragraph \ref{2.1} for a
description of maximal irreducible tori in $Sp(W)$.

Let us say some more words about the content of this article.

First, in our case $S$ being a compact subgroup of $Sp(W)$, it admits
a lifting into $Mp(W)$. Thus, we may consider $S$ as a subgroup of
$Mp(W)$. As $\pi(\epsilon)=-Id$, it is clear that our problem is the
same as studying the restriction of $\pi$ to $S$. More precisely, $S$
is admissible if and only if the Weil representation restricts to $S$
with finite multiplicities, a character $\chi$ of $\tilde{S}$ such
that $\chi(\epsilon)=-1$ appears in $\pi_{\vert\tilde{S}}$ if and only
if $\chi_{\vert S}$ appears in $\pi_{\vert S}$, and in any case we
have $m(\chi)=m(\chi_{\vert S})$.

The proof of Theorem \ref{theo3.1} which says that conditions
(i) to (iv) above are equivalent goes as follows. We remark
that if $W$ contains two transverse Lagrangian subspaces which are
$S$-stable, then $S$ is not admissible and none of the conditions (ii)
to (iv) is true. It follows that we may suppose that $W$ decomposes
irreducibly under $S$ as direct sum of mutually orthogonal symplectic
subspaces $W=\oplus_{1\leq i\leq n} W_{i}$ (see Propositions
\ref{pr3.1} and \ref{pr3.2}). Thus, $S$ is embedded in a product of
maximal irreducible tori $T_{i}\subset Sp(W_{i})$, $1\leq i\leq n$. In
fact the tori $T_{1}\times\ldots\times T_{n}$ obtained so are the
elliptic maximal tori in $Sp(W)$. Contrary to the real case, it turns
out that anisotropic tori in $Sp(W)$ are not always contained in an
elliptic maximal torus (see appendix \ref{apB}).

At this stage, we compute the momentum map and obtain formula
\eqref{eq12}, which is an analogue of formula \eqref{eqI1} in the real
case. In paragraph \ref{3.6}, we prove the equivalence of conditions
(ii), (iii) and (iv).

In paragraph \ref{3.7} and \ref{3.8}, we prove the equivalence of
conditions (i) and (ii) : see Proposition \ref{pr3.4}, Corollary
\ref{co3.1} and Theorem \ref{theo3.2}. For doing that, we need to know
the restriction of the Weil representation to maximal irreducible
tori.  Theorem \ref{theo2.1} asserts that the restriction of the Weil
representation to a maximal irreducible torus is multiplicity free and
describes the characters which appear therein. This is an easy
consequence of the work of Yang \cite{yang-1998} which concerns the
case of $SL_{2}(k)$.

%Moreover, our proof of theorem \ref{theo3.2}
%requires to use the exponential map for the maximal irreducible tori
%in $Sp(W)$, mainly to manage with the embedding of the torus $S$ in a
%maximal elliptic torus. For doing that and also in view of
%forthcoming papers on the subject, in paragraph \ref{2.3} we give some
%precise results about the domain of definition of the exponential map
%and the logarithm.

Once known the conditions for an anisotropic torus of $Sp(W)$ to be
admissible, one question arises naturally : given a maximal
irreducible torus $T$ in $Sp(W)$, does there exist a proper subtorus
$S\subset T$ which still is admissible ?

If $T\subset Sp(W)$ is a
maximal irreducible torus, there exists field extensions $k\subset
k'\subset k''$, with
\begin{itemize}
\item  $k''$ of degree $2$ over $k'$
\item  $k'$ of degree
  $\frac{1}{2}\dim W$ over $k$
\end{itemize}
  and a $k$-linear isomorphism from $k''$
over $W$ such that $T$ is the subgroup of $k''^{\times}$ of elements
of norm $1$ over $k'$ (acting by multiplication on $k''$). Then we
show that if $k''$ is unramified over $k'$ and the order of
ramification of $k'$ over $k$ is odd, or if $k''$ is ramified over
$k'$ and the modular degree of $k'$ over $k$ is odd, there is no
proper subtorus of $T$ which is admissible (see Proposition
\ref{pr4.1}). The same is true if $k''$ is ramified over $k'$, $k'$ is
unramified of even degree over $k$ and there is no uniformizer of
$k''$ whose square is a uniformizer of $k$ (see Proposition
\ref{pr4.2}).

Then, in the case where $k''$ is ramified over $k'$, $k'$ is
unramified of even degree over $k$ and $k''$ admits a uniformizer
whose square is a uniformizer of $k$ (case A) we are only able to give
examples of proper subtori of $T$ which are admissible : see
Proposition \ref{pr4.3}. In the case where $k''$ is unramified over
$k'$, $k'$ is totally and tamely ramified of even degree over $k$
(case B) we obtain a nice and complete description of the admissible
subtori of $T$ : see Proposition \ref{pr4.6}.

For doing that, we must understand the structure of the
group $\mathbf{X}(T)$ of rational characters of the torus $T$ as a
module over the $\mathbb{Z}$-algebra of the Galois group of a finite
Galois extension of $k$ over which $T$ splits. As is well known,
subtori of $T$ correspond to submodules $M$ of $\mathbf{X}(T)$ such
that the quotient $\mathbf{X}(T)/M$ is free over $\mathbb{Z}$. Our
study of these submodules relies on Theorem \ref{theoap} of the
appendix describing for a finite cyclic group $\Gamma$ the
$\mathbb{Z}[\Gamma]$-submodules $M$ of $\mathbb{Z}[\Gamma]$ such that
$\mathbb{Z}[\Gamma]/M$ is free over $\mathbb{Z}$ and minimal for this
property. For this theorem we are indebted to François
Courtès$^{\dag}$.

%Precisely, we study more accurately the cases where $k''$ is ramified
%over $k'$, $k'$ is unramified of even degree prime to $p$ over $k$ and
%the square of a uniformizer of $k''$ is equal to a uniformizer of $k$
%(case A) or $k''$ is unramified over $k'$ and $k'$ is totally and
%tamely ramified over $k$ of even degree prime to $p$ (case B).

%Contrary to case A where we are only able to give examples of such
%tori (see proposition \ref{pr4.3}), in case B, we obtain a nice and
%complete description of the admissible subtori $S$ of $T$ (see
%proposition \ref{pr4.6}).

In both cases A and B, with the restriction in case A that the degree
of $k'$ over $k$ is prime to $p$, given an admissible subtorus
$S\subset T$, we compute the multiplicity of the unitary characters of
$S$ in the Weil representation (Propositions \ref{pr5.2} and
\ref{pr5.3}), describe the structure of the symplectic reduction
$S\backslash\phi^{-1}(g)$, $g\in\im\phi\subset\mathfrak{s}^{*}$ (Lemma
\ref{lem5.15.5}, Corollary \ref{co5.3}, and Lemma \ref{lem5.16.5}), which
turns out to be a manifold and compute its volume (Corollary
\ref{co5.4}). Then to each unitary character $\chi$ of $S$ appearing
in the restriction to $S$ of the metaplectic representation
corresponds a family of linear forms on $\mathfrak{s}$ which
determines $\chi$ on a neighbourhood of $1$. If $\chi$ is non-trivial
and $g$ is such a linear form which is contained in the image of
$\phi$, Theorems \ref{theo5.1} and \ref{theo5.2} assert that the
multiplicity of $\chi$ in the Weil representation is equal to the
volume of the symplectic reduction $S\backslash\phi^{-1}(g)$,
generalizing in these particular cases known results in the real case
(see again \cite[Introduction]{guillemin-1997}).

On the way to treat the general case, next step would be to address
the case where $S$ is a subtorus of a maximal irreducible subtorus
with the only restriction that $k'$ is tamely ramified and that its
modular degree is prime to $p$. An other interesting case in this
direction, would be the one where $S$ is contained in a maximal
elliptic torus which is the product of copies of a one dimensional
anisotropic torus, a case analogous to the general real case.

We would like to thank Paul Broussous for interesting discussions on
the subject. We also thank the referee for his helpfull suggestions.

\section{Generalities}\label{1}

\subsection{}\label{1.1}
Along this paper, $k$ denotes a non-archimedean local field of zero
characteristic and residual characteristic $p\neq2$. Let $\mathcal{O}$
be the ring of integers of $k$, $\mathcal{P}$ its maximal ideal, and
$\varpi$ a uniformizer or a prime element of $k$, that is
$\mathcal{P}=\varpi\mathcal{O}$. Then, the residual class field
$\mathcal{O}/\mathcal{P}=\mathbb{F}_{q}$ is a finite field of
characteristic $p$ and cardinal $q$, a power of $p$. Moreover, $k$ is
a finite extension of $\mathbb{Q}_{p}$. Let $v :
k^{\times}\rightarrow\mathbb{Z}$ be the discrete valuation such that
$v(\varpi)=1$ and $\vert\,\vert$ be the absolute value normalized by :
$\vert x\vert=q^{-v(x)}$, $x\in k^{\times}$. We extend $v$ to $0$ by
putting $v(0)=+\infty$. This allows us to write $\vert0\vert=q^{-v(0)}$.

We fix a non-trivial unitary character $\psi$ of the additive group $k$, and
denote by $\lambda_{\psi}$ its conductor, that is $\lambda_{\psi}$ is
the smallest integer such that $\mathcal{P}^{\lambda_{\psi}}$ is contained in
the kernel of $\psi$.

\subsection{}\label{1.2}
Let $(W,\beta)$ be a symplectic space over $k$ of dimension $2r$. We
denote by $Sp(W)$ the related symplectic group and $Mp(W)$ the related
metaplectic group, that is the non-trivial double covering of
$Sp(W)$. Thus, we have the exact sequence :
\begin{equation}\label{eq1}
1\rightarrow\{1,\epsilon\}\rightarrow Mp(W)
\rightarrow Sp(W)\rightarrow1
\end{equation}
where $\{1,\epsilon\}$ is the multiplicative group with two elements.

According to \cite[Section 4]{maktouf-torasso-2016}, the exact
sequence \eqref{eq1} splits canonically (and even uniquely if
$q\geq4$) above any maximal compact subgroup of $Sp(W)$. Thus, every
compact subgroup of $Sp(W)$ admits a lifting in $Mp(W)$, which in
general is neither canonical, nor unique.

\subsection{}\label{1.3}
Let $\pi$ be the so-called metaplectic or Weil representation\! of \!$Mp(W)$
associated to the character $\psi$. See for instance \cite[chapitre
  2]{moeglin-vigneras-waldspurger-1987}. Then $\pi$ is a genuine
representation of $Mp(W)$.

Let $K\subset Sp(W)$ be a compact subgroup. Choosing a lifting of $K$
into $Mp(W)$, we may consider it as a subgroup of $Mp(W)$ and then
study the restriction of $\pi$ to $K$.

\subsection{}\label{1.4}
Let $\ell$ and $\ell'$ be two Lagrangian subspaces of $W$ such that
$W=\ell\oplus\ell'$. The symplectic form $\beta$ induces a duality
between $\ell$ and $\ell'$. We denote by $^{t}a$ the transpose of the
homomorphism $a$ with respect to this duality. Then, $Sp(W)$ is the set
of endomorphisms of $W$ which, in the decomposition $W=\ell\oplus\ell'$,
are of the form
$$
\begin{pmatrix}
  a & b\\
  c & d
\end{pmatrix},
$$ with $a\in L(\ell)$, $d\in L(\ell')$,
$b\in L(\ell',\ell)$ and $c\in L(\ell,\ell')$ such that
$a\,^{t}d-b\,^{t}c=Id_{\ell}$, $a\,^{t}b=b\,^{t}a$, and $c\,^{t}d=d\,^{t}c$.

Let $Sp(W)_{\ell,\ell'}$ be the stabilizer in $Sp(W)$ of the
subspaces $\ell$ and $\ell'$. %Its elements are of the form
%$$
%\begin{pmatrix}
%  a & 0\\
%  0 & ^{t}a^{-1}
%\end{pmatrix}
%$$ with $a\in GL(\ell)$.
The subgroup $Sp(W)_{\ell,\ell'}$ admits a lifting into
$Mp(W)$. Moreover, {\color{black}once chosen a Haar measure on
  $\ell$,} the representation $\pi$ admits a realization in
  $\mathrm{L}^{2}(\ell)$ {\color{black}in which its restriction to
    $Sp(W)_{\ell,\ell'}$ is given by}
\begin{equation}\label{eq2}
\pi(g)\phi(x)=\chi_{\ell,\ell'}(a)
\vert\sideset{}{_{\ell}}\det a\vert^{-\frac{1}{2}}
\phi(a^{-1}x), \phi\in\mathrm{L}^{2}(\ell), x\in\ell
\end{equation}
where $\chi_{\ell,\ell'}$ is a unitary character of $GL(\ell)$, and
\begin{equation*}
g=\begin{pmatrix}  a & 0\\
 0 & ^{t}a^{-1}
\end{pmatrix}\in Sp(W)_{\ell,\ell'}\mbox{, }a\in GL(\ell).
\end{equation*}

\section{Restriction of the Weil representation to maximal irreducible
  tori}\label{2}

\subsection{}\label{2.1}
%Let $\mathbf{Sp}(W)$ be the algebraic group defined over $k$ such that
%$Sp(W)$ is the group of $k$-points. A torus 
In the following, a torus is the group of $k$-points of an algebraic
torus defined over $k$.

\begin{defi}\label{defi2.1}
A torus $T\subset Sp(W)$ is irreducible if $W$ is irreducible over $k$ as
a $T$-module.
\end{defi}

We put $2r=\dim W$. Let $T\subset Sp(W)$ be a maximal irreducible
torus. According to \cite{howe-1973} there exist an extension $k'$ of
$k$ of degree $r$, a quadratic extension $k''=k'(u)$ of $k'$ with
$\tr_ {k''/k'}u=0$, and a linear isomorphism of $k''$ onto $W$ such
that, identifying $W$ with $k''$ through this isomorphism and denoting
by $\tau$ the non-trivial element of the Galois group of $k''$ over
$k'$, we have :

(i) $\beta(x,y)=\frac{1}{2}\tr_{k''/k}uxy^{\tau}\mbox{, }x,y\in W$,

(ii) $T=\{\xi\in k''\vert\Norm_{k''/k'}(\xi)=1\}$ acting in $W$ by
multiplication.

We denote by $\mathcal{O}'$ (resp. $\mathcal{O}''$) the ring of
integers of $k'$ (resp. $k''$), $\mathcal{P}'$ (resp. $\mathcal{P}''$)
its maximal ideal, $\varpi'$ (resp. $\varpi''$) one of its
uniformizers, $v'$ (rep. $v''$) the associated valuation, and $q'$
(resp. $q''$) the cardinal of its residual class field.

Let $e$ be the order of ramification of $k'$ over $k$ and $\delta$ the
differental exponent of $k'$ over $k$, that is $\delta$ is the largest
integer such that $\tr_{k'/k}x\in\mathcal{O}$ for all
$x\in\mathcal{P}'^{-\delta}$.

We put $\mu=e\lambda_{\psi}-\delta-v''(u)$.

\subsection{}\label{2.2}
We suppose first that the extension $k''/k'$ is unramified so that we may
take $\varpi''=\varpi'$. In this case, according to \cite[Theorem 1.6
  and Proposition 1.10]{morris-1991}, up to conjugation of the maximal
tori in $Sp(W)$, we may suppose that $v''(u)\in\{0,1\}$. Let us write
$u=\varpi'^{v''(u)}\nu$, with $\nu\in\mathcal{O}''^{\times}$ and
$\tr_{k''/k'}\nu=0$. We put $d=\nu^{2}$ which is an element of
$\mathcal{O}'^{\times}$ and not a square.

Then we have
$$ T=\{\xi+\eta\nu\vert\xi,\eta\in\mathcal{O}'\mbox{ and }
\xi^{2}-d\eta^{2}=1\},
$$
and the two choices $v''(u)=0$ or $1$ correspond to the two
distinct conjugacy classes in $Sp(W)$ of tori isomorphic to $T$.

Now, we suppose that $k''/k'$ is ramified. In this case we may choose
$\varpi'$ and $\varpi''$ such that $\varpi''^{2}=\varpi'$, and we put
$u=\nu=\varpi''$. Then we have
$$
T=\{\xi+\eta\varpi''\vert\xi,\eta\in\mathcal{O}'\mbox{ and }
\xi^{2}-\varpi'\eta^{2}=1\}.
$$

In both situations, for $j\in\mathbb{N}$ we denote by $T_{j}$ the
$j$-th congruence subgroup of $T$, that is
$$
T_{j}=\{g\in T\vert g-1\in\varpi''^{j}\mathcal{O}''\}.
$$
Note that $T_{1}\subsetneq T_{0}=T$.

If $l\geq1$ we denote by $\mu_{l}$ the cyclic group of $l$-th roots of
unity. For the following lemma, see \cite{yang-1998}.

\begin{lem}\label{lem2.1}
  a) Suppose that $k''$ is unramified over $k'$. Then one has
  $T=\mu_{q'+1}\times T_{1}$. Moreover, $(T_{j})_{j\in\mathbb{N}}$ is a
  strictly decreasing sequence of open subgroups of $T$ and a basis of
  neighbourhoods of $1$.

  b) Suppose that $k''$ is ramified over $k'$. Then one has
  $T=\{1,-1\}\times T_{1}$ and $T_{2j}=T_{2j+1}$, $j>0$. Moreover,
  $(T_{2j+1})_{j\in\mathbb{N}}$ is a strictly decreasing sequence of
  open subgroups of $T$ and a basis of neighbourhoods of $1$.
\end{lem}

Let $l\geq0$ be an integer. The conductor of a unitary character
$\chi$ of $T_{l}$ is the smallest integer $j\geq l$ such that $\chi$
is trivial on $T_{j}$. We shall denote it by $c(\chi)$.

The trivial character of $T$ is the only one whose conductor is $0$.

The characters of $T$ whose conductor is less or equal to $1$ are the
characters of $T/T_{1}$. When $k''$ is unramified over $k'$ they form
a cyclic group of order $q'+1$. We denote by $\eta_{0}$ the unique
element among them which is of order $2$.

When $k''$ is ramified over $k'$, there is only one character with the
conductor $1$ : it is the non-trivial character of the two element
group $T/T_{1}$. It is of order $2$, and we also denote it $\eta_{0}$. All
other characters of $T$ have an even conductor.

The following proposition is an easy reformulation and generalization
of \cite[Proposition 1.3]{yang-1998}.

\begin{pr}\label{pr2.1}
(i) Suppose that $k''$ is unramified over $k'$. Let $j, l$ be integers
  such that $0<\frac{j}{2}\leq l<j$ and $\chi$ be a character of
  $T_{l}$. The following assertions are equivalent :

a) the conductor of $\chi$ is $j$,

b) there exists $b\in\mathcal{O}'^{\times}$ such that
\begin{equation}\label{eq3}
\chi(x)=\psi(-\frac{1}{4}\sideset{}{_{k''/k}}\tr\varpi'^{\mu-j}bux)\mbox{,
}x\in T_{l},
\end{equation}

c)  there exists $a\in\mathcal{O}''^{\times}$ such that
\begin{equation}\label{eq3b}
  \chi(x)=\psi(-\frac{1}{4}\sideset{}{_{k''/k}}\tr\varpi'^{\mu-j}
  \sideset{}{_{k''/k'}}\Norm(a)
  ux)\mbox{,
}x\in T_{l}.
\end{equation}

(ii) Suppose that $k''$ is ramified over $k'$. Let $j, l$ be
integers such that $0<j\leq l<2j$, and $\chi$ be a character
of $T_{l}$.  The following assertions are equivalent :

a) the conductor of $\chi$ is $2j$,

b) there exists $b\in\mathcal{O}'^{\times}$ such that
\begin{equation}\label{eq4}
  \chi(x)=\psi(-\frac{(-1)^{\mu-j}}{4}\sideset{}{_{k''/k}}
    \tr\varpi'^{\mu-j}b\varpi''x)\mbox{, }x\in T_{l}.
\end{equation}
\end{pr}

In the following we shall denote by $\chi_{b,l,j}$
(resp. $\chi_{b,l,2j}$) the character of $T_{l}$ defined by formula
\eqref{eq3} (resp. \eqref{eq4}).

\subsection{}\label{2.3}
For the following lemma see \cite[p. 154]{gouvea-1993}

\begin{lem}\label{lem2.2}
  (i) Let $x\in k''^{\times}$ such that $v''(x-1)>0$. Then the series
$$
\ln x=\sum_{l\geq1}(-1)^{l-1}\frac{(x-1)^{l}}{l}
$$ converges, and moreover $v''(\ln x)=v''(x-1)$ if
$v''(x-1)>\frac{v''(p)}{p-1}$.

(ii)  Let $X\in k''$ such that $v''(X)>\frac{v''(p)}{p-1}$. Then the series
$$
\exp X=\sum_{l\geq0}\frac{X^{l}}{l!}
$$ converges, and moreover $v''(\exp X-1)=v''(X)$.

(iii) For $j>\frac{v''(p)}{p-1}$, $\exp$ induces an
isomorphism from the additive group $\varpi''^{j}\mathcal{O}''$ onto
the multiplicative group $1+\varpi''^{j}\mathcal{O}''$, the inverse
isomorphism being the restriction of $\ln$ to this last group.
\end{lem}

Remark that the Lie algebra of $T$ is $\mathfrak{t}=k'\nu$. Then the
following corollary follows readily from Lemma \ref{lem2.2}.

\begin{co}\label{co2.1}
(i) Suppose $k''$ is unramified over $k'$. Then for
  $j>\frac{v''(p)}{p-1}$, $\exp$ induces an isomorphism from
  the additive group $\varpi''^{j}\mathcal{O}'\nu$ onto the
  multiplicative group $T_{j}$, the inverse isomorphism being the
  restriction of $\ln$ to this last group.

(ii) Suppose $k''$ is ramified over $k'$. Then for $j$ such that
  $2j+1>\frac{v''(p)}{p-1}$, $\exp$ induces an
  isomorphism from the additive group $\varpi''^{2j+1}\mathcal{O}'$
  onto the multiplicative group $T_{2j}=T_{2j+1}$, the inverse
  isomorphism being the restriction of $\ln$ to this last group.
\end{co}

\subsection{}\label{2.4}
\begin{pr}\label{pr2.2}   
  (i) Suppose $k''$ is unramified over $k$. Let $l_{k''}$ be the
  smallest integer such that $l_{k''}>\frac{v''(p)}{p-1}$.
  
  Let $\chi$ be a unitary character of $T$ and $j\geq l_{k''}$. Then the
  following assertions are equivalent

  a) the conductor of $\chi$ is $j$,

  b) there exists $b\in\mathcal{O}'^{\times}$ such that
\begin{equation}\label{eq5}
\chi(\exp
X)=\psi(-\frac{1}{4}\sideset{}{_{k''/k}}\tr\varpi'^{\mu-j}buX)\mbox{,
}X\in\varpi'^{l_{k''}}\mathcal{O}'\nu,
\end{equation}

c) there exists $a\in\mathcal{O}''^{\times}$ such that
\begin{equation}\label{eq6}
\chi(\exp
X)=\psi(-\frac{1}{4}\sideset{}{_{k''/k}}\tr\varpi'^{\mu-j}
\sideset{}{_{k''/k}}\Norm(a)uX)\mbox{,
}X\in\varpi'^{l_{k''}}\mathcal{O}'\nu.
\end{equation}

(ii) Suppose $k''$ is ramified over $k$. Let $l_{k''}$ be the
  smallest integer such that $2l_{k''}+1>\frac{v''(p)}{p-1}$.

Let $\chi$ be a unitary
  character of $T$ and $j\geq l_{k''}$. Then the following assertions
  are equivalent

  a) the conductor of $\chi$ is $2j$,

  b) there exists $b\in\mathcal{O}'^{\times}$ such that
\begin{equation}\label{eq7}
\chi(\exp
X)=\psi(-\frac{(-1)^{\mu-j}}{4}\sideset{}{_{k''/k}}
\tr\varpi'^{\mu-j}b\varpi''X)\mbox{,
}X\in\varpi'^{l_{k''}}\mathcal{O}'\nu.
\end{equation}
\end{pr}
\begin{dem}
The equivalence of assertions a) and b) in both cases (i) and (ii)
readily follows from Corollary \ref{co2.1}

The equivalence of assertions b) and c) in case (i) follows from the
fact that $\Norm_{k''/k'}(\mathcal{O}''^{\times})=\mathcal{O}'^{\times}$.
\end{dem}

In the following, $[x]$ denotes the integral part of the real number
$x$.

\begin{pr}\label{pr2.3}
(i) Suppose $k''$ is unramified over $k'$ and let $l,j$ be integers such
  that $l\geq\frac{2v''(p)}{p-1}$, $j\in\{2l-1,2l\}$, and
  $b\in\mathcal{O}'^{\times}$. Then we have
\begin{equation}\label{eq8}
  \chi_{b,l,j}(\exp X)=\psi(-\frac{1}{4}\sideset{}{_{k''/k}}\tr\varpi'^{\mu-j}buX)
  \mbox{, }X\in\varpi'^{l}\mathcal{O}'\nu.
\end{equation}
(ii) Suppose $k''$ is ramified over $k'$ and let $l$ be an integer such
  that $l\geq\frac{2v''(p)}{p-1}$ and
  $b\in\mathcal{O}'^{\times}$. Then we have
\begin{equation}\label{eq9}
  \chi_{b,l,2l}(\exp X)=\psi(-\frac{(-1)^{\mu-l}}{4}\sideset{}{_{k''/k}}
  \tr\varpi'^{\mu-l}b\varpi''X)
  \mbox{, }X\in\varpi'^{[\frac{l}{2}]}%\footnote{}
\mathcal{O}'\varpi''.
  \end{equation}
\end{pr}
\begin{dem}
(i) Let $l>\frac{v''(p)}{(p-1)}$, $j\in\{2l-1,2l\}$ and
  $X\in\varpi'^{l}\mathcal{O}'\nu$. 

 Then the series $\exp X$ converges,
  and one has
\begin{eqnarray*}
  \chi_{b,l,j}(\exp X)&=&\psi\big(-\frac{(-1)^{\mu-j}}{4}\sideset{}{_{k''/k}}
\tr\varpi'^{\mu-j}bu\exp X\big)\\
  &=&\psi\big(-\frac{(-1)^{\mu-j}}{4}\sideset{}{_{k''/k}}
\tr\varpi'^{\mu-j}bu(\varpi''X\!+\!\sum_{m\geq1}\frac{X^{2m+1}}{(2m+1)!})\big)
\end{eqnarray*}
due to the fact that $\sum_{m\geq0}\frac{X^{2m}}{(2m)!}\in
  k'$. Then equality \eqref{eq8} is equivalent to
  \begin{equation*}
\psi\big(-\frac{(-1)^{\mu-j}}{4}\sideset{}{_{k''/k}}
\tr\varpi'^{\mu-j}bu(\sum_{m\geq1}\frac{X^{2m+1}}{(2m+1)!})\big)=1
  \end{equation*}
which is true as soon as $v'(\frac{X^{2m+1}}{(2m+1)!})\geq 2l$,
$m\geq1$.

If $m\geq1$ is an integer,  we define the number $s(m)$ as follows : 
 write     $m=a_{0}+a_{1}p+\cdots+a_{t}p^{t}$ with $a_{i}$ integers such
      that $0\leq a_{i}\leq p-1$, $0\leq i\leq t$ and put
      $s(m)=a_{0}+a_{1}+\cdots+a_{t}$. Then it is well known and easy
      to check that
      $v''(m!)=\sum_{i\geq1}[\frac{m}{p^{i}}]=\frac{(m-s(m))v''(p)}{p-1}$. Remark
      that $s(m)\geq1$.

Then, for $m\geq 1$ one has $v'(\frac{X^{2m+1}}{(2m+1)!})\geq
v'(\frac{\varpi'^{(2m+1)l}}{(2m+1)!})$ with equality if
$X\in\varpi'^{l}\mathcal{O}'^{\times}\nu$. But
$v'(\frac{\varpi'^{(2m+1)l}}{(2m+1)!})=
(2m+1)l-\frac{(2m+1-s(2m+1))v'(p)}{p-1}=
2m(l-\frac{v'(p)}{p-1})+\frac{s(2m+1)-1}{p-1}v'(p)+l\geq
2(l-\frac{v'(p)}{p-1})+l$. It is clear that this last expression is
greater or equal to $l$ as soon as
$l\geq\frac{2v'(p)}{p-1}=\frac{2v''(p)}{p-1}$.

(ii) In this case, for $m\geq1$ we have
    $v'(\frac{X^{2m+1}}{(2m+1)!}\varpi'')\geq
    v'(\frac{\varpi'^{(2m+1)[\frac{l}{2}]}\varpi'^{m+1}}{(2m+1)!})$
    with equality if
    $X\in\varpi'^{[\frac{l}{2}]}\mathcal{O}'^{\times}\varpi''$. But
    $v'(\frac{\varpi'^{(2m+1)[\frac{l}{2}]}\varpi'^{m+1}}{(2m+1)!})
    =(2m+1)[\frac{l}{2}]-v'((2m+1)!)+m+1
    =(2m+1)[\frac{l}{2}]-\frac{(2m+1-s(2m+1))v'(p)}{p-1}+m+1
    \geq2m([\frac{l}{2}]-\frac{v'(p)}{p-1})+[\frac{l}{2}]+m+1
    =m(2[\frac{l}{2}]+1-\frac{v''(p)}{p-1})+[\frac{l}{2}]+1
    \geq(l-\frac{v''(p)}{p-1})+[\frac{l}{2}]+1$.
 Again, equality \eqref{eq9}
    will be true as soon as this last expression will be greater or
    equal to $l$, condition which is satisfied as soon as
    $l\geq\frac{2v''(p)}{p-1}$.
\end{dem}
%\begin{co}\label{co2.2}
%We suppose $k''$ is ramified over $k'$. Let $\chi$ be a character of
%$T$ and $j\geq l_{k''=$ be an integer.  The the assertions are equivalent

 %(i) $\chi$ appears 
%\end{co}

\subsection{}\label{2.5}
In this paragraph we examine the restriction of the Weil
representation of $Mp(W)$ to a maximal irreducible torus $T$.

To embed $T$ into $Mp(W)$ we proceed as follows. First we consider $T$
as a subgroup of $SL_{2}(k')$. It turns out that $T$ is contained in a maximal
compact subgroup of $SL_{2}(k')$ and as such possesses a natural
lifting into the metaplectic group $ML_{2}(k')$ over $SL_{2}(k')$
which is described in \cite{yang-1998}.

Now, we embed $SL_{2}(k')$ into $Sp(W)=Sp(k'')$ as follows. 
We have $k''=k'\oplus k'\nu$, and if
$g=\begin{pmatrix} a&b\\c&d\end{pmatrix}\in SL_{2}(k')$, we consider it as the linear transformation such that for $x,y\in
k'$, $g.(x+y\nu)=(ax+by)+(cx+dy)\nu$. Then one checks easily that this
embedding lifts to an embedding of $ML_{2}(k')$ into $Mp(W)$ such that
the restriction of the Weil representation of $Mp(W)$ to $ML_{2}(k')$
is its own Weil representation associated to the character
$\psi\circ\tr_{k'/k}$ of $k'$.

Let us denote by $\sigma$ the embedding of $T$ into $Mp(W)$ obtained
by composing the embedding into $ML_{2}(k')$ followed by the embedding
of $ML_{2}(k')$ into $Mp(W)$ described above.

Then to embed $T$ into $Mp(W)$ we use $\sigma$ when $k''$ is ramified
over $k'$ or $\mu$ is even and $\eta_{0}\sigma$ if $k''$ is unramified
over $k'$ and $\mu$ is odd (recall that $\eta_{0}$ is
the unique character of $T/T_{1}$ of order $2$. Thus, $\eta_{0}$ takes
its values in $\{\pm1\}$, the center of $Mp(W)$).

In \cite{yang-1998}, Yang gave the description of the restriction of
the Weil representation of $ML_{2}(k')$ to $T$. The following theorem
is a direct consequence of Yang's results and the above considerations.

To be more precise, with the notations in \cite[p. 2394]{yang-1998},
we take $E=k''$, $F=k'$, $\psi\circ\tr_{k'/k}$ as the character of
$F=k'$, and $\alpha=\varpi'^{v''(u)}$, $\delta=\nu$ in case $k''/k'$
is unramified, or $\alpha=1$, $\delta=\varpi''$ otherwise. Then, one
checks that the character $\psi'$ of $E$ is given by
$\psi'(x)=\psi(\frac{1}{4}\tr_{k''/k}ux)$, $x\in E=k''$ and that its
conductor is $\mu$ in case $k''/k'$ is unramified or $2\mu$
otherwise. Moreover, it is clear that, with these choices, the spaces
of the metaplectic representation as described in
\cite[p. 2394]{yang-1998} and in \cite[chapitre 2,
  p. 42]{moeglin-vigneras-waldspurger-1987} are the same.

\begin{theo}\label{theo2.1}
(i) The restriction of the Weil representation to $T$ is a direct sum
  of unitary characters, each one appearing with multiplicity $1$.

(ii) Suppose $k''$ is unramified over $k'$. Then the characters
  appearing in the restriction of the Weil representation to $T$ are
  the ones with conductor having the same parity as $\mu$.

(iii) Suppose $k''$ is ramified over $k'$. Then the characters
  appearing in the restriction of the Weil representation to $T$ are
  the ones with even conductor, and if this conductor is $2j>0$, such
  that their restriction to $T_{j}$ is of the form
  $\chi_{\Norm_{k''/k'}(a),j,2j}$, with $a\in\mathcal{O}''^{\times}$.
\end{theo}

{\color{black}
\begin{pr}\label{pr2.4}
We suppose that $k''$ is ramified over $k'$. Let $\chi$ be a character
of $T$ appearing in the restriction of the Weil representation to $T$
and $j\geq \frac{2v''(p)}{p-1}$ be an integer. Then the following
assertions are equivalent

a) the conductor of $\chi$ is $2j$,

b) there exists $a\in\mathcal{O}''^{\times}$ such that
\begin{equation*}
\chi(\exp
X)=\psi(-\frac{(-1)^{\mu-j}}{4}\sideset{}{_{k''/k}}
\tr\varpi'^{\mu-j}\sideset{}{_{k''/k'}}\Norm(a)\varpi''X)\mbox{,
}X\in\varpi'^{l_{k''}}\mathcal{O}'\nu.
\end{equation*}
\end{pr}
\begin{dem}
According to Proposition \ref{pr2.2} and Theorem \ref{theo2.1}, on one
hand there exists $b\in\mathcal{O}'$ satisfying formula \eqref{eq7}
and on the other hand, there exists $a'\in\mathcal{O}''^{\times}$ such
that $\chi_{\vert T_{j}}=\chi_{\Norm_{k''/k'}(a'),j,2j}$. It follows
from this and Proposition \ref{pr2.3} that
$b-\Norm_{k''/k'}(a')\in\varpi'^{j-[\frac{j}{2}]}\mathcal{O}'$ or
equivalently that
$b\in\Norm_{k''/k'}(a')(1+\varpi'^{j-[\frac{j}{2}]}\mathcal{O}')$. But
according to \cite[Proposition 8, p. 32]{weil-1967}, one has
$(1+\varpi'^{j-[\frac{j}{2}]}\mathcal{O}')^{2}=1+\varpi'^{j-[\frac{j}{2}]}\mathcal{O}'$. Then
we have
$b\in\Norm_{k''/k'}(a')(1+\varpi'^{j-[\frac{j}{2}]}\mathcal{O}')^{2}\subset
\im\Norm_{k''/k'}(\mathcal{O}''^{\times})$. Thus, the proposition.
\end{dem}}

In this section we gave the description of the maximal irreducible
tori in $Sp(W)$, described their structure, recalled the properties
for them of the exponential map. In Theorem \ref{theo2.1} we stated
that the multiplicity of their unitary characters appearing in the
restriction of the Weil representation is always $1$ and gave the list
of these characters. All these results will be used in Section 3 to
prove Theorem \ref{theo3.1} and in the other sections to compute the
multiplicities of the unitary characters of an admissible subtorus of
a maximal irreducible one appearing in the restriction of the Weil
representation.

\section{Momentum map and admissibility of anisotropic tori}\label{3}
\subsection{}\label{3.1}
In this section we consider an anisotropic torus $S\subset
Sp(W)$. Then $S$ is a compact subgroup of $Sp(W)$, and, as such,
possesses a lifting $S\subset Mp(W)$.

\begin{defi}\label{def3.1}
We say that $S$ is admissible (with respect to $\pi$ if we need to be
precise) if the restriction of the metaplectic representation to $S$
decomposes with finite multiplicities.
\end{defi}

It follows easily from the fact that two liftings of $S$ differ from
each other by multiplication by a character of $S$ taking its values into
$\{\pm1\}$ that the property to be admissible for $S$ doesn't depend
upon the chosen lifting.

We will relate the admissibility of $S$ to geometric properties
of the momentum map induced by the symplectic action of $S$ in $W$.

\subsection{}\label{3.2}
\begin{defi}\label{def3.2}
Let $S\subset Sp(W)$ be an anisotropic torus with Lie algebra
$\mathfrak{s}$. The momentum map of the symplectic action of $S$ in
$W$ is the application $\phi:W\rightarrow\mathfrak{s}^{*}$ such that
\begin{equation}\label{eq10}
\langle\phi(w),X\rangle=-\frac{1}{2}\beta(X\cdot w,w)\mbox{, }w\in W\mbox{,
}X\in\mathfrak{s}.
\end{equation}
\end{defi}

In this section we will prove the following theorem :

\begin{theo}\label{theo3.1}
Let $S\subset Sp(W)$ be an anisotropic torus with momentum map
$\phi$. Then the following assertions are equivalent :

(i) $S$ is admissible,

(ii) $0\notin\phi(W\smallsetminus\{0\})$,

(iii) $\phi$ is a proper map,

(iv) the restriction of $\phi$ to $W\smallsetminus\{0\}$ is proper on its
image.
\end{theo}
  
\subsection{}\label{3.3}

\begin{lem}\label{lem3.1}
(i) Let $V\subset W$ be a non-zero $S$-stable subspace which is
  irreducible under the action of $S$. Then either $V$ is an isotropic
  subspace or $V$ is a symplectic subspace.

(ii) Let $V_{1}$ and $V_{2}$ be $S$-stable isotropic subspaces of $W$,
  both irreducible under the action of $S$ and such that
  $\beta(V_{1},V_{2})\neq0$. Then $V_{1}\oplus V_{2}$ is a symplectic
  subspace of $W$.
\end{lem}
\begin{dem}
(i) The subspace $V^{\perp_{\beta}}=\{w\in W\vert\beta(w,v)=0\mbox{,
  }v\in V\}$ is clearly $S$-stable. As $V$ is irreducible under the
  action of $S$, one has either $V\cap V^{\perp_{\beta}}=V$, or $V\cap
  V^{\perp_{\beta}}=0$.

(ii) For $i\neq j\in\{1,2\}$,  $V_{i}\cap
  V_{j}^{\perp_{\beta}}$ is $S$-stable and thus reduced to zero as
  $V_{i}$ is $S$-irreducible and $\beta(V_{1},V_{2})\neq0$. Let
  $v_{i}\in V_{i}$, $i=1,2$, be such that
  $\beta(v_{1}+v_{2},V_{1}+V_{2})=0$. Then $\beta(v_{1},V_{2})=0$ and
  $\beta(v_{2},V_{1})=0$. Hence $v_{1}=v_{2}=0$ and $V_{1}\oplus
  V_{2}$ is a symplectic subspace of $W$.
\end{dem}

\begin{pr}\label{pr3.1}
We may decompose $W$ as a direct sum
\begin{equation}\label{eq11}
W=W_{1}\oplus\cdots\oplus W_{n}
\end{equation}
 of $S$-stable symplectic and mutually orthogonal
subspaces such that for $1\leq i\leq n$

(i) either $W_{i}$ is $S$-irreducible,

(ii) or $W_{i}$ is the direct sum of two $S$-stable and
$S$-irreducible Lagrangian subspaces.
\end{pr}
\begin{dem}
The proof goes by induction on $\dim W$. Consider the set of non-zero
$S$-stable subspaces of $W$ and minimal with respect to the inclusion for
this property. Then according to Lemma \ref{lem3.1}, either we have that one
of this subspaces, say $W_{1}$, is symplectic, or else they all are
isotropic.

In the first case we conclude by applying the induction hypothesis to
$W_{1}^{\perp_{\beta}}$.

In the second case, as $S$ is the group of $k$-points of a reductive
algebraic group defined over $k$, $W$ is a direct sum of $S$-irreducible
isotropic subspaces. Let $V_{1}$ be one of them. Then there exists
another one, say $V_{2}$, such that $\beta(V_{1},V_{2})\neq0$. Thus,
according to Lemma \ref{lem3.1} $W_{1}=V_{1}\oplus V_{2}$ is a $S$-stable
symplectic subspace. We conclude by induction as in the previous
case.
\end{dem}

\subsection{}\label{3.4}
\begin{pr}\label{pr3.2}
  Assume that $W$ contains two transverse and $S$-stable Lagrangian
  subspaces. Then
  
  (i) the torus $S$ is not admissible,

  (ii) none of the assertions (ii) to (iv) of Theorem \ref{theo3.1} is true.
\end{pr}
\begin{dem}
(i) Write $W=\ell\oplus\ell'$, where $\ell$ and $\ell'$ are transverse
  and $S$-stable Lagrangian subspaces.Then $S$ is a subgroup of
  $Sp(W)_{\ell,\ell'}$ and as such identifies itself with a compact
  subgroup of $GL(\ell)$. We realize the Weil representation $\pi$
  into $\mathrm{L}^{2}(\ell)$ as in paragraph \ref{1.4}. Then the
  restriction of $\pi$ to $S$ is given by \eqref{eq2}. Now, as the
  torus $S$ is a compact subgroup of $GL(\ell)$, there exists a basis
  $e_{1},\ldots,e_{r}$ of $\ell$ such that the lattice
  $\mathfrak{a}=\mathcal{O}e_{1}\oplus\cdots\oplus\mathcal{O}e_{r}$ is
  $S$-stable. Moreover, we have $\vert\det_{\ell}g\vert=1$, $g\in
  S$. It follows that the characteristic functions
  $1_{\varpi^{m}\mathfrak{a}}$ of $\varpi^{m}\mathfrak{a}$,
  $m\in\mathbb{Z}$, are eigenvectors of $S$ for the same character
  $\chi_{\ell,\ell'\vert S}$. But the subspace of
  $\mathrm{L}^{2}(\ell)$ generated by these vectors is infinite
  dimensional. Thus, $S$ is not admissible.

(ii) One has $\ell\subset\phi^{-1}(0)$ from what follows that none of
the assertions (ii) to (iv) of Theorem \ref{theo3.1} is true.
\end{dem}

\subsection{}\label{3.5}
In this paragraph, as in Proposition \ref{pr3.1} we consider a
decomposition
$$ W=W_{1}\oplus\cdots\oplus W_{n}.
$$ For $1\leq i\leq n$ let $\pi_{i}$ be
the Weil representation of $Mp(W_{i})$ and
$\gamma_{i}:Sp(W_{i})\subset Sp(W)$ the injection naturally induced by
the above decomposition : if $g\in Sp(W_{i})$, $\gamma_{i}(g)$ acts on
$W$ as $g$ on $W_{i}$ and trivially on $W_{i}^{\perp_{\beta}}$. Then
$\gamma_{i}$ uniquely lifts as an injection
$\tilde{\gamma}_{i}:Mp(W_{i})\subset Mp(W)$. Moreover, the subgroups
$Sp(W_{i})$ (resp.$Mp(W_{i})$), $1\leq i\leq n$, of $Sp(W)$
(resp. $Mp(W)$) commute each other. Let $Sp(W_{1})\cdots Sp(W_{n})$
(resp. $Mp(W_{1})\cdots Mp(W_{n})$) be their product in $Sp(W)$
(resp. $Mp(W)$). The mapping $\gamma_{1}\cdots\gamma_{n}$
(resp. $\tilde{\gamma}_{1}\cdots\tilde{\gamma_{n}}$) is an isomorphism
(resp. epimorphism) from the direct product
$Sp(W_{1})\times\cdots\times Sp(W_{n})$
(resp. $Mp(W_{1})\times\cdots\times Mp(W_{n})$) onto $Sp(W_{1})\cdots
Sp(W_{n})$ (resp. $Mp(W_{1})\cdots Mp(W_{n})$).

Now consider the representation $\pi_{1}\otimes\cdots\otimes\pi_{n}$
of $Mp(W_{1})\times\cdots\times Mp(W_{n})$. It is trivial on
$\ker\tilde{\gamma}_{1}\cdots\tilde{\gamma_{n}}$ and as such is a
genuine representation of $Mp(W_{1})\cdots Mp(W_{n})$. Moreover, one
has $\pi_{1}\otimes\cdots\otimes\pi_{n}=\pi_{\vert Mp(W_{1})\cdots
  Mp(W_{n})}$.

Let for $1\leq i\leq n$ $S_{i}=\gamma_{i}^{-1}(S)$ and
$\tilde{S}=S_{1}\cdots S_{n}$ which is a subgroup of $ Sp(W_{1})\cdots
Sp(W_{n})$. Then $S_{i}$, $1\leq i\leq n$, and $\tilde{S}$ are the
groups of $k$-points of anisotropic tori over $k$. Moreover, $S$ is a
subgroup of $\tilde{S}$. For $1\leq i\leq n$, we choose a lifting
$\alpha_{i}$ of $S_{i}$ into $Mp(W_{i})$. Then
$(\tilde{\gamma}_{1}\circ\alpha_{1})\cdots(\tilde{\gamma}_{n}\circ\alpha_{n})$
is a lifting of $\tilde{S}$ into $Mp(W)$ which induces the chosen
lifting from $S$ into $Mp(W)$.

\begin{lem}\label{lem3.2}
  (i) If $S$ is admissible, then $\tilde{S}$ is admissible.

  (ii) $\tilde{S}$ is admissible if and only if for any $1\leq i\leq n$
  $S_{i}$ is admissible with respect to $\pi_{i}$.
\end{lem}
\begin{dem}
  The first statement is immediate since $S$ is a subgroup of $\tilde{S}$.

  For the second statement, we consider the decompositions
\begin{equation*}
  \pi_{\vert\tilde{S}}=\oplus_{\chi\in\hat{\tilde{S}}}m(\chi)\chi\mbox{,
  } (\pi_{i})_{\vert
    S_{i}}=\oplus_{\chi\in\hat{S_{i}}}m(\chi)\chi\mbox{, }1\leq
  i\leq n.
\end{equation*}
If $\chi\in\hat{\tilde{S}}$, there exist $\chi_{i}\in\hat{S_{i}}$,
$1\leq i\leq n$, uniquely determined, such that
$\chi(\gamma_{1}(x_{1})\cdots\gamma_{n}(x_{n}))=
\chi_{1}(x_{1})\cdots\chi_{n}(x_{n})$, $x_{i}\in S_{i}$, $1\leq i\leq
n$.  As $\pi_{1}\otimes\cdots\otimes\pi_{n}$, viewed as a
representation of $Mp(W_{1})\cdots Mp(W_{n})$, is the restriction of
$\pi$, it follows easily that $m(\chi)=\prod_{1\leq i\leq
  n}m(\chi_{i})$. Thus, the result.
\end{dem}

\subsection{}\label{3.6}
It follows from Propositions \ref{pr3.1}, \ref{pr3.2} and Lemma
\ref{lem3.2} that to prove Theorem \ref{theo3.1} we may suppose that
all the symplectic subspaces appearing in decomposition \eqref{eq11}
are $S$-irreducible. Thus, we do so until the end of the section.

In this paragraph, we study the momentum map $\phi :
W\rightarrow\mathfrak{s}^{*}$, and prove the equivalence of assertions
(ii), (iii), (iv) of Theorem \ref{theo3.1}.

With our hypothesis, for $1\leq i\leq n$ the torus $S_{i}\subset
Sp(W_{i})$ is irreducible, thus we may choose a maximal irreducible torus
$T^{i}$ in $Sp(W_{i})$ containing $S_{i}$. For $1\leq i\leq n$, let
$\delta_{i}: S\rightarrow Sp(W_{i})$ be the injection such that
$\delta_{i}(s)=s_{\vert W_{i}}$ and let $k'_{i}$, $k''_{i}$,
$u_{i}$, $\beta_{i}$, etc. the objects related to the torus $T^{i}$ in
$Sp(W_{i})$ introduced in Section \ref{2}. Also $\mathfrak{t}^{i}$
denotes the Lie algebra of $T^{i}$.

Let $1\leq i\leq n$. Then $\delta_{i}$, in fact the differential of
$\delta_{i}$ at $1$, maps $\mathfrak{s}$ into
$\mathfrak{t}^{i}=k'u_{i}=k'\nu_{i}$.

If $k'$ is a finite extension of $k$, we consider the non-degenerate
$k$-bilinear form $\langle\, ,\,\rangle$ (or $\langle\,
,\,\rangle_{k'/k}$ if we need to be more precise) on $k'$ such that
\begin{equation*}
\langle x,y\rangle=\sideset{}{_{k'/k}}\tr xy\mbox{, }x,y\in k'
\end{equation*}
and, if $a\in k'$, we denote by $L_{a}$ the $k$-linear form on
$k'$ such that
\begin{equation*}
L_{a}(x)=\langle a,x\rangle\mbox{, } x\in k'.
\end{equation*}

\begin{lem}\label{lem3.3}
The momentum map is given by
\begin{equation}\label{eq12}
\langle\phi(w),\xi\rangle=\frac{1}{2}\sum_{1\leq i\leq
  n}L_{\sideset{}{_{k''_{i}/k'_{i}}}\Norm(u_{i}w_{i})}(\frac{\delta_{i}(\xi)}{u_{i}})
\mbox{, } w\in W\mbox{, }\xi\in\mathfrak{s}
\end{equation}
where for $w\in W$ we write $w=\sum_{1\leq i\leq n}w_{i}$ with $w_{i}\in
W_{i}=k''_{i}$, $1\leq i\leq n$.
\end{lem}
\begin{dem}
Let $\xi\in\mathfrak{s}$ and $w\in W$. Write $w=\sum_{1\leq i\leq
  n}w_{i}$ where $w_{i}\in W_{i}$, $1\leq i\leq n$. As the subspaces
$W_{i}$ are mutually orthogonal with respect to $\beta$ and as
$\beta_{\vert W_{i}}=\beta_{i}$, $1\leq i\leq n$, we have
\begin{align*}
  \langle\phi(w),\xi\rangle &=-\frac{1}{2}\beta(\xi\cdot w,w)\\ &=
  -\frac{1}{2}\sum_{1\leq i\leq n} \beta_{i}(\delta_{i}(\xi)\cdot
  w_{i},w_{i})\\ &= -\frac{1}{4}\sum_{1\leq i\leq
    n}\sideset{}{_{k''_{i}/k}}\tr(u_{i}\delta_{i}(\xi)w_{i}w_{i}^{\tau})
  \\ &=\frac{1}{2}\sum_{1\leq i\leq
    n}\sideset{}{_{k'_{i}/k}}\tr(\frac{\delta_{i}(\xi)}{u_{i}}
  \sideset{}{_{k''_{i}/k'_{i}}}\Norm(u_{i}w_{i}))
\end{align*}
thus the result.
\end{dem}

\begin{pr}\label{pr3.3}
The assertions (ii), (iii) and (iv) in Theorem \ref{theo3.1} are equivalent.
\end{pr}
\begin{dem}
Assume that $0\notin\phi(W\smallsetminus\{0\})$. Suppose that $\phi$
is not a proper map : there exists a compact subset $C$ of
$\mathfrak{s}^{*}$ such that $\phi^{-1}(C)$ is not a compact subset of
$W=\prod_{1\leq i\leq n}k''_{i}$. As $\phi^{-1}(C)$ is closed, it is
not bounded, that is there exists a sequence
$(w_{m})_{m\in\mathbb{N}}$ in $W\smallsetminus\{0\}$ such that, using the
notations introduced in Lemma \ref{lem3.3} :
\begin{equation}\label{eq13}
 \Vert w_{m}\Vert=\Vert(w_{m,1},\ldots,w_{m,n})\Vert:=\sup_{1\leq i\leq n}
q_{i}''^{-v''_{i}(w_{m,i})} 
\end{equation}
tends to infinity with $m$, and
\begin{equation}\label{eq14}
  \phi(w_{m})\in C\mbox{ for all }m\in\mathbb{N}.
\end{equation}

Let $\tilde{e}_{i}$ be the order of ramification of $k''_{i}$ over
$k$. %Then there exists $\eta_{i}\in\mathcal{O}_{i}''^{\times}$ such
%that $\varpi_{i}''^{\tilde{e}_{i}}=\varpi\eta_{i}$.
In case
$w_{m,i}\neq0$, we may write
$v''_{i}(w_{m,i})=\alpha_{m,i}\tilde{e}_{i}+\epsilon_{m,i}$ with
$\alpha_{m,i}$ an integer and $0\leq\epsilon_{m,i}<\tilde{e}_{i}$. Let
$\alpha_{m}=\inf_{1\leq i\leq n\vert w_{m,i}\neq0}\alpha_{m,i}$.

It follows from \eqref{eq13} that
$\lim_{m\rightarrow+\infty}\alpha_{m}=-\infty$. Moreover, one has
\begin{equation*}
  0<\inf_{1\leq i\leq n}q_{i}''^{-\tilde{e}_{i}}\leq
  \sup_{1\leq i\leq n\vert w_{m,i}\neq0}
  q_{i}''^{(\alpha_{m}-\alpha_{m,i})\tilde{e}_{i}-\epsilon_{m,i}}=
  \Vert\varpi^{-\alpha_{m}}w_{m}\Vert\leq1
\end{equation*}
showing that the sequence $(\varpi^{-\alpha_{m}}w_{m})_{m}$ remains in a
compact subset of $W\smallsetminus\{0\}$. Thus, taking a subsequence if
necessary, we may suppose that this sequence converges to an element
$a\neq0$ in $W$. Then we have
\begin{equation*}
  \phi(a)=\lim_{m\rightarrow+\infty}\phi(\varpi^{-\alpha_{m}}w_{m})=
 \lim_{m\rightarrow+\infty}\varpi^{-2\alpha_{m}}\phi(w_{m})=0 
\end{equation*}
and thus, $0\in\phi(W\smallsetminus\{0\})$.

The implication (iii)$\Rightarrow$(iv) is immediate.

In order to prove (iv)$\Rightarrow$(ii), assume that
$0\in\phi(W\smallsetminus\{0\})$, and let $a\in W\smallsetminus\{0\}$
such that $\phi(a)=0$. Then one has
$k^{\times}a\subset\phi^{-1}(\{0\})\cap(W\smallsetminus\{0\})$, showing
that $\phi_{\vert W\smallsetminus\{0\}}$ is not proper on its image.
\end{dem}

\subsection{}\label{3.7}
The end of this section is devoted to the proof of the equivalence of
assertions (i) and (ii) of Theorem \ref{theo3.1}.

For $1\leq i\leq n$ let $\hat{T}^{i}$ be the set of unitary characters
of $T^{i}$ and $\hat{T}^{i}_{0}$ be the subset of the ones which occur
in $\pi_{i\vert T^{i}}$ and let, similarly, $\hat{S}_{0}$ be the set
of unitary characters of $S$ which occur in $\pi_{\vert S}$.

Let $\mathfrak{u}\subset\mathfrak{s}$ be a lattice. We may suppose
that $\mathfrak{u}$ is small enough in order that for
$X\in\mathfrak{u}$ and $1\leq i\leq n$ , $v''_{i}(\delta_{i}(X))>
\frac{v_{i}''(p)}{p-1}$. It follows from Lemma \ref{lem2.2} that
$\exp$ is defined on $\mathfrak{u}$. Then, for any natural integer $m$
we put $U_{m}=\exp\varpi^{m}\mathfrak{u}$. The $U_{m}$,
$m\in\mathbb{N}$, are open subgroups of $S$ which form a basis of
neighbourhoods of $1$.

If $A$ is a set, we denote by $\vert A\vert$ its cardinal.

If $\chi_{i}\in\hat{T}^{i}$, $1\leq i\leq n$, we denote by
$\chi_{1}\cdots\chi_{n}$ the unitary character of
$T_{1}\times\cdots\times T_{n}=T_{1}\cdots T_{n}$ such that
$\chi_{1}\cdots\chi_{n}(x_{1},\ldots
x_{n})=\chi_{1}(x_{1})\cdots\chi_{n}(x_{n})$, $(x_{1},\ldots x_{n})\in
T_{1}\times\cdots\times T_{n}$.

For $U\subset S$  any open subgroup let
\begin{equation*}
\Xi_{U}=\big\{(\chi_{1},\ldots,\chi_{n})\in\hat{T}^{1}_{0}\times\cdots
    \times\hat{T}^{n}_{0}\vert (\chi_{1}\cdots\chi_{n})_{\vert
      U}=1\big\}.
\end{equation*}

\begin{pr}\label{pr3.4}
  The following assertions are equivalent

  (i) the torus $S$ is admissible,

  (ii) for any sufficiently small open subgroup $U$ of $S$,
  $\big\vert\Xi_{U}\big\vert \in \mathbb{N}$,

  (iii) for any $m\in\mathbb{N}$,
    $\big\vert\Xi_{U_{m}}\big\vert\in \mathbb{N}$.
  
\end{pr}
\begin{dem}
It follows from Theorem \ref{theo2.1} (i) that one has
\begin{equation*}
\pi_{\vert T_{1}\cdots
  T_{n}}=\bigoplus_{(\chi_{1},\ldots,\chi_{n})\in\hat{T}^{1}_{0}\times\cdots
  \times\hat{T}^{n}_{0}}\chi_{1}\cdots\chi_{n}.
\end{equation*}   

Therefore $\pi_{\vert S}=\oplus_{\chi\in\hat{S}_{0}}m(\chi)\chi$ where
\begin{equation*}
m(\chi)=\big\vert\big\{(\chi_{1},\ldots,\chi_{n})\in\hat{T}^{1}_{0}\times\cdots
\times\hat{T}^{n}_{0}\vert(\chi_{1}\cdots\chi_{n})_{\vert
  S}=\chi\big\}\big\vert.
\end{equation*}
Now, let $U$ be an open subgroup of $S$. We have
\begin{equation*}
  \sum_{\chi\in\hat{S}_{0},\chi_{\vert U}=1}m(\chi)=
 \big\vert\Xi_{U}\big\vert.
\end{equation*}
The equivalence between (i) and (ii) follows easily from these
considerations and from the fact that on one hand, for any open
subgroup $U$ of $S$, $S/U$ is a finite abelian group and on the other hand,
since $S$ is a totally disconnected group, for any unitary character
$\chi$ of $S$, $\ker\chi$ is an open subgroup.

The equivalence between (ii) and (iii) follows from the fact that
$(U_{m})_{m\in\mathbb{N}}$ is a basis of neighbourhoods of $1$ in $S$
while $\Xi_{U}\subset\Xi_{V}$ if $V\subset U$ are open subgroups of
$S$.
\end{dem}

%\subsection{}
%Let $\kappa:\mathbb{N}\rightarrow\mathbb{N}$ be the map such that
%\begin{equation*}
 % \kappa(j)=\begin{cases}\frac{j}{2}\mbox{ if $j$ is even}\\
  %\frac{j+1}{2}\mbox{ if $j$ is od}
  %\end{cases}
%\end{equation*}

For $U\subset S$ an open subgroup, let
$\kappa_{U}:\Xi_{U}\rightarrow\mathbb{N}^{n}$ be the map defined by
\begin{equation*}
  \kappa_{U}(\chi_{1},\ldots,\chi_{n})=
  (c(\chi_{1}),\ldots,c(\chi_{n})).
\end{equation*}
\begin{co}\label{co3.1}
  The following assertions are equivalent

  (i) the torus $S$ is admissible,

  (ii) for any open subgroup $U$ of $S$,
  $\big\vert\kappa_{U}(\Xi_{U})\big\vert\in\mathbb{N}$,

  (iii) for any $m\in\mathbb{N}$,
  $\big\vert\kappa_{U_{m}}(\Xi_{U_{m}})\big\vert\in\mathbb{N}$.
\end{co}
\begin{dem}
This follows easily from Proposition \ref{pr3.4} and the fact that for
any open subgroup $U$ of $S$, the fibers of $\kappa_{U}$ are finite.
\end{dem}

\subsection{}\label{3.8}
After permutation of indices if necessary, there
exists an integer $0\leq n'\leq n$ such that $k''_{i}$ is
unramified over $k'_{i}$ if and only if $1\leq i\leq n'$.

It follows from the results of the two preceding
paragraphs that proving Theorem \ref{theo3.1} reduces to proving
the following

\begin{theo}\label{theo3.2}
  The following assertions are equivalent

  (i) $0\notin\phi(W\smallsetminus\{0\})$,

  (ii) for any sufficiently small subgroup $U$ of $S$,
  $\big\vert\kappa_{U}(\Xi_{U})\big\vert\in\mathbb{N}$.
\end{theo}
\begin{dem}
Assume (i). Suppose that there exists an open subgroup $U\subset S$
such that $\big\vert\kappa_{U}(\Xi_{U})\big\vert=\infty$. As
$\Xi_{U}\subset\Xi_{V}$ for any open subgroup $V\subset U$, we may
suppose that $U=\exp\mathfrak{u}$ with
$\mathfrak{u}\subset\mathfrak{s}$ a lattice such that for
$X\in\mathfrak{u}$ and $1\leq i\leq n$, $v''_{i}(\delta_{i}(X))\geq
l_{i}$ where $l_{i}$ is the smallest integer such that
$l_{i}\geq\frac{v''_{i}(p)}{p-1}$, $1\leq i\leq n'$ and
$l_{i}\geq\frac{2v''_{i}(p)}{p-1}$, $n'+1\leq i\leq n$.

Let $(J_{m})_{m\in\mathbb{N}}=((j_{m,1},\ldots,j_{m,n}))_{m\in\mathbb{N}}$ be a
sequence of elements all distinct of $\kappa_{U}(\Xi_{U})$. It follows
from Theorem \ref{theo2.1} that for $m\in\mathbb{N}$, $j_{m,i}$ has the
same parity as $\mu_{i}$ for $1\leq i\leq n'$ and is even for
$n'+1\leq i\leq n$.

Let $m\in\mathbb{N}$.  Then for $1\leq i\leq n$, there exists a
unitary character $\chi_{m,i}$ of $T^{i}$ with conductor $j_{m,i}$
which appears in the restriction of the Weil representation to $T^{i}$
and such that $(\chi_{m,1}\cdots\chi_{m,n})_{\vert U}=1$. Moreover, it
follows from Theorem {\color{black}\ref{theo2.1} and Propositions
  \ref{pr2.2}, \ref{pr2.3} and \ref{pr2.4}} that there exists
$a_{m,i}\in\mathcal{O}_{i}''^{\times}\cup\{0\}$ with $a_{m,i}\neq0$ if
$j_{m,i}\geq l_{i}$ and such that
\begin{eqnarray}
  \chi_{i}(\exp X)&=&\psi(-\frac{1}{4}\sideset{}{_{k''_{i}/k}}
  \tr\varpi_{i}'^{\mu_{i}-j_{m,i}}\sideset{}{_{k''_{i}/k}}\Norm(a_{m,i})u_{i}X),\notag\\
  &\phantom{=}&
  X\in\varpi_{i}'^{l_{i}}\mathcal{O}'_{i}\nu_{i}\mbox{, }1\leq i\leq
  n',\label{eq15}\\ \chi_{i}(\exp
  X)&=&\psi(-\frac{(-1)^{\mu_{i}-j_{m,i}}}{4}
  \sideset{}{_{k''_{i}/k}}\tr\varpi_{i}'^{\mu_{i}-j_{m,i}}
  \sideset{}{_{k''_{i}/k}}\Norm(a_{m,i})\varpi''_{i}X),
  \notag\\ &\phantom{=}&
  X\in\varpi_{i}'^{[\frac{l_{i}}{2}]}\mathcal{O}'_{i}\varpi''_{i}\mbox{, }n'+1\leq
  i\leq n.\label{eq16}
\end{eqnarray}

Now let $w_{m}=(w_{m,1},\ldots,w_{m,n})\in W$ such that
\begin{equation*}
  w_{m,i}=\begin{cases}
  \varpi_{i}'^{\frac{\mu_{i}-j_{m,i}}{2}}a_{m,i}\mbox{, }1\leq i\leq n',\\
  \varpi_{i}''^{\mu_{i}-j_{m,i}}a_{m,i}\mbox{, }n'+1\leq i\leq n.
  \end{cases}
\end{equation*}
Then one  easily checks that
\begin{equation*}
(\chi_{1}\cdots\chi_{n})(\exp
  X)=\psi(\langle\phi(w_{m}),X\rangle)\mbox{, }
  X\in\mathfrak{u}.
\end{equation*}
As $(\chi_{1}\cdots\chi_{n})_{\vert U}=1$, one deduces that
$\langle\phi(w_{m}),X\rangle\in\varpi^{\lambda_{\psi}}\mathcal{O}$,
$X\in\mathfrak{u}$. This means that $\phi(w_{m})\in\mathfrak{u}^{*}$
where $\mathfrak{u}^{*}=\{\lambda\in\mathfrak{s}^{*}
\vert\langle\lambda,X\rangle\in\varpi^{\lambda_{\psi}}\mathcal{O}\mbox{,
}X\in\mathfrak{u}\}$ is the dual lattice in $\mathfrak{s}^{*}$ of the
lattice $\mathfrak{u}$ with respect to the additive character $\psi$.

For $m\in\mathbb{N}$ let us write
\begin{eqnarray*}
  \frac{j_{m,i}-\mu_{i}}{2}&=&\alpha_{m,i}\tilde{e}_{i}+\epsilon_{m,i}
  \mbox{ with } \alpha_{m,i}\in\mathbb{Z}\\ &\phantom{=}&\mbox{ and }0\leq
  \epsilon_{m,i}<\tilde{e}_{i}\mbox{, }1\leq i\leq n',\\
  j_{m,i}-\mu_{i}&=&\alpha_{m,i}\tilde{e}_{i}+\epsilon_{m,i}
  \mbox{ with } \alpha_{m,i}\in\mathbb{Z}\\ &\phantom{=}&\mbox{ and }0\leq
  \epsilon_{m,i}<\tilde{e}_{i}\mbox{, }n'+1\leq i\leq n.
\end{eqnarray*}
and let $\alpha_{m}=\sup_{1\leq i\leq n}\alpha_{m,i}$,
$\overline{w}_{m}=\varpi^{\alpha_{m}}w_{m}$.

As $\vert\kappa_{U}(\Xi_{U})\vert=\infty$, the sequence
$(\alpha_{m})_{m\in\mathbb{N}}$ is not bounded and there exists $1\leq
i_{0}\leq n$ such that $a_{m,i_{0}}\neq0$ and
$\alpha_{m}=\alpha_{m,i_{0}}$ for infinitely many $m$.

Then, taking a subsequence if necessary, we may assume that these
conditions are satisfied for all $m\in\mathbb{N}$ and that
$\lim_{m\rightarrow+\infty}\alpha_{m}=+\infty$. It thus follows that
one has
\begin{equation*}
  1\leq
  \sup_{1\leq i\leq n\vert w_{m,i}\neq0}
  q_{i}''^{\epsilon_{m,i}-(\alpha_{m}-\alpha_{m,i})\tilde{e}_{i}}=
  \Vert\overline{w}_{m}\Vert\leq\sup_{1\leq i\leq n}q''^{\epsilon_{m,i}}
\end{equation*}
showing that the sequence $(\overline{w}_{m})_{m}$ remains in a
compact subset of $W\smallsetminus\{0\}$. Thus, we may assume that there
exists $a\in W\smallsetminus\{0\}$ such that
$a=\lim_{m\rightarrow+\infty}\overline{w}_{m}$. But at the same time
$\phi(\overline{w}_{m})=\varpi^{2\alpha_{m}}\phi(w_{m})
\in\varpi^{2\alpha_{m}}\mathfrak{u}^{*}$, which implies
$\phi(a)=\lim_{m\rightarrow+\infty}\phi(\overline{w}_{m})=0$. This
shows that $0\in\phi(W\smallsetminus\{0\})$.

Conversely, assume (ii) and suppose
$0\in\phi(W\smallsetminus\{0\})$. Let
$w=(w_{1},\ldots,{\color{black}w_{n}})\in W\smallsetminus\{0\}$ such
that $\phi(w)=0$. For $1\leq i\leq n$ let us put
{\color{black}$w_{i}=\varpi''^{\upsilon_{i}}a_{i}$ with $\upsilon_{i}\in\mathbb{Z}$
  and $a_{i}\in\mathcal{O}''^{\times}$ if $w_{i}\neq0$, $\upsilon_{i}=0$ and
  $a_{i}=0$ if not,} and for $m\in\mathbb{N}$ let $w_{m}=\varpi^{m}w$.

Remark that for $1\leq i\leq n$, there exists
$\rho_{i}\in\mathcal{O}_{i}''^{\times}$ such that
$\varpi=\varpi_{i}''^{\tilde{e_{i}}}\rho_{i}$.

For $m\in\mathbb{N}$ let us put
\begin{eqnarray*}
  j_{m,i}&=&\begin{cases}
    \mu_{i}-2(m\tilde{e}_{i}-{\color{black}\upsilon_{i}})\mbox{,
    }1\leq i\leq
    n',\\ \mu_{i}-m\tilde{e}_{i}-{\color{black}\upsilon_{i}}\mbox{,
    }n'+1\leq i\leq n,
  \end{cases}\\
  a_{m,i}&=&
\rho_{i}^{m}u_{i}^{-1}a_{i}\mbox{, }1\leq i\leq n.
\end{eqnarray*}
Let $\mathfrak{u}\subset\mathfrak{s}$ be a lattice such that for
$X\in\mathfrak{u}$ and $1\leq i\leq n$, $v''_{i}(\delta_{i}(X))\geq
l_{i}$ and let $U=\exp(\mathfrak{u})$.
Then, if $m$ is big enough, there exist characters
$\chi_{m,i}\in\hat{T}^{i}$ satisfying formula \eqref{eq15}
(resp. \eqref{eq16}) for $1\leq i\leq n'$ (resp. $n'+1\leq i\leq
n$). Then one checks easily that these characters satisfy
\begin{equation*}
(\chi_{m,1}\cdots\chi_{m,n})(\exp
  X)=\psi(\langle\varpi^{m}\phi(w),X\rangle)=1\mbox{, }X\in\mathfrak{u}
\end{equation*}
and moreover, the terms of the sequence
$(c(\chi_{m,1}),\ldots, c(\chi_{m,n}))_{m}$ are all
distinct. It follows easily that $\vert\kappa_{U}(\Xi_{U})\vert=\infty$.
\end{dem}

\section{Subtori of maximal irreducible tori}\label{4}
\subsection{}\label{4.1}
Let $T\subset Sp(W)$ be a maximal irreducible torus. Recall the
notations introduced in paragraphs \ref{2.1} and \ref{3.2}. On one
hand it follows from Theorem \ref{theo2.1} that $T$ is admissible, and
on the other hand by Lemma \ref{lem3.3} the momentum map for the
symplectic action of $T$ in $W$ is given by
\begin{equation}\label{eq17}
  \langle\phi(w),X\rangle=\frac{1}{2}L_{\sideset{}{_{k''/k'}}\Norm(uw)}(\frac{X}{u})
    \mbox{, }w\in W\mbox{, }X\in\mathfrak{t}.
\end{equation}
This shows that $T$ and its momentum map satisfy assertions (i) to
(iv) of Theorem \ref{3.1}.

In this section we are interested in the following question
\begin{qu}\label{qu1}
Does there exist a strict subtorus $S$  of $T$ which
%acts irreducibly on $W$, that is to say is irreducible, and
is admissible ?
\end{qu}

First, in this paragraph we show that under certain conditions the
answer is negative.

Recall that $k'$ is equipped with the non-degenerate $k$-bilinear form
$\langle x,y\rangle=\tr_{k'/k} xy$. We consider the orthogonality with
respect to it.

Let us begin with a lemma.
\begin{lem}\label{lem4.1}
  Let $S\subset T$ be a torus. The following assertions are equivalent

  (i) $S$ is admissible,

  (ii) $\sideset{}{_{k''/k'}}\Norm(k''^{\times})\cap
  (\mathfrak{s}/u)^{\perp}=\emptyset$.
\end{lem}
\begin{dem}
This follows readily from the fact that the momentum map $\phi$ of the
symplectic action of $S$ in $W$ is also given by formula \eqref{eq17}
where $\mathfrak{t}$ is replaced by $\mathfrak{s}$.
\end{dem}

\begin{pr}\label{pr4.1}
  Let $S$ be a subtorus of the maximal irreducible torus $T$ which is
  admissible.

  (i) If $k''$ is unramified over $k'$ and the order of
  ramification of $k'$ over $k$ is odd, then $S=T$.

  (ii) If $k''$ is ramified over $k'$ and the modular degree of
  $k'$ over $k$ is odd, then $S=T$.
\end{pr}
\begin{dem}
  (i) In this case we have
$\Norm_{k''/k'}(k''^{\times})=\{x\in k'\vert
  v'(x)\equiv0\pmod{2}\}$
As $v'(\varpi)\equiv1\pmod{2}$, every
  non-zero $k$-subspace of $k'$ contains at least an element with even
  valuation : if $0\neq x\in k'$ then $v'(x)$ or $v'(\varpi x)$ is
  even. Thereby condition (ii) of Lemma \ref{lem4.1} is satisfied if
  and only if $(\mathfrak{s}/u)^{\perp}=0$, that is
  $\mathfrak{s}=\mathfrak{t}$.

  (ii) As $p\neq2$ it follows from \cite[Proposition 8,
    p. 32]{weil-1967} that
  $(1+\varpi'\mathcal{O}')^{2}=1+\varpi'\mathcal{O}'$. If
  $a\in\mathcal{O}''^{\times}$, we have $a=\xi+\eta\varpi''$ with
  $\xi\in\mathcal{O}'^{\times}$ and $\eta\in\mathcal{O}'$ and thus
  $\Norm_{k''/k'}(a)=
  \xi^{2}(1-(\eta\xi^{-1})^{2}\varpi')\in(\mathcal{O}'^{\times})^{2}$.
 
Thus we have
\begin{equation}\label{eq18}
  \sideset{}{_{k''/k'}}\Norm(k''^{\times})=\cup_{m\in\mathbb{Z}}(-\varpi')^{m}
(\mathcal{O}'^{\times})^{2}.
\end{equation}
Let $\gamma$ be a primitive $(q-1)$-th root of unity ($\gamma$ is
necessarily in $\mathcal{O}^{\times}$). Then there exists $1\leq l\leq
q-1$ an integer prime to $q-1$ such that
$\gamma=\gamma'^{\frac{q'-1}{q-1}l}$, where $\gamma'$ is a primitive
root of unity in $k'$. Remark that $l$ is odd. Let us denote by $f$
the modular degree of $k'$ over $k$. We have
\begin{equation*}
\frac{q'-1}{q-1}=\frac{q^{f}-1}{q-1}=q^{f-1}+\cdots+q+1.
\end{equation*}
Thus, $\gamma$ is a square in $\mathcal{O}'$ if and only if $f$ is
even.

As $f$ is odd, it follows that $\gamma\in\mathcal{O}'^{\times}
\smallsetminus(\mathcal{O}'^{\times})^{2}$. Now let $x\in k'$ be a non-zero
element, and write $x=(-\varpi')^{v'(x)}\alpha$ with
$\alpha\in\mathcal{O}'^{\times}$. Then $\alpha$ or $\alpha\gamma$ is a
square in $\mathcal{O}'$, and $x$ or $x\gamma$ belongs to
$\sideset{}{_{k''/k'}}\Norm(k''^{\times})$. Thus, every non-zero
$k$-subspace of $k'$ meets
$\sideset{}{_{k''/k'}}\Norm(k''^{\times})$. We conclude as in case
(i).
\end{dem}

Until the end of this paragraph we suppose that $k''$ is ramified
over $k'$ and $k'$ is unramified over $k$ of even degree, say
$2f$. Moreover, we may suppose that
$\varpi''^{2}=\varpi'=\varpi\gamma$ where $\gamma\in\mu_{q'-1}$.

We denote by $\Gamma'$ the Galois group of $k'$ over $k$ which is
cyclic of order $2f$ and by $\tau$ the non-trivial element of the
Galois group of $k''$ over $k'$. Recall that $\Gamma'$ is generated by
the Frobenius automorphism $\theta$ : $\theta$ is the unique
$k$-automorphism of $k'$ such that $\theta(u)=u^{q}$, $u\in\mu_{q'-1}$
(here $q'=q^{2f}$).

\begin{lem}\label{lem4.2}
The extension $k''$ is Galois over $k$. The kernel of the natural
morphism from the Galois group $\Gamma$ of $k''$ over $k$ onto
$\Gamma'$ is $\{1,\tau\}$. Moreover, if $\gamma$ is a square in $k'$,
$\Gamma'$ is naturally identified with a subgroup of $\Gamma$ and
$\Gamma$ is the direct product of $\{1,\tau\}$ and of $\Gamma'$, and,
if $\gamma$ is not a square in $k'$, $\Gamma$ is cyclic of order $4f$.
\end{lem}
\begin{dem}
Let $\Omega$ be an algebraic closure of $k$ containing $k''$, and
$\sigma$ be an embedding over $k$ of $k''$ into $\Omega$. Then the
restriction to $k'$ of $\sigma$ is a power, say $l$, of the Frobenius
of $k'$ over $k$. It follows that one has
$\sigma(\varpi''^{2})=\varpi\gamma^{q^{l}}=\varpi''^{2}\gamma^{q^{l}-1}$. Thus
we have $\sigma(\varpi'')=\pm\varpi''\gamma^{\frac{q^{l}-1}{2}}\in
k''$, showing that $k''$ is Galois over $k$. It follows that
$\{1,\tau\}$ is an invariant subgroup of $\Gamma$, and that it is the
kernel of the natural morphism from $\Gamma$ onto $\Gamma'$.

If $\gamma$ is a square in $k'$, we may take $\varpi'=\varpi$. Then,
we may embed $\Gamma'$ into $\Gamma$ by deciding that its elements
act trivially on $\varpi''$. Thus in this case one has
$\Gamma=\{1,\tau\}\times\Gamma'$.

Suppose now that $\gamma$ is not a square in $k'$, and let
$\sigma\in\Gamma$ such that its restriction to $k'$ is the
Frobenius. We may suppose that $\gamma$ is a $(q'-1)$-primitive root of
the unity. According to what we saw above, we may choose $\sigma$ such
that $\sigma(\varpi'')=\gamma^{\frac{q-1}{2}}\varpi''$. Then we see
by induction that
$\sigma^{l}(\varpi'')=\gamma^{\frac{q^{l}-1}{2}}\varpi''$,
$l\in\mathbb{N}$. In particular, if $l=2f$, we get
$\sigma^{2f}(\varpi'')=\gamma^{\frac{q'-1}{2}}\varpi''
=-\varpi''$. As $\sigma^{2f}_{\vert k'}=Id$, it follows that
$\sigma^{2f}=\tau$. The lemma
follows easily.
\end{dem}

\begin{pr}\label{pr4.2}
Suppose that $k''$ is ramified over $k'$, $k'$ is unramified over $k$
of degree $2f$ and that $\varpi''^{2}=\gamma\varpi$ with
$\gamma\in\mu_{q'-1}$ not a square. Then no proper subtorus of $T$ is
admissible.
\end{pr}
\begin{dem}
Let $V$ be a non-zero subspace of $k'$ over $k$. It follows from
formula \eqref{eq18} that
$(\mathcal{O}'^{\times})^{2}\cup-\gamma\varpi(\mathcal{O}'^{\times})^{2}\subset
\Norm_{k''/k'}(k''^{\times})$. Let $x$ be a non-zero element of
$V$. Multiplying it by a power of $\varpi$ , we may suppose that
$x\in\mathcal{O}'^{\times}\cap V$. If
$x\in(\mathcal{O}'^{\times})^{2}$, then $x\in\im\Norm_{k''/k'}$. If
not, $x\in\gamma(\mathcal{O}'^{\times})^{2}$, and one has $-\varpi x\in
V\cap\im\Norm_{k''/k'}$. The proposition follows.
\end{dem}

\subsection{}\label{4.2}
In the remaining paragraphs of this section we show that the answer to
question \ref{qu1} may be positive when the conditions in
Propositions \ref{pr4.1} and \ref{pr4.2} are not satisfied.

In the sequel, when we speak of a torus $S$ defined over $k$, we
mean an algebraic torus defined over $k$ as well as the group of
$k$-points of such a torus. The choice between the two meanings will
be clear from the context. For instance, the group of rational
characters of $S$ is the group of rational character of the algebraic
torus $S$ while the unitary characters of $S$ are the ones of its
group of $k$-points.

If $S$ is a torus defined over $k$, let us denote by $\mathbf{X}(S)$
its group of rational characters which is a free $\mathbb{Z}$-module
isomorphic to $\mathbb{Z}^{\dim S}$. Let $K$ be a finite Galois
extension of $k$ over which $S$ splits and let $\Gamma$ be the Galois
group of $K$ over $k$. Then $\mathbf{X}(S)$ is also a module over the
group ring $\mathbb{Z}[\Gamma]$.

The mapping $U\mapsto\{\chi\in\mathbf{X}(S)\vert\chi_{\vert U}=1\}$ is
a reversing inclusion bijection from the set of subtori of $S$ onto
the set of $\mathbb{Z}[\Gamma]$-submodules $M\subset\mathbf{X}(S)$
such that the quotient module $\mathbf{X}(S)/M$ is torsion-free as a
$\mathbb{Z}$-module. Under this bijection, maximal subtori of $S$
correspond to minimal $\mathbb{Z}[\Gamma]$-submodules with this
property. Moreover, in all cases we have $\dim
U=\dim_{\mathbb{Z}}\mathbf{X}(S)/M$.

In the sequel, identifying them with their differential at $1$, we
also consider the elements of $\mathbf{X}(S)$ as $k$-linear maps from
the Lie algebra $\mathfrak{s}$ of $S$ into $K$. Then,
$\mathbf{X}(S)$ becomes a $\mathbb{Z}$-submodule of
$\mathfrak{s}^{*}_{K}$, the dual space over $K$ of
$\mathfrak{s}_{K}=K\otimes_{k}\mathfrak{s}$, whose
rank is $\dim S$.

Now we look at $T$ a maximal irreducible torus in $Sp(W)$. It is a
subtorus of the torus whose group of $k$-points is $k''^{\times}$. Let
$\tilde{k''}$ be a Galois closure of $k''$ over $k$. Then
$k''^{\times}$ and $T$ split over $\tilde{k''}$. Let $\Gamma$
(resp. $\Gamma_{0}$, $\Gamma_{1}$) be the Galois group of
$\tilde{k}''$ over $k$ (resp. $k''$, $k'$). Then
$\mathbf{X}(k''^{\times})$ is the $\mathbb{Z}[\Gamma]$-module
$\mathbb{Z}[\Gamma/\Gamma_{0}]$ (see \cite[p. 54,
  Example]{platonov-rapinchuk-1994}). Moreover, one has
$\Gamma_{1}/\Gamma_{0}=\{1,\tau\}\subset\Gamma/\Gamma_{0}$ and
$\mathbf{X}(T)=\mathbb{Z}[\Gamma/\Gamma_{0}]/\mathbb{Z}[\Gamma](1+\tau)$.

In the following, if $k'$ is an extension of finite degree of $k$, we
consider adjointness of $k$-linear endomorphisms of $k'$ and
orthogonality with respect to the $k$-bilinear form
$\langle\,,\,\rangle_{k'/k}$. Also, if $\sigma$ is an automorphism of
$k'$ over $k$, we shall denote by $k'^{\sigma}$ the subfield of the
elements of $k'$ which are invariant under $\sigma$.

\begin{lem}\label{lem4.3}
  Let $S$ be an algebraic subtorus of $T$ over $k$ with Lie algebra
  $\mathfrak{s}$. The restriction to $\mathfrak{s}/u$ of the
  bilinear form $\langle\, ,\,\rangle_{k'/k}$ is non-degenerate. Thus,
  we have
\begin{equation*}
k'=\mathfrak{s}/u\oplus(\mathfrak{s}/u)^{\perp}.
\end{equation*}
\end{lem}
\begin{dem}
Letting $k''$ act on itself by multiplication, naturally embed it in the
Lie algebra $\End_{k}(k'')$. In that case $\tr_{k''/k}$ is the restriction
to $k''$ of the natural trace in $\mathfrak{gl}_{k}(k'')$. Then the lemma
follows easily from \cite[p. 127, Ex. 11]{bourbaki-b-1960}.
\end{dem}

\subsection{}\label{4.3}
In this paragraph, we suppose that $k''$ is ramified over $k'$, $k'$
is an unramified extension of $k$ of degree $2f$, and
$\varpi''^{2}=\varpi'$. We keep the notations introduced in Lemma
\ref{lem4.2} and just before it. According to this Lemma, we may
embed $\Gamma'$ into $\Gamma$ by deciding that its elements act
trivially on $\varpi''$ and, doing so, we have
$\Gamma=\Gamma'\times\{1,\tau\}$.

According to what we have seen in paragraph \ref{4.2}, we have
$\mathbf{X}(T)=\mathbb{Z}[\Gamma]/\mathbb{Z}[\Gamma](1+\tau)$. %where
%$\mathbb{Z}[\Gamma]_{+}=
%\oplus_{l\in\mathbb{Z}/2f\mathbb{Z}}\mathbb{Z}\theta^{l}(1+\tau)$.
Moreover, one has the decomposition in $\mathbb{Z}[\Gamma']$-modules
$\mathbb{Z}[\Gamma]=\mathbb{Z}[\Gamma']
\oplus\mathbb{Z}[\Gamma'](1+\tau)$. It follows that, there is a
natural isomorphism of $\mathbb{Z}[\Gamma]$-modules between
$\mathbf{X}(T)$ and $\mathbb{Z}[\Gamma']$, $\mathbb{Z}[\Gamma']$
acting naturally on the left on itself and $\tau$ acting as $-Id$. We
identify $\mathbf{X}(T)$ and $\mathbb{Z}[\Gamma']$ via this
isomorphism.

In the following, it may be also convenient to consider the elements
of $\mathbb{Z}[\Gamma']$ as endomorphisms of the vector-space
$k'\varpi''$ over $k$.

Then, the $\mathbb{Z}[\Gamma]$-submodules $M$ of
$\mathbf{X}(T)=\mathbb{Z}[\Gamma']$ such that $\mathbf{X}(T)/M$ is
torsion free over $\mathbb{Z}$ (resp. torsion free over $\mathbb{Z}$
and minimal for this property) are the
$\mathbb{Z}[\Gamma']$-submodules of $\mathbb{Z}[\Gamma']$ such that
$\mathbb{Z}[\Gamma']/M$ is torsion free over $\mathbb{Z}$
(resp. torsion free over $\mathbb{Z}$ and minimal for this
property). These submodules are the ones described just below.

For $d$ a divisor (resp. $\underline{d}$ a subset of the set of
divisors) of $2f$, we consider the $\mathbb{Z}[\Gamma']$-submodule
(resp. $\mathbb{Z}[\Gamma']$-submodules) $M_{d}$
(resp. $M_{\underline{d}}$ and $\overline{M}_{\underline{d}}$) of
$\mathbb{Z}[\Gamma']$ introduced in the appendix (see definitions
\ref{defA1} and \ref{defA2}).  Then we denote by $S_{\underline{d}}$ the
algebraic subgroup over $k$ of $T$ such that
$S_{\underline{d}}=\cap_{\chi\in
  \overline{M}_{\underline{d}}}\ker\chi$ and by
$\mathfrak{s}_{\underline{d}}$ is Lie algebra. We clearly have
$\mathfrak{s}_{\underline{d}}=\cap_{\chi\in
  M_{\underline{d}}}\ker\chi$. If $d$ is a divisor of $2f$ we simply
put $S_{d}=S_{\{d\}}$ and $\mathfrak{s}_{d}=\mathfrak{s}_{\{d\}}$. The
following lemma follows from the above considerations.

\begin{lem}\label{lem4.4}
The $S_{d}$ for $d$ a divisor of $2f$ are the maximal proper subtori
of $T$ and the $S_{\underline{d}}$ for $\underline{d}$ a subset of the
set divisors of $2f$ are the proper subtori of $T$. The map
$\underline{d}\mapsto S_{\underline{d}}$ is a bijection from the set
of divisors of $2f$ onto the set of subtori of $T$ defined over
$k$. In particular, we have $S_{\emptyset}=T$ and
$S_{\underline{d}}=\{1\}$ if $\underline{d}$ is the set of divisors of
$2f$.
\end{lem}

Now, until the end of this paragraph, we suppose that $f$ is an odd
prime number. It follows that the $\mathbb{Z}[\Gamma']$-submodules $M$
of $\mathbf{X}(T)=\mathbb{Z}[\Gamma']$ such that $\mathbf{X}(T)/M$ is
torsion free and minimal for this property are the $M_{d}$ with
$d=1,2,f,2f$ and the maximal subtori of $T$ are the $S_{d}$, for the
same values of $d$. We recall that $\dim M_{d}=\varphi(d)$ where $\varphi$ is
the Euler function.

\begin{pr}\label{pr4.3}
The proper subtori $S_{2}$ and $S_{2f}$ of $T$ are admissible if and
only if $q\equiv1\pmod{4}$.  Moreover, $S_{2}$ (resp. $S_{2f}$) is of
codimension $1$ (resp. $f-1$) in $T$.
\end{pr}
\begin{dem}
(i) case $S_{2}$. We have $M_{2}=\mathbb{Z}(Y_{0}-Y_{1})$ with
  $Y_{i}=\sum_{l\in\mathbb{Z}/f\mathbb{Z}}\theta^{2l+i}$,
  $i=0,1$. Then the Lie algebra of $S_{2}$ is
  $\mathfrak{s}_{2}=\ker(Y_{0}-Y_{1})$, and
  $\dim\mathfrak{s}_{2}=\dim_{\mathbb{Z}}
  \mathbb{Z}[\Gamma']/M_{2}=2f-1$ so that $S_{2}$ is of codimension
  $1$ in $T$.  Moreover, considering $Y_{0}-Y_{1}$ as an
    endomorphism of $k'$, we have
    $\mathfrak{s}_{2}=\ker(Y_{0}-Y_{1})\varpi''$. But, $Y_{0}-Y_{1}$
    is a self-adjoint $k$-endomorphism of $k'$. Thus, one has
    $(\mathfrak{s}_{2}/\varpi'')^{\perp} =\im(Y_{0}-Y_{1})$.

  But, as $\theta(Y_{0}-Y_{1})=-(Y_{0}-Y_{1})$, we get
  $\im(Y_{0}-Y_{1})\subset k'^{\theta^{2}}$. Yet, $k'^{\theta^{2}}$ is
  the extension of degree two of $k$ contained in $k'$ :
  $k'^{\theta^{2}}=k(x)$ where $x=y^{\frac{q+1}{2}}$ with
  $y\in\mu_{q'-1}$ a primitive root of order $q^{2}-1$. Now, it is
  clear that $\im(Y_{0}-Y_{1})=kx$.

  Let $z\in k'$ be a $(q'-1)$-th primitive root of unity. Then
  $t=z^{\frac{q'-1}{q-1}}$ is a $(q-1)$-th primitive root of unity, so
  that $t\in\mathcal{O}^{\times}$, and as
  $\frac{q'-1}{q-1}=1+q+\cdots+q^{2f-1}\equiv0\pmod{2}$, $t$ is a
  square in $\mathcal{O}'$. It follows that
  $(\mathcal{O}'^{\times})^{2}\cap k=\mathcal{O}^{\times}$. Then as
  $x$ is a square in $\mathcal{O}'^{\times}$ if and only if
  $q\equiv3\pmod{4}$, it follows from \eqref{eq18} that $S_{2}$ is
  admissible if and only if $q\equiv1\pmod{4}$.

    %Now, $x\in\ker\tr_{k'/k}$, since $\tr_{k(x)/k}x=x+x^{q}=0$, and
    %thus, $1\in(kx)^{\perp}=\mathfrak{s}_{2}/\varpi''$. Thereby,
    %$k(\mathfrak{s}_{2})$ contains the subspace
    %$\mathfrak{s}_{2}\oplus\varpi''\mathfrak{s}_{2}$ of dimension
    %$4f-2$. Thus $k(\mathfrak{s}_{2})=k''$ and $S_{2}$ is irreducible.

(ii) case $S_{2f}$.  Here $M_{2f}=\oplus_{0\leq l\leq
      f-2}\mathbb{Z}Y\theta^{l}$ with
    $Y=(\theta+1)(\theta^{f}-1)$. The Lie algebra $\mathfrak{s}_{f}$
    of $S_{f}$ has dimension $f+1$, and
    $\mathfrak{s}_{f}/\varpi''=\ker Y$, so that $S_{f}$ is of
    codimension $f-1$ in $T$. Now consider $k'^{\theta^{f}}$ which is
    the extension of degree $f$ of $k$ contained in $k'$. With the
    notations introduced in the preceding case, we have
    $k'=k'^{\theta^{f}}(x)$.

    Moreover, one has $\theta(x)=x^{q-1}x= y^{\frac{q^{2}-1}{2}}x=-x$,
    and thus $\theta^{f}(x)=-x$. Let $\lambda\in k'$. Writing
    $\lambda=\alpha+\beta x$ with $\alpha,\beta\in k'^{\theta^{f}}$,
    we see that $\lambda\in\ker Y$ if and only if
    $\theta\beta=\beta$. It follows that
    $\mathfrak{s}_{f}/\varpi''=k'^{\theta^{f}}\oplus kx$ and
    $(\mathfrak{s}_{f}/\varpi'')^{\perp}=x\ker\tr_{k'^{\theta^{f}}/k}$.

    However, as $k'$ is a quadratic extension of $k'^{\theta^{f}}$, we
    have $(\mathcal{O}'^{\times})^{2}\cap
    k'^{\theta^{f}}=\mathcal{O}_{k'^{\theta^{f}}}^{\times}$ (of course
    $\mathcal{O}_{k'^{\theta^{f}}}$ stands for the ring of integers of
    the field $k'^{\theta^{f}}$).

    But as $x$ is not a square in $\mathcal{O}''$ if and only if
    $q\equiv1\pmod{4}$, it follows that $S_{2f}$ is
    admissible if and only if $q\equiv1\pmod{4}$.
   %Then one checks readily that the $k$-subspace of dimension $2f+2$
   %$k'^{\theta^{f}}\oplus kx\oplus(k'^{\theta^{f}}\oplus kx)\varpi''$
   %is contained in $k(\mathfrak{s}_{f})$. It follows that $S_{2f}$
   %is irreducible.
\end{dem}

\begin{rem}\label{rem4.1}
We have not included the proof, but it may be proven that neither
$S_{1}$, nor $S_{f}$ is admissible.
\end{rem}

\begin{rem}\label{rem4.2}
Contrary to the case examined in the following paragraph we are not
able to characterize among the sets of divisors of $2f$, with $f$ any
positive integer, the ones such that the subtorus $S_{\underline{d}}$
is admissible.
\end{rem}

\subsection{}\label{4.4}

In this paragraph, we suppose that $k''$ is unramified over $k'$ and
that $k'$ is tamely and totally ramified over $k$ with even order of
ramification $2e$. The fact that $k'$ is tamely ramified over $k$
implies that $2e$ is prime to $p$ and that we may choose the
uniformizers $\varpi$ and $\varpi'$ such that $\varpi'^{2e}=\varpi$.
Moreover, $\varpi'$ is also a uniformizer for $k''$.

Let $k''_{1}\subset k''$ be the maximal unramified extension of $k$,
$\gamma\in k''$ be a primitive element of $\mu_{q^{2}-1}$, and let
$\nu=\gamma^{\frac{q+1}{2}}$. Then $\gamma$ belongs to $k''_{1}$,
$\nu^{2}$ is in $k$, $k''_{1}=k(\nu)$ and
$k''=k''_{1}(\varpi')=k'(\nu)$.

Let $\Omega$ be an algebraic closure of $k''$ and let $\zeta\in\Omega$
be a primitive $2e$-th root of the unity. Then
$\tilde{k}''=k''(\zeta)$ (resp. $\tilde{k}'=k'(\zeta)$) is the minimal
Galois extension over $k$ of $k''$ (resp. $k'$) contained in
$\Omega$. Moreover, $\tilde{k}''_{1}=k''_{1}(\zeta)$
(resp. $\tilde{k}'_{1}=k(\zeta)$) is the maximal unramified extension
of $k$ contained in $\tilde{k}''$ (resp. $\tilde{k}'$). Then we have
$\tilde{k}''_{1}=\tilde{k}'_{1}(\nu)$.

One easily checks that the degree of $\tilde{k}''_{1}$
(resp. $\tilde{k}'_{1}$) over $k$ is $2f$ (resp. $g$) where $f$
(resp. $g$) is the smallest positive integer such that $q^{2f}-1$
(resp.  $q^{g}-1$) is divisible by $2e$ : in both cases, such an
integer exists because $q$ and $2e$ are relatively prime numbers from
which follows that the image of $q$ in the finite ring
$\mathbb{Z}/2e\mathbb{Z}$ is invertible. As the order of $\gamma$ is
$q^{2}-1$, it follows that $\tilde{k}''_{1}=\tilde{k}'_{1}$ if and
only if $g\geq 2$ that is if and only if $e$ does not divide
$\frac{q-1}{2}$. Moreover one has $\tilde{k}''=k''$
(resp. $\tilde{k}'=k'$) if and only if $2e$ divides $q^{2}-1$
(resp. $q-1$).

Let $\Gamma$ (resp. $\Gamma_{0}$, $\Gamma_{1}$, $\Gamma'$) be the
Galois group of $\tilde{k}''$ over $k$ (resp. $k''$, $k'$,
$\tilde{k}''_{1}$).

It is clear that $\tilde{k}''$ over $k'$ is unramified of even degree,
say $2f$. Thus $\Gamma_{1}$ is a cyclic group of order $2f$ generated
by the Frobenius automorphism $\theta$ of $\tilde{k}''$ over $k'$. On
the other side, $\tilde{k}''$ over $\tilde{k}''_{1}$ is totally
ramified of order $2e$ and its Galois group $\Gamma'$ is cyclic of
order $2e$ generated by the element $\sigma$ such that
$\sigma(\varpi')=\zeta\varpi'$.

As in paragraph \ref{4.2}, we consider the homogeneous space
$\Gamma/\Gamma_{0}$. It contains $\Gamma_{1}/\Gamma_{0}$, which is the
order $2$ Galois group of $k''$ over $k'$ whose non trivial element is
$\tau$. Then we have
$\mathbf{X}(T)=\mathbb{Z}[\Gamma/\Gamma_{0}]/\mathbb{Z}[\Gamma](1+\tau)$.

\begin{pr}\label{pr4.4}
  We have $\Gamma=\Gamma_{1}\Gamma'$,
  $\Gamma_{1}\cap\Gamma'=\{1\}$, $\Gamma'$ is an invariant
  subgroup of $\Gamma$, and we have the commutation relation
  \begin{equation}\label{eq20}
\theta\sigma\theta^{-1}=\sigma^{q}.
  \end{equation}
  
Moreover, the homogeneous space $\Gamma/\Gamma_{0}$ is isomorphic to
$\Gamma'\times\{1,\tau\}$, $\Gamma'$ acting by left translation
  on the first factor and trivially on the second one, while $\theta$
  acts by $\sigma^{l}\mapsto\sigma^{lq}$ on $\Gamma'$ and permutes
  $1$ and $\tau$.
\end{pr}
\begin{dem}
As $\tilde{k}''_{1}$ is Galois over $k$, the map
$\alpha\mapsto\alpha_{\vert\tilde{k}''_{1}}$ is an epimorphism from
$\Gamma$ onto the Galois group $\Gamma_{2}$ of $\tilde{k}''_{1}$ over
$k$ with kernel $\Gamma'$. As $\tilde{k}''$ is generated by
$\tilde{k}''_{1}$ and $k'$, the restriction of this morphism to
$\Gamma_{1}$ is injective. But the extensions $\tilde{k}''$ over
$\tilde{k}''_{1}$ and $k'$ over $k$  are both totally ramified
and of degree $2e$. It follows that the extensions
$\tilde{k}''_{1}$ over $k$ and $\tilde{k}''$ over $k'$ have same
degree. All this implies that the map
$\alpha\mapsto\alpha_{\vert\tilde{k}''_{1}}$ induces an isomorphism
from $\Gamma_{1}$ onto $\Gamma_{2}$. The commutation relation
\eqref{eq20} follows from the fact that the Frobenius $\theta$
satisfies $\theta(x)=x^{q}$ for $x\in\tilde{k}''$ a root of unity of
order prime to $p$. Thus the first part of the proposition.

Then, the last part of the proposition follows from the fact that
$\Gamma_{1}/\Gamma_{0}=\{1,\tau\}$.
\end{dem}

\begin{co}\label{co4.1}
As a $\mathbb{Z}[\Gamma]$-module, $\mathbf{X}(T)$ is isomorphic to
$\mathbb{Z}[\Gamma']$ on which the action of $\Gamma'$ is the
natural one and $\theta$ acts by $\theta.\sigma^{l}=-\sigma^{ql}$.

The $\mathbb{Z}[\Gamma]$-submodules $M$ of $\mathbf{X}(T)$
which are such that $\mathbf{X}(T)/M$ is torsion free and are minimal for
this property are the $\mathbb{Z}[\Gamma']$-submodules $M$ of
$\mathbb{Z}[\Gamma']$ which are such that
$\mathbb{Z}[\Gamma']/M$ is torsion free and are minimal for this
property.
\end{co}
\begin{dem}
The first assertion follows from Proposition \ref{pr4.4}.

Let $M$ be a $\mathbb{Z}[\Gamma']$-submodule of
$\mathbb{Z}[\Gamma']$ which is such that $\mathbb{Z}[\Gamma']/M$
is torsion free and minimal for this property. It follows from Theorem
\ref{theoap} of the appendix that the eigenvalues of $\sigma$ in
$M\otimes_{\mathbb{Z}}\mathbb{C}$ are the $d$-th primitive roots of
unity for a certain divisor $d$ of $2e$, all with multiplicity
$1$. As $q$ is invertible in the ring $\mathbb{Z}/2e\mathbb{Z}$,
it follows that $M$ is stable under $\theta$. Thus the second assertion.
\end{dem}

As in the previous case, for $d$ a divisor (resp. $\underline{d}$ a
subset of the set divisors) of $2e$, we consider the
$\mathbb{Z}[\Gamma']$-submodule
(resp. $\mathbb{Z}[\Gamma']$-submodules) $M_{d}$
(resp. $M_{\underline{d}}$ and $\overline{M}_{\underline{d}}$) of
$\mathbb{Z}[\Gamma']$ introduced in the appendix.  Then we denote by
$S_{\underline{d}}$ the algebraic subgroup over $k$ of $T$ such that
$S_{\underline{d}}=\cap_{\chi\in
  \overline{M}_{\underline{d}}}\ker\chi$ and by
$\mathfrak{s}_{\underline{d}}$ is Lie algebra. We clearly have
$\mathfrak{s}_{\underline{d}}=\cap_{\chi\in
  M_{\underline{d}}}\ker\chi$ (here we consider elements of
$\mathbf{X}(T)$ as $k$-linear maps from $k'\nu$ into
$\tilde{k}''$). If $d$ is a divisor of $2e$ we simply put
$S_{d}=S_{\{d\}}$ and $\mathfrak{s}_{d}=\mathfrak{s}_{\{d\}}$.

\begin{co}\label{co4.2}
The $S_{d}$ for $d$ a divisor of $2e$ are the maximal proper subtori
of $T$ and the $S_{\underline{d}}$ for $\underline{d}$ a subset of the
set divisors of $2e$ are the proper subtori of $T$. The map
$\underline{d}\mapsto S_{\underline{d}}$ is a bijection from the set
of divisors of $2e$ onto the set of subtori of $T$ defined over
$k$. In particular, we have $S_{\emptyset}=T$ and
$S_{\underline{d}}=\{1\}$ if $\underline{d}$ is the set of divisors of
$2e$.

\end{co}
\begin{dem}
The result follows from Theorem \ref{theoap} and Corollary \ref{coap2}.
\end{dem}

With the notations of the appendix, for any divisor $d$ of $2e$, we have
\begin{equation}\label{eq21}
M_{d}=\oplus_{0\leq l\leq\varphi(d)}\mathbb{Z}P_{2e,d}(\sigma)\sigma^{l}
\end{equation}
where $\varphi$ is the Euler phi function and
$P_{2e,d}(X)=\frac{X^{2e}-1}{\Phi_{d}(X)}$, $\Phi_{d}$ being the
$d$-th cyclotomic polynomial.

If $d$ is a divisor of $2e$ we put $I_{d}=\{0\leq l\leq2e-1\vert
P_{2e,d}(\zeta^{l})=0\}$.

\begin{rem}\label{rem4.3}
  We have
  \begin{equation}\label{eq22}
I_{d}=\{0,\ldots,2e-1\}\smallsetminus\{\frac{2e}{d}l\vert
l\in(\mathbb{Z}/d\mathbb{Z})^{\times}\}.
\end{equation}
\end{rem}

If $\underline{d}$ is a set of divisors of $2e$, we put
$I_{\underline{d}}=\cap_{d\in\underline{d}}I_{d}$. We clearly have
$I_{d}=I_{\{d\}}$.

For $\underline{d}$ a set of divisors of $2e$, we denote by
$I'_{\underline{d}}$ the complementary subset of $I_{\underline{d}}$
into $\{0,\ldots,2e-1\}$.

\begin{pr}\label{pr4.5}
Let $\underline{d}$ be a subset of the set of divisors of $2e$. Then we have
\begin{equation}\label{eq23}
\mathfrak{s}_{\underline{d}}=\oplus_{l\in I_{\underline{d}}}k\varpi'^{l}\nu.
\end{equation}
\end{pr}
\begin{dem}
Let $\chi=P(\sigma)$, with $P(X)=\sum_{0\leq m\leq
  2e-1}a_{m}X^{m}\in\mathbb{Z}[X]$, be an element of
$\mathbb{Z}[\Gamma']=\mathbf{X}(T)$ viewed as a linear map from
$\mathfrak{t}$ into $\tilde{k}''$. Let $x=\sum_{0\leq l\leq
  2e-1}\alpha_{l}\varpi'^{l}\nu$ be an element of
$k'\nu=\mathfrak{t}$, where the $\alpha_{l}$ are in $k$.
Then one has
\begin{equation*}
\chi(x)=\sum_{0\leq l\leq 2e-1}\alpha_{l}P(\zeta^{l})\varpi'^{l}\nu.
\end{equation*}
As the numbers $\alpha_{l}P(\zeta^{l})$ belong to $k(\zeta)=\tilde{k}'_{1}$,
the valuation of the non-zero terms appearing in the sum above are all
distinct modulo $2e$. Thus we see that $\chi(x)=0$ if and only if
$\alpha_{l}P(\zeta^{l})=0$ for $0\leq l\leq2e-1$.  But, it follows
from \eqref{eq21} that
$\mathfrak{s}_{\underline{d}}=\cap_{d\in\underline{d}}\ker
P_{2e,d}(\sigma)$. All this implies the proposition.
\end{dem}

\begin{pr}\label{pr4.6}
Let $\underline{d}$ be a set of divisors of $2e$.
  
a) Suppose that $v''(u)=0$. Then the following assertions are equivalent

(i) $S_{\underline{d}}$ is admissible,

(ii) $\oplus_{0\leq l\leq e-1} k\varpi'^{2l}\nu\subset\mathfrak{s}_{\underline{d}}$,

(iii) all the numbers $\frac{2e}{d}$, $d\in\underline{d}$ are odd,

(iv) $I'_{\underline{d}}\subset 2\mathbb{N}+1$.

b) Suppose that $v''(u)=1$. Then the following assertions are equivalent

(i) $S_{\underline{d}}$ is admissible,

(ii) $\oplus_{0\leq l\leq e-1} k\varpi'^{2l+1}\nu\subset\mathfrak{s}_{\underline{d}}$,

(iii) all the numbers $\frac{2e}{d}$, $d\in\underline{d}$ are even.

(iv) $I'_{\underline{d}}\subset 2\mathbb{N}$.
\end{pr}
\begin{dem}
As $k''$ is unramified over $k'$, we have
$\Norm_{k''/k'}(k''^{\times})=\{x\in k'\vert
v'(x)\equiv0\pmod{2}\}$. Then according to Lemma \ref{lem4.1},
$S_{\underline{d}}$ is admissible if and only if
$(\mathfrak{s}_{\underline{d}}/u)^{\perp}\subset\oplus_{0\leq l\leq
  e-1}k\varpi'^{2l+1}$.

One checks that for $l\in\mathbb{Z}$
\begin{equation*}%\label{eq24}
  \sideset{}{_{k'/k}}\tr\varpi'^{l}=
  \sideset{}{_{m=0}^{2e-1}}\sum\sigma^{m}(\varpi'^{l})=
  \begin{cases}
  2e\varpi^{\frac{l}{2e}}\mbox{ if }l\in2e\mathbb{Z},\\
  0\mbox{ if }l\notin2e\mathbb{Z}.
\end{cases}
\end{equation*}
As the sets $I_{\underline{d}}$ and $I'_{\underline{d}}$ are invariant
modulo $2e$ under the transformation $l\mapsto 2e-l$, it follows that
\begin{equation}
  (\mathfrak{s}_{\underline{d}}/u)^{\perp}=
\begin{cases}\label{eq24}
  \oplus_{l\in I'_{\underline{d}}}k\varpi'^{l}\mbox{ if }v''(u)=0,\\
  \oplus_{l\in I'_{\underline{d}}}k\varpi'^{l+1}\mbox{ if }v''(u)=1.
\end{cases}
\end{equation}
As $I'_{\underline{d}}=\cup_{d\in\underline{d}}\{\frac{2e}{d}l\vert
l\in(\mathbb{Z}/d\mathbb{Z})^{\times}\}$, the proposition follows.
\end{dem}

\section{Multiplicity of characters and geometry of the momentum map}

In this section, in some particular cases, we compute the multiplicities of
the characters appearing in the restriction of the Weil representation
to an admissible anisotropic torus and establish a relation between
these multiplicities and some geometric objects related to the
momentum map. Our philosophy is largely inspired by what is known in
the real case.

More precisely we are interested in the restriction of the metaplectic
representation to an admissible subtorus of a maximal irreducible
torus $T\subset Sp(W)$ such that, using the notations introduced in
paragraphs \ref{2.1}, \ref{2.2} and \ref{3.2}, $k''$ is ramified over
$k'$ and $k'$ is unramified over $k$ or $k''$ is unramified over $k'$
and $k'$ is tamely and totally ramified over $k$. According to
Propositions \ref{pr4.1} and \ref{pr4.2}, in the first case $k'$ has
even degree over $k$, say $2f$ and we may suppose that
$\varpi''^{2}=\varpi'=\varpi$, while in the second case $k'$ has
degree $2e$ over $k$ with $e$ prime to $p$.

\subsection{}\label{5.1}
Along this paragraph and until paragraph \ref{5.4}, we suppose that
$k''$ is ramified over $k'$, $k'$ is unramified of degree $2f$ over
$k$ and $\varpi''^{2}=\varpi'=\varpi$. Also $S$ denotes a subtorus of $T$.

We denote by $\Gamma'$ the Galois group of $k'$ over $k$ which is
cyclic of order $2f$ and by $\tau$ the non-trivial element of the
Galois group $\Gamma$ of $k''$ over $k'$. In this case, according to
Lemma \ref{lem4.2}, we have $\Gamma=\{1,\tau\}\times\Gamma'$.

%\subsection{}\label{5.2}

We keep the notations of the preceding paragraph and of paragraph
\ref{2.2}, and study more closely the structure of $T$ and $S$.

\begin{lem}\label{lem5.2}
For all $\eta\in\mathcal{O}'\varpi''$, there exists a unique element
$\rho(\eta)\in T_{1}$ such that $\rho(\eta)-\eta\in k'$. The map
$\rho$ so defined satisfies $\rho(-\eta)=\rho(\eta)^{-1}$ and, for
$j\in\mathbb{N}$, $\rho(\eta)-\eta-1\in\varpi^{2j+1}\mathcal{O}'$ if
$\eta\in\varpi^{j}\mathcal{O}'\varpi''$.

Moreover, for $j\in\mathbb{N}$, the map $\rho$ induces a bijection
from $\varpi^{j}\mathcal{O}'\varpi''$ onto $T_{2j+1}$ and also a
group isomorphism $\rho_{j}$ from
$\varpi^{j}\mathcal{O}'\varpi''/\varpi^{2j+1}\mathcal{O}'\varpi''$ onto
$T_{2j+1}/T_{4j+3}$.
\end{lem}
\begin{dem}
Let $j$ be a natural integer. Recall that $T_{2j+1}=\{g\in T\vert
g-1\in\varpi''^{2j+1}\mathcal{O}''=\varpi^{j}\mathcal{O}''\varpi''\}$.
Thus, if $g\in T_{2j+1}$, there exists
$\eta\in\varpi^{j}\mathcal{O}'\varpi''$ and $\xi\in k'$ such
that $g=1+\eta+\xi$, and moreover,
\begin{equation*}
  (1+\xi)^{2}-\varpi^{2j+1}(\frac{\eta}{\varpi^{j}\varpi''})^{2}=1.
\end{equation*}
This equation forces $\xi$ to be in $\mathcal{O}'$.  But it follows
from Hensel's lemma that there exists a unique $\xi\in
\varpi\mathcal{O}'$ which satisfies it and moreover, we see that
$\xi\in\varpi^{2j+1}\mathcal{O}'$. The first assertion and the fact
that $\rho$ induces a bijection from $\varpi^{j}\mathcal{O}'\varpi''$
onto $T_{2j+1}$ follow easily.

  As for $\eta\in\varpi^{j}\mathcal{O}'\varpi''$ one has
  $\rho(\eta)\in1+\eta+\varpi^{2j+1}\mathcal{O}''$, it
  follows that for $\eta,\eta'\in\varpi^{j}\mathcal{O}'\varpi''$,
  $\rho(\eta)\rho(\eta')\in
  1+\eta+\eta'+\varpi^{2j+1}\mathcal{O}''$. We then have
$\rho(\eta)\rho(\eta')\rho(-(\eta+\eta'))\in1+\varpi^{2j+1}\mathcal{O}''$
  The last
  assertion follows easily.
\end{dem}

For $j\in\mathbb{N}$ we put $S_{j}=S\cap T_{j}$ and call it the $j$-th
congruence subgroup of $S$.

\begin{lem}\label{lem5.3}
  The restriction
    of $\rho_{j}$ to $\mathfrak{s}\cap\varpi^{j}\mathcal{O}'\varpi''
    /\mathfrak{s}\cap\varpi^{2j+1}\mathcal{O}'\varpi''$ is a group
    isomorphism onto $S_{2j+1}/S_{4j+3}$. We have $S_{2j}=S_{2j+1}$,
    $j\neq0$.

   The sequence of congruence subgroups $(S_{2j+1})_{j\in\mathbb{N}}$ is
    strictly decreasing and form a basis of neighbourhood of the
    identity in $S$. Moreover, one has $S=\{\pm1\}\cap S\times S_{1}$.  
\end{lem}
\begin{dem}
Recall that $k''$ (resp. $k'$) is a Galois extension of $k$ with
Galois group $\Gamma$ (resp. $\Gamma'$). It follows from paragraph
\ref{4.2} that $\mathbf{X}(T)$ is the set of restrictions to $T$ of
elements of $\mathbb{Z}[\Gamma]$. But, if $\sigma\in\Gamma$, we have
$\sigma(\mathcal{O}')=\mathcal{O}'$. Moreover, as $\Gamma$ is the
direct product of $\Gamma'$ and $\{1,\tau\}$, it follows that
$\sigma(T)=T$ and also that $\sigma(\mathfrak{t})=\mathfrak{t}$ for
every $\sigma\in\Gamma$.

Now, let $\sigma\in\Gamma$. If $\eta\in\mathcal{O}'\varpi''$, applying
$\sigma$ to the relation $\rho(\eta)-\eta-1\in k'$ we get that
$\sigma(\rho(\eta))=\rho(\sigma(\eta))$. From this follows that
$\sigma(T_{1})=T_{1}$. As every element of $\mathbf{X}(T)$ is the
restriction to $T$ of a product of a finite number of elements of
$\Gamma$, it follows that for $\chi\in\mathbf{X}(T)$,
$\chi(T_{1})\subset T_{1}$.

Let $\sigma_{1},\sigma_{2}\in\Gamma$ and $j\in\mathbb{N}$. It then
follows from Lemma \ref{lem5.2} that, for
$\eta\in\varpi^{j}\mathcal{O}'\varpi''$, we have
$\sigma_{1}(\rho(\eta))\sigma_{2}(\rho(\eta))=\rho(\sigma_{1}(\eta)+
\sigma_{2}(\eta))$ up to multiplication by an element of
$T_{4j+3}$. Under the same conditions, it follows readily by induction
that if $\sigma_{1},\ldots,$ $\sigma_{m}\in\Gamma$ we have
$\sigma_{1}(\rho(\eta))\cdots\sigma_{m}(\rho(\eta))=\rho(\sigma_{1}(\eta)+
\cdots+\sigma_{m}(\eta))$ up to multiplication by an element of
 $T_{4j+3}$. Put another way, if
$\chi\in\mathbf{X}(T)$ is a rational character of $T$, we have
$\chi(\rho(\eta))=\rho(\diff\chi(\eta))$ up to multiplication by an
element of $T_{4j+3}$, where $\diff\chi$ stands for
the differential of $\chi$ at $1$.

This last result shows that, under the same conditions, $\rho_{j}$
induces a group isomorphism from
$\ker\diff\chi\cap\varpi^{j}\mathcal{O}'\varpi''/
\ker\diff\chi\cap\varpi^{2j+1}\mathcal{O}'\varpi''$ onto
$\ker\chi\cap T_{2j+1}/\ker\chi\cap T_{4j+3}$.

From this result and the fact that $S$ is the intersection of the
kernel of a family of rational characters of $T$, follows that, under
the same conditions, the restriction of $\rho_{j}$ to
$\mathfrak{s}\cap\varpi^{j}\mathcal{O}'\varpi''
/\mathfrak{s}\cap\varpi^{2j}\mathcal{O}'\varpi''$ is a group
isomorphism onto $S_{2j+1}/S_{4j+3}$.

As $\mathfrak{s}$ is a non-zero subspace of $k''$ over $k$, it is
clear that for $j\in\mathbb{N}$,
$\mathfrak{s}\cap\varpi^{j+1}\mathcal{O}'\varpi''
\subsetneq\mathfrak{s}\cap\varpi^{j}\mathcal{O}'\varpi''$.

The last assertions of the lemma follow easily.
\end{dem}

\subsection{}\label{5.3}
 In order to study the restriction of the unitary characters of $T$ to
 $S$, we have to extend the orthogonal decomposition in Lemma
 \ref{lem4.3} to $\mathcal{O}'$.

We begin by the following result.

Let $n>0$ be an integer prime to $p$. Then the roots of $X^{n}-1$ in
an algebraic closure of $k$ are all simple. Thus if $\mathcal{I}$ is
the set of irreducible polynomials in $k[X]$ which divide $X^{n}-1$, one has 
$X^{n}-1=\prod_{P\in\mathcal{I}}P$.

Let $k'$ be an unramified extension of degree $n$ of $k$ and let
$\theta$ be the Frobenius of $k'$ over $k$ or the one of
$\mathbb{F}_{q'}$ over $\mathbb{F}_{q}$ (of course, here
$q'=q^{n}$). We equip $k'$ (resp $\mathbb{F}_{q'}$) with the
non-degenerate $k$-bilinear form $\langle x,y\rangle_{k'/k}=\tr_{k'/k}
xy$, $x,y\in k'$ (resp. $\mathbb{F}_{q}$-bilinear form $\langle
x,y\rangle_{\mathbb{F}_{q'}/\mathbb{F}_{q}}=\tr_{\mathbb{F}_{q'}/\mathbb{F}_{q}}
xy$, $x,y\in \mathbb{F}_{q'})$.

As $k$ is the quotient field of the factorial ring $\mathcal{O}$, we
have $\mathcal{I}\subset\mathcal{O}[X]$. If $P\in\mathcal{O}[X]$ let
us denote by $\overline{P}$ the element of $\mathbb{F}_{q}[X]$
obtained by reduction modulo $\mathcal{P}$ of the coefficients of
$P$. Then we clearly have
$X^{n}-1=\prod_{P\in\mathcal{I}}\overline{P}$ in
$\mathbb{F}_{q}[X]$. As $n$ is prime to $p$, the roots of $X^{n}-1$ in
an algebraic closure of $\mathbb{F}_{q}$ are all simple. It follows
that the $\overline{P}$, $P\in\mathcal{I}$ are all distinct and, by
Hensel's lemma, that the above factorisation is the decomposition into
irreducible factors of $X^{n}-1$ in $\mathbb{F}_{q}[X]$. In
particular, for $P\in\mathcal{I}$, $P$ and $\overline{P}$ have same
degree.

For $P\in\mathcal{I}$ let $W_{P}=\ker P(\theta)$
(resp. $\overline{W}_{P}=\ker \overline{P}(\theta)$), which is a
$k$-linear subspace of $k'$ (resp. $\mathbb{F}_{q}$-linear subspace of
$\mathbb{F}_{q'}$).

\begin{lem}\label{lem5.4}
The subspaces $W_{P}$ (resp.$\overline{W}_{P}$), $P\in\mathcal{I}$ are
invariant under $\theta$, minimal for this property and their
dimension is the degree of $P$. We have the decomposition as vector-spaces
\begin{equation}\label{eq25}
k'=\oplus_{P\in\mathcal{I}}W_{P}\mbox{,
}\mathbb{F}_{q'}=\oplus_{P\in\mathcal{I}}\overline{W}_{P},
\end{equation}
and also as $\mathcal{O}$-modules
\begin{equation}\label{eq26}
\mathcal{O}'=\oplus_{P\in\mathcal{I}}\mathcal{O'}\cap W_{P}.
\end{equation}
Let $W\subset k'$ be a $\theta$-invariant $k$-linear subspace such
that the restriction to $W$ of $\langle\, ,\,\rangle_{k'/k}$ is
non-degenerate. Then
$\mathcal{O}'=\mathcal{O}'\cap W\oplus\mathcal{O}'\cap W^{\perp}$.
\end{lem}
\begin{dem}
The first assertion and relation \eqref{eq25} follow easily from the
decomposition into irreducibles of $X^{n}-1$ which is both the minimal and
the characteristic polynomial of the $k$-endomorphism or
$\mathbb{F}_{q}$-endomorphism $\theta$.

Let $P\in\mathcal{I}$. As $\mathcal{O}'$ is an open neighbourhood of
$0$ in $k'$, we have $\mathcal{O}'\cap W_{P}\neq0$ and thus,
$\mathcal{O}'^{\times}\cap W_{P}\neq\emptyset$. It follows that the
image of $\mathcal{O}'\cap W_{P}$ under the natural projection
$\mathsf{p}:\mathcal{O}'\rightarrow \mathbb{F}_{q'}$ is a non-zero
$\theta$-invariant subspace of $\overline{W}_{P}$ and thus, that
$\mathsf{p}(\mathcal{O}'\cap W_{P})=\overline{W}_{P}$. In particular,
there exists a subset $b_{P}$ of $\mathcal{O}'^{\times}\cap W_{P}$
such that $\mathsf{p}(b_{P})$ is a basis of $\overline{W}_{P}$ over
$\mathbb{F}_{q}$.

Now let $b=\cup_{P\in\mathcal{I}}b_{P}$. Then $\mathsf{p}(b)$ is a
basis of $\mathbb{F}_{q'}$ over $\mathbb{F}_{q}$. It follows that $b$
is a basis of $\mathcal{O}'$ over $\mathcal{O}$. Thus the relation
\eqref{eq26}.

Let $W$ be as in the last assertion of the lemma. Then $W^{\perp}$ is
also a $\theta$-invariant subspace. It follows easily that both $W$
and $W^{\perp}$ are the direct sum of the subspaces $W_{P}$,
$P\in\mathcal{I}$ they contain. The lemma follows easily.
\end{dem}

\begin{pr}\label{pr5.1}
  We suppose that $f$ and $p$ are coprime. Then we have
\begin{equation*}
  \mathcal{O}'=\mathcal{O}'\cap(\mathfrak{s}/\varpi'')
  \oplus\mathcal{O}'\cap(\mathfrak{s}/\varpi'')^{\perp}.
\end{equation*}
\end{pr}
\begin{dem}
It is clear that $\mathfrak{s}/\varpi''$ is a $\Gamma'$-invariant
$k$-linear subspace of $k'$. The proposition then follows from
Lemmas \ref{lem4.3} and \ref{lem5.4}.
\end{dem}

Let us put %$\varepsilon_{S}=\vert S/S_{1}\vert$. Then,
$\varepsilon_{S}=2$ if $S\cap\{\pm1\}=\{1\}$ and $1$ otherwise.

\begin{rem}\label{rem5.1}
If
$\chi=\sum_{\sigma\in\Gamma'}a_{\sigma}\sigma\in\mathbb{Z}[\Gamma']$
is a character of $T$, let us put
$a(\chi)=\sum_{\sigma\in\Gamma'}a_{\sigma}$. Then one has
$\chi(-1)=(-1)^{a(\chi)}$. Let $\theta$ be a generator of
$\Gamma'$. Let $d$ be a divisor of $2f$. Let us put
$\chi_{d}=\sum_{0\leq l\leq\frac{2f}{d}-1}\theta^{ld}Q_{d}(\theta)$
(see the appendix for the definition of $Q_{d}$). Then it follows from
Theorem \ref{theoap} of the appendix that
$(\chi_{d}\theta^{l})_{0\leq l\leq\varphi(d)-1}$ is a basis of
$M_{d}$. We have $a(\chi_{d}\theta^{l})=\frac{2f}{d}Q_{d}(1)$, $0\leq
l\leq\varphi(d)-1$. But $Q_{d}(1)=0$ if $d>1$ and $Q_{1}(1)=1$. Thus,
in any case we have $\chi(-1)=1$, $\chi\in M_{d}$, showing that $-1\in
S_{d}$ and $\varepsilon_{S_{d}}=1$.

Nevertheless, it follows from the results of Lemma \ref{lem5.14} below that
it happens in many cases that $\varepsilon_{S}=2$ if $S$ is a non
maximal subtorus of $T$.

For instance, if $f$ is odd and $\underline{d}=\{1,2\}$, we have
$\varepsilon_{S_{\underline{d}}}=2$.
\end{rem}

\subsection{}\label{5.4}
From now on we suppose that $f$ and $p$ are coprime.  As for $T$, we
call conductor of a unitary character $\chi$ of $S$ the smallest
integer $j$ such that $S_{j}$ is contained into $\ker\chi$.

In this paragraph, keeping the notations of the preceding one, we
suppose that $S$ is admissible. Then, we
determine the characters of $S$ appearing in the restriction of the
metaplectic representation and compute their multiplicities therein.

\begin{lem}\label{lem5.6}
Le $\chi$ be a unitary character of $T$ appearing in the restriction
to $T$ of the metaplectic representation. Then $\chi_{\vert S}$ has
the same conductor as $\chi$.
\end{lem}
\begin{dem}
Let $\chi\in\hat{T}$ appearing in the restriction to $T$ of the
metaplectic representation. Clearly, we may suppose that its conductor
is strictly positive. Then according to Theorem \ref{theo2.1} and
Proposition \ref{pr2.1}, its conductor is even, say $2j$, and there
exists $a\in\mathcal{O}''^{\times}$ such that
\begin{equation}\label{eq28}
  \chi(\rho(\eta))\!
  =\!\psi(-\frac{(-1)^{\mu-j}}{4}\sideset{}{_{k''/k}}\tr\varpi^{\mu-j}
  \sideset{}{_{k''/k'}}\Norm(a)\varpi''\eta)\mbox{,
  }\eta\in\varpi^{[\frac{j}{2}]}\mathcal{O}'\varpi''.
\end{equation}
Suppose that $\chi_{\vert S}$ has conductor less or equal to
$2j-1$. Then according to Lemma \ref{lem5.3}, we have
$\chi(\rho(\varpi^{j-1}\eta\varpi''))=1$ for all
$\eta\in(\mathfrak{s}/\varpi'')\cap\mathcal{O}'$.  But, taking in
account formula \eqref{eq28} and Lemma \ref{lem5.3}, this condition is
equivalent to
\begin{equation*}
  \sideset{}{_{k'/k}}\tr(\sideset{}{_{k''/k'}}\Norm(a)\eta)
  \in\varpi\mathcal{O}\mbox{,
  }\eta\in(\mathfrak{s}/\varpi'')\cap\mathcal{O}'.
\end{equation*}
Then writing $\Norm_{k''/k'}(a)=\alpha+\beta$ with
$\alpha\in\mathfrak{s}/\varpi''$
and $\beta\in(\mathfrak{s}/\varpi'')^{\perp}$, we get
\begin{equation}\label{eq29}
  \sideset{}{_{k'/k}}\tr(\alpha\eta)
  \in\varpi\mathcal{O}\mbox{,
  }\eta\in(\mathfrak{s}/\varpi'')\cap\mathcal{O}'.
\end{equation}
But it follows from Proposition \ref{pr5.1} that the relation
\eqref{eq29} is still valid for any $\eta\in\mathcal{O}'$.  We deduce
from this that
$\alpha\in(\mathfrak{s}/\varpi'')\cap\varpi\mathcal{O}'$ and thus,
$\beta=\Norm_{k''/k'}(a)-\alpha
\in(\mathfrak{s}/\varpi'')^{\perp}\cap\Norm_{k''/k'}(\mathcal{O}''^{\times})$, 
in contradiction with the hypothesis that $S$ is admissible.
\end{dem}

If $\alpha\in(\mathbb{F}_{q'}^{\times})^{2}$, we denote by
$\Theta(\alpha)$ the set of
$s\in\mathsf{p}(\mathcal{O}'\cap(\mathfrak{s}/\varpi'')^{\perp})$ such
that $\alpha+s \in(\mathbb{F}_{q'}^{\times})^{2}$. Then, as it
contains $0$, $\Theta(\alpha)$ is a non-empty finite set. We also
recall that
$\Norm_{k''/k'}(\mathcal{O''}^{\times})=(\mathcal{O'}^{\times})^{2}$,
thus if $a\in\mathcal{O''}^{\times}$ and
$\alpha=\mathsf{p}(\Norm_{k''/k'}(a))$, then $\Theta(\alpha)$ is
defined.

\begin{pr}\label{pr5.2}
(i) A unitary character of $S$ appears in the restriction to $S$ of
  the metaplectic representation if and only if it is trivial or its
  conductor is $2j$ with $j>0$ and there exists
  $a\in\mathcal{O}''^{\times}$ such that %its restriction to $S_{j}$
  %viewed as a character of $S_{j}/S_{2j}$ satisfies the relation
\begin{equation}\label{eq30}
  \chi_{\vert S_{j}}(\rho_{j}(\eta))=
  \psi(-\frac{(-1)^{\mu-j}}{4}\sideset{}{_{k''/k}}\tr\varpi^{\mu-j}
  \sideset{}{_{k''/k'}}\Norm(a)\varpi''\eta)\mbox{, }\!\eta\in
  \varpi^{[\frac{j}{2}]}\mathcal{O}'\varpi''\cap\mathfrak{s}.
\end{equation}

(ii) The trivial character of $S$ appears in the
restriction of the metaplectic representation with multiplicity $1$.

(iii) Let $\chi$ be a unitary character of $S$ satisfying the relation
\eqref{eq30}. Then $\chi$ appears in the restriction of the
metaplectic representation with multiplicity
\begin{equation}\label{eq31}
  m(\chi)=
  \varepsilon_{S}\vert\Theta(\mathsf{p}(\sideset{}{_{k''/k'}}\Norm(a)))\vert
  q^{\codim(\mathfrak{s})(j-1)}.
\end{equation}
\end{pr}
\begin{dem}
Assertions (i) and (ii) follow readily from Lemma \ref{lem5.6} and
Theorem \ref{theo2.1}.

Now let $\chi$ be a unitary character of $S$ satisfying the relation
\eqref{eq30}. Then $m(\chi)$ is the number of unitary characters of
$T$ with the conductor $2j$ appearing in the restriction to $T$ of the
metaplectic representation and having $\chi$ as restriction to $S$.

First we will compute the number of unitary characters of $T_{j}$ with
the conductor $2j$ appearing in the restriction to $T_{j}$ of the
metaplectic representation and having same restriction as $\chi$ to
$S_{j}$. Such a character is of the form
$\chi_{\Norm_{k''/k'}(b),j,2j}$ with
$b\in\mathcal{O}''^{\times}$. Thus, we must have
$(\chi_{\Norm_{k''/k'}(b),j,2j})_{\vert S_{j}}(\chi_{\vert S_{j}})^{-1} =1$ that is,
for any $\eta\in\varpi^{[\frac{j}{2}]}\mathcal{O}'\varpi''\cap\s$
\begin{equation*}
\psi(-\frac{(-1)^{\mu-j}}{4}\sideset{}{_{k''/k}}\tr\varpi^{\mu-j}
(\sideset{}{_{k''/k'}}\Norm(b)-\sideset{}{_{k''/k'}}\Norm(a))\varpi''\eta)=1.
%\mbox{,
%}\eta\in\varpi^{[\frac{j}{2}]}\mathcal{O}'\varpi''\cap\s
\end{equation*}
Denoting by $p_{S}$ the orthogonal projection from
$k'$ onto $\mathfrak{s}/\varpi''$, this last relation may be rewritten as
\begin{equation*}
  v(p_{S}(\sideset{}{_{k''/k'}}\Norm(b)-\sideset{}{_{k''/k'}}\Norm(a)))
  \geq j-[\frac{j}{2}].
\end{equation*}
It follows from Proposition \ref{pr5.1} that this is equivalent to
\begin{equation*}
  \sideset{}{_{k''/k'}}\Norm(b)-\sideset{}{_{k''/k'}}\Norm(a)\in
  \mathcal{O}'\cap(\mathfrak{s}/\varpi'')^{\perp}+
  \varpi^{j-[\frac{j}{2}]}\mathcal{O}'\cap(\mathfrak{s}/\varpi'').
\end{equation*}
Applying the projection $\mathsf{p}$ to this last relation we obtain
that
\begin{equation*}
  \mathsf{p}(\sideset{}{_{k''/k'}}\Norm(b))-
 \mathsf{p}(\sideset{}{_{k''/k'}}\Norm(a))  
   \in\Theta(\mathsf{p}(\sideset{}{_{k''/k'}}\Norm(a))).
\end{equation*}
It follows that the character $\chi_{\Norm_{k''/k'}(b),j,2j}$ of
$T_{j}$ has the same restriction to $S_{j}$ as $\chi$ if and only if there
exists $\alpha\in\Theta(\mathsf{p}(\sideset{}{_{k''/k'}}\Norm(a)))$ %and
%$x\in\mathsf{p}^{-1}(\alpha)\cap(\mathfrak{s}/\varpi'')^{\perp}$
such that
\begin{equation}\label{eq32}
  \sideset{}{_{k''/k'}}\Norm(b)\in\sideset{}{_{k''/k'}}\Norm(a)+
  \mathsf{p}^{-1}(\alpha)\cap(\mathfrak{s}/\varpi'')^{\perp}+
  \varpi^{j-[\frac{j}{2}]}\mathcal{O}'\cap(\mathfrak{s}/\varpi'').
\end{equation}
Now let $\alpha\in\Theta(\mathsf{p}(\Norm_{k''/k'}(a)))$. Due to the
  definition of $\Theta(\mathsf{p}(\Norm_{k''/k'}(a)))$ and Hensel's
  lemma, there exists $y\in\mathsf{p}^{-1}(\alpha)$ such that
  $\Norm_{k''/k'}(a)+y\in(\mathcal{O}'^{\times})^{2}=
  \Norm_{k''/k'}(\mathcal{O}''^{\times})$. Then one has
\begin{eqnarray}\label{eq33}
  \sideset{}{_{k''/k'}}\Norm(a)\!+\!
  \mathsf{p}^{-1}(\alpha)\cap(\mathfrak{s}/\varpi'')^{\perp}
  \!+\!\varpi\mathcal{O}'&=&\sideset{}{_{k''/k'}}\Norm(a)\!+\!y\!+
  \!\varpi\mathcal{O}'\\
&\subset&\sideset{}{_{k''/k'}}\Norm(\mathcal{O}''^{\times}).\notag
\end{eqnarray}
It follows that the set appearing in the right side of \eqref{eq32}, which
depends only on $\alpha$, is contained in
$\Norm_{k''/k'}(\mathcal{O}''^{\times})$. Moreover, if
$b,b'\in\mathcal{O}'^{\times}$ we have
$\chi_{\Norm_{k''/k'}(b),j,2j}=\chi_{\Norm_{k''/k'}(b'),j,2j}$ if and
only if $\Norm_{k''/k'}(b)-\Norm_{k''/k'}(b')
\in\varpi^{j-[\frac{j}{2}]}\mathcal{O}'$.

  But one has $\mathsf{p}^{-1}(\alpha)\cap(\mathfrak{s}/\varpi'')^{\perp}=
  x+\varpi\mathcal{O}'\cap(\mathfrak{s}/\varpi'')^{\perp}$ for any $x\in
 \mathsf{p}^{-1}(\alpha)\cap(\mathfrak{s}/\varpi'')^{\perp}$, while 
\begin{multline*}
%\begin{split}
\big(\varpi\mathcal{O}'\cap(\mathfrak{s}/\varpi'')^{\perp}+
\varpi^{j-[\frac{j}{2}]}\mathcal{O}'\cap(\mathfrak{s}/\varpi'')\big)/
\varpi^{j-[\frac{j}{2}]}\mathcal{O}'\\
=\varpi\mathcal{O}'\cap(\mathfrak{s}/\varpi'')^{\perp}/
\varpi^{j-[\frac{j}{2}]}\mathcal{O}'\cap(\mathfrak{s}/\varpi'')^{\perp}.
%\end{split}
\end{multline*}
This last set is finite and has cardinal
$q^{\codim(\mathfrak{s})(j-[\frac{j}{2}]-1)}$. It follows that the
number of unitary characters of $T_{j}$ with the conductor $2j$
appearing in the restriction to $T_{j}$ of the metaplectic
representation and having same restriction as $\chi$ to $S_{j}$ is
$\vert\Theta(\mathsf{p}(\Norm_{k''/k'}(a)))\vert
q^{\codim(\mathfrak{s})(j-[\frac{j}{2}]-1)}$.

Now, given a unitary character $\vartheta$ of $T_{j}$ with the
conductor $2j$ appearing in the restriction to $T_{j}$ of the
metaplectic representation and having same restriction as $\chi$ to
$S_{j}$ we will compute the number of unitary characters of $T$ having
restriction $\vartheta$ to $T_{j}$ and $\chi$ to $S$. Consider the
subgroup $ST_{j}$ of $T$. We define a unitary character
$\tilde{\vartheta}$ of $ST_{j}$ by deciding that
$\tilde{\vartheta}_{\vert T_{j}}=\vartheta$ and $\tilde{\vartheta}_{\vert
  S}=\chi$. The expected number is the number of unitary characters of
$T$ whose restriction to $ST_{j}$ is $\tilde{\vartheta}$, that is
$\vert T/ST_{j}\vert$. But it follows easily from Lemma \ref{lem5.3}
that $\vert T/ST_{j}\vert=\varepsilon_{S}
q^{\codim(\mathfrak{s})[\frac{j}{2}]}$. As $m(\chi)=\vert
T/ST_{j}\vert\vert\Theta(\mathsf{p}(\Norm_{k''/k'}(a)))\vert
q^{\codim(\mathfrak{s})(j-[\frac{j}{2}]-1)}$, formula \eqref{eq31}
follows.
\end{dem}

\subsection{}\label{5.5}
Along this paragraph and until paragraph \ref{5.8}, we suppose that
$k''$ is unramified over $k'$, $k'$ is tamely and totally ramified
over $k$ of degree $2e$ with $e$ prime to $p$. Also we use the
notations introduced in paragraph \ref{4.4}.

In the following we consider $S$ a subtorus of $T$ defined over $k$
with Lie algebra $\mathfrak{s}$ and $\underline{d}$ a set of divisors
of $2e$ such that $S=S_{\underline{d}}$. We also put
$M=M_{\underline{d}}$, $\overline{M}=\overline{M}_{\underline{d}}$,
$I=I_{\underline{d}}$  and $I'=I'_{\underline{d}}$.

\subsection{}\label{5.6}
In this paragraph we study more closely the structure of $T$ and
$S$. We use the notations of paragraph \ref{2.2}. The proof of the
following lemma, analogous to the one of Lemma \ref{lem5.2}, is left
to the reader.

\begin{lem}\label{lem5.7}
For all $\eta\in\varpi'\mathcal{O}'\nu$, there exists a unique element
$\rho(\eta)\in T_{1}$ such that $\rho(\eta)-\eta\in k'$. The map $\rho$
so defined satisfies $\rho(-\eta)=\rho(\eta)^{-1}$ and, for any
integer $j>0$, $\rho(\eta)-\eta-1\in\varpi'^{2j}\mathcal{O}''$ if
$\eta\in\varpi'^{j}\mathcal{O}'\nu$.

Moreover, for any integer $j>0$, the map $\rho$ induces a bijection
from $\varpi'^{j}\mathcal{O}'\nu$ onto $T_{j}$ and also a group
isomorphism $\rho_{j}$ from
$\varpi'^{j}\mathcal{O}'\nu/\varpi'^{2j}\mathcal{O}'\nu$ onto
$T_{j}/T_{2j}$.
\end{lem}

For $j\in\mathbb{N}$, we put $S_{j}=S\cap T_{j}$ and call it the
$j$-th congruence subgroup of $S$.

\begin{lem}\label{lem5.8}
  Let $j>0$ be an integer. Then we have
\begin{equation}\label{eq34}
  ST_{j}=\{x\in T\vert\chi(x)\in1+\varpi'^{j}\tilde{\mathcal{O}}''\mbox{, }
  \chi\in\overline{M}\}.
\end{equation}
\end{lem}
\begin{dem}
First of all, we clearly have the inclusion
\begin{equation*}
  ST_{j}\subset\{x\in
  T\vert\chi(x)\in1+\varpi'^{j}\tilde{\mathcal{O}}''\mbox{, }
  \chi\in\overline{M}\}.
\end{equation*}
We will prove the reverse inclusion.
  
Recall that $T=\mu_{q+1}\times T_{1}$. It is clear that if
$x\in\mu_{q+1}$ (resp. $y\in T_{1}$) and $\chi\in\mathbf{X}(T)$, then
$\chi(x)\in\mu_{q+1}$
(resp. $\chi(y)\in1+\varpi'\tilde{\mathcal{O}}''$). As
$\mu_{q+1}\cap(1+\varpi'\tilde{\mathcal{O}}'')=\{1\}$, it follows that
if $x\in\mu_{q+1}$, $y\in T_{1}$, and $\chi\in\mathbf{X}(T)$,
$\chi(xy)\in1+\varpi'\tilde{\mathcal{O}}''$ if and only if $\chi(x)=1$
and $\chi(y)\in1+\varpi'\tilde{\mathcal{O}}''$.

%Thus, if
%$x\in\mu_{q+1}$, $y\in T_{1}$, and $\chi\in\mathbf{X}(T)$, are such
%that $\chi(xy)\in1+\varpi'\tilde{\mathcal{O}}''$, then one has
%$\chi(x)\in\mu_{q+1}\cap(1+\varpi'\tilde{\mathcal{O}}'')=\{1\}$ and
%$\chi(y)\in1+\varpi'\tilde{\mathcal{O}}''$. It follows that, under the
%same conditions, $\chi(xy)\in1+\varpi'\tilde{\mathcal{O}}''$ if and
%only if $\chi(x)=1$ and $\chi(y)\in1+\varpi'\tilde{\mathcal{O}}''$.

We deduce that it suffices to prove that
\begin{equation*}
  \{x\in T_{1}\vert\chi(x)\in1+\varpi'^{j}\tilde{\mathcal{O}}''\mbox{, }
  \chi\in\overline{M}\}\subset S_{1}T_{j}.
\end{equation*}
Thus, let $x\in T_{1}$ such that
$\chi(x)\in1+\varpi'^{j}\tilde{\mathcal{O}}''$ for every $\chi\in\overline{M}$.

Let $N$ be the direct sum of the $M_{d}$ such that
$d\notin\underline{d}$. Then $N$ is clearly orthogonal to
$\overline{M}$ with respect to the bilinear form $\langle\,,\,\rangle$
of Corollary \ref{coap1}. As this bilinear form is positive definite,
it follows that $N\cap\overline{M}=\{0\}$. Thus $N$ identifies, under
the natural projection from $\mathbf{X}(T)$ onto
$\mathbf{X}(S)=\mathbf{X}(T)/\overline{M}$, to a
$\mathbb{Z}[\Gamma]$-submodule of $\mathbf{X}(S)$ having the same
rank.

Now let $\alpha_{x}:N\rightarrow\tilde{k''}$ the map such that
$\alpha_{x}(\eta)=\eta(x)$. We clearly have
$\alpha_{x}(\eta)\in1+\varpi'\tilde{\mathcal{O}}''$. Moreover $\alpha_{x}$ is
clearly $\Gamma$-equivariant.

Now let $\theta\in\mathbf{X}(S)$. Then there exists
$a\in\mathbb{Z}\smallsetminus\{0\}$ such that $a\theta\in N$ and,
according to Corollary \ref{coap3}, we may choose $a$ to be a divisor
of $2e$ and thus to be prime to $p$. Then there exists a unique
element $\theta(x)\in1+\varpi'\tilde{\mathcal{O}}''$ such that
$\theta(x)^{a}=(a\theta)(x)$. It is clear that $\theta(x)$ doesn't
depend upon the choice of $a$. Then we extend $\alpha_{x}$ to
$\mathbf{X}(S)$ putting $\alpha_{x}(\theta)=\theta(x)$. We easily
check that $\alpha_{x}\in
\Hom_{\Gamma}(\mathbf{X}(S),\tilde{k}'')$. As $\alpha_{x}$ takes it
values in $1+\varpi'\tilde{\mathcal{O}}''$, it follows that there exists $y\in
S_{1}$ such that $\alpha_{x}(\theta)=\theta(y)$,
$\theta\in\mathbf{X}(S)$.

Then $xy^{-1}\in T_{1}$ and satisfies $\chi(xy^{-1})=1$, $\chi\in N$
and thus $\chi(xy^{-1})\in1+\varpi'^{j}\tilde{\mathcal{O}}''$, $\chi\in
M+N$. Now, let $\chi\in\mathbf{X}(T)$. There exists
$a\in\mathbb{Z}\smallsetminus\{0\}$ prime to $p$ such that $a\chi\in
M+N$. Thus $\chi(xy^{-1})^{a}\in1+\varpi'^{j}\tilde{\mathcal{O}}''$, and
$\chi(xy^{-1})\in1+\varpi'^{j}\tilde{\mathcal{O}}''$ as $a$ is prime to
$p$. This implies that $xy^{-1}\in T_{j}$.
\end{dem}

\begin{lem}\label{lem5.9}
  Let $j$ and $m$ be integers such that $0<j<m\leq2j$. Then for
  $\eta\in\varpi'^{j}\mathcal{O}'$ we have
\begin{equation*}
\rho(\eta\nu)\in ST_{m} \Longleftrightarrow
\eta\nu\in(\mathfrak{s}\cap\varpi'^{j}\mathcal{O}'\nu)+\varpi'^{m}\mathcal{O}''.
\end{equation*}    
\end{lem}
\begin{dem}
Let $\eta\in\varpi'^{j}\mathcal{O}'$. According to Lemma
\ref{lem5.8}, the condition $\rho(\eta\nu)\in ST_{m}$ is equivalent
to $\chi(\rho(\eta))\in1+\varpi'^{m}\tilde{\mathcal{O}}''$,
$\chi\in\overline{M}$.

But, according to Corollary \ref{coap4}, the order of any element in
the torsion $\mathbb{Z}$-module $\overline{M}/M$ is prime to $p$. It
follows that the last condition just above is equivalent to
$\chi(\rho(\eta))\in1+\varpi'^{m}\tilde{\mathcal{O}}''$, $\chi\in M$. But as
$M_{d}=\oplus_{0\leq
  l\leq\varphi(d)-1}\mathbb{Z}P_{2e,d}(\sigma)\sigma^{l}$ for any
divisor $d$ of $2e$, this is also equivalent to
\begin{equation}\label{eq35}
P_{2e,d}(\sigma)(\rho(\eta))\in1+\varpi'^{m}\tilde{\mathcal{O}}''\mbox{, }
d\in\underline{d}.
\end{equation}
Now, let us write $j=2e[\frac{j}{2e}]+l$ with $0\leq l\leq 2e-1$ an
integer. Then we have
\begin{equation*}
  \varpi'^{j}\mathcal{O}'=
  \varpi^{[\frac{j}{2e}]}(\oplus_{0\leq i\leq l-1}\varpi\mathcal{O}\varpi'^{i}
  \oplus\oplus_{l\leq i\leq 2e-1}\mathcal{O}\varpi'^{i}).
\end{equation*}
Thus, using equations \eqref{eq23} and \eqref{eq24}, we deduce that on
one hand
\begin{equation}\label{eq36}
  \varpi'^{j}\mathcal{O}'=\mathfrak{s}/\nu\cap\varpi'^{j}\mathcal{O}'
  \oplus(\mathfrak{s}/\nu)^{\perp}\cap\varpi'^{j}\mathcal{O}', 
\end{equation}
and on the other hand, there exists $\epsilon'_{i}\in\{0,1\}$, $i\in
I'_{\underline{d}}$ such that
\begin{equation}\label{eq37}
  (\mathfrak{s}/\nu)^{\perp}\cap\varpi'^{j}\mathcal{O}'=
  \varpi^{[\frac{j}{2e}]}(\oplus_{i\in I'_{\underline{d}}}
  \varpi^{\epsilon'_{i}}\mathcal{O}\varpi'^{i}). 
\end{equation}
Let $d\in\underline{d}$ and write $P_{2e,d}(X)=\sum_{0\leq
  t\leq2e-\varphi(d)-1}a_{t}X^{t}$, where
$a_{t}\in\mathbb{Z}$. According to \eqref{eq36} we may write
$\eta=\eta_{1}+\eta_{2}$ with
$\eta_{1}\in\mathfrak{s}/\nu\cap\varpi'^{j}\mathcal{O}'$ and
$\eta_{2}\in(\mathfrak{s}/\nu)^{\perp}\cap\varpi'^{j}\mathcal{O}'$. Moreover,
according to Lemma \ref{lem5.7}, we have $\rho(\eta)=1+\eta\nu+\xi$
with $\xi\in\varpi'^{2j}\mathcal{O}'$. Then, one has
\begin{equation*}
  P_{2e,d}(\sigma)(\rho(\eta))=
  \prod_{t=0}^{2e-\varphi(d)-1}(\sigma^{t}(\rho(\eta))^{a_{t}}=
\prod_{t=0}^{2e-\varphi(d)-1}(1+\sigma^{t}(\eta\nu)+\sigma^{t}(\xi))^{a_{t}}.
\end{equation*}
But, for each $t$, one has
$\sigma^{t}(\eta\nu)\in\varpi'^{j}\tilde{\mathcal{O}}''$ and
$\sigma^{t}(\xi))\in\varpi'^{2j}\tilde{\mathcal{O}}''$. From this, we
deduce that
\begin{equation*}
  P_{2e,d}(\sigma)(\rho(\eta))\in
  1+P_{2e,d}(\sigma)(\eta\nu)+\varpi'^{2j}\tilde{\mathcal{O}}''.
\end{equation*}
Moreover, one clearly has
$P_{2e,d}(\sigma)(\eta\nu)=P_{2e,d}(\sigma)(\eta_{2}\nu)$. Thus the
condition \eqref{eq35} implies
$P_{2e,d}(\sigma)(\eta_{2}\nu)\in\varpi'^{m}\tilde{\mathcal{O}}''$,
$d\in\underline{d}$.

Now, according to \eqref{eq37}, we write $\eta_{2}=
\varpi^{[\frac{j}{2e}]}\sum_{i\in
  I'_{\underline{d}}}\alpha_{i}\varpi^{\epsilon'_{i}}\varpi'^{i}$,
with $\alpha_{i}\in\mathcal{O}$. Then, we have
\begin{equation}\label{eq38}
P_{2e,d}(\sigma)(\eta_{2}\nu)=\varpi^{[\frac{j}{2e}]}\sum_{i\in
  I'_{\underline{d}}}\alpha_{i}\varpi^{\epsilon'_{i}}P_{2e,d}(\zeta^{i})\varpi'^{i}
\in\varpi'^{m}\tilde{\mathcal{O}}''
\mbox{, }d\in\underline{d}.
\end{equation}
The roots of the polynomials $P_{2e,d}(X)$ all appear with
multiplicity one and moreover are $2e$-th roots of unity. But if
$\alpha\neq\beta$ are two $2e$-th roots of unity in $\tilde{k}''$, one
has $\vert\alpha-\beta\vert=1$. It follows that for $d\in\underline{d}$ and
$i\in I'_{\underline{d}}$, one has
$P_{2e,d}(\zeta^{i})\in\tilde{\mathcal{O}}''^{\times}$ or $P_{2e,d}(\zeta^{i})=0$.

Now, let $i\in I'_{\underline{d}}$. Then, there exists
$d\in\underline{d}$ such that $P_{2e,d}(\zeta^{i})\neq0$. We choose
such a $d$. It follows from the above considerations, that the
non-zero terms appearing inside the sum in the right-hand side of
\eqref{eq38} have distinct valuations modulo $2e$.Thus we have
$\varpi^{[\frac{j}{2e}]}\alpha_{i}\varpi^{\epsilon'_{i}}\varpi'^{i}
\in\varpi'^{m}\tilde{\mathcal{O}}''$. It follows that
$\eta_{2}\in\varpi'^{m}\mathcal{O}'$ and that $\eta\nu\in
(\mathfrak{s}\cap\varpi'^{j}\mathcal{O}'\nu)+\varpi'^{m}\mathcal{O}''$,
showing the direct implication of the lemma.

Conversely, let
$\eta\nu\in
(\mathfrak{s}\cap\varpi'^{j}\mathcal{O}'\nu)+\varpi'^{m}\mathcal{O}''$. Then
for all $\chi\in\overline{M}$, we have
$\chi(\rho(\eta))\in1+\varpi'^{m}\tilde{\mathcal{O}}''$. Then it
follows from Lemma \ref{lem5.8} that $\rho(\eta)\in ST_{m}$.
\end{dem}

\begin{co}\label{co5.1}
Let $j>0$  be an integer. The restriction
of $\rho_{j}$ to $\mathfrak{s}\cap\varpi'^{j}\mathcal{O}'\nu/
\mathfrak{s}\cap\varpi'^{2j}\mathcal{O}'\nu$ is a group isomorphism
onto $S_{j}/S_{2j}$.
\end{co}
\begin{dem}
The result follows from the preceding lemma and from the fact
that $S_{j}T_{2j}/T_{2j}=S_{j}/S_{2j}$.
\end{dem}

\begin{lem}\label{lem5.10}
  We have
  \begin{equation*}
S=S\cap\mu_{q+1}\times S_{1}.
  \end{equation*}
\end{lem}
\begin{dem}
This follows from the fact that for each
$\chi\in\mathbf{X}(T)=\mathbb{Z}[\Gamma']$, one has
$\chi(\mu_{q+1})\subset\mu_{q+1}$ and
$\chi(T_{1})\cap\mu_{q+1}\subset
(1+\varpi'\tilde{\mathcal{O}}'')\cap\mu_{q+1}=\{1\}$.
\end{dem}

\subsection{}\label{5.7}
In this paragraph and until the end of the section we suppose that
$S$ is admissible. It is clear that the sequence of congruence
subgroups $(S_{j})_{j\in\mathbb{N}}$ of $S$ is decreasing.

\begin{lem}\label{lem5.11}
a)  For $j\in\mathbb{N}$ the following assertions are equivalent :

  (i) $j\in 2e\mathbb{N}+I'$,

  (ii) $\varpi'^{j}\mathcal{O'}\cap(\mathfrak{s}/\nu)=
  \varpi'^{j+1}\mathcal{O'}\cap(\mathfrak{s}/\nu)$,

  (iii) $S_{j}=S_{j+1}$.

  b) If $j$ is a positive integer which is not in $2e\mathbb{N}+I'$, then
  $\vert S_{j}/S_{j+1}\vert=q$.

  c)  For $j\in\mathbb{N}\smallsetminus\{0\}$, one has
\begin{equation*}
\begin{cases}
  S_{2j+1}\subsetneq S_{2j}\subset S_{2j-1}\mbox{ if } v''(u)=0,\\
  S_{2j+1}\subset S_{2j}\subsetneq S_{2j-1}\mbox{ if } v''(u)=1.
\end{cases}
\end{equation*} 
Moreover $(S_{2j})_{j\in\mathbb{N}}$ and $(S_{2j+1})_{j\in\mathbb{N}}$
are strictly decreasing sequences of compact open subgroups of $S$
forming a neighbourhood basis of $1$.
 
\end{lem}
\begin{dem}
Assertion c) is an easy consequence of a) and b).
  
The equivalence of the two last assertions in a) follows from Corollary
\ref{co5.1}.

We will prove the equivalence of (i) and (ii) and also assertion b) in
the case where $v''(u)=0$, the proof in the other case being
similar. According to \ref{pr4.6}, in this case $I$ contains all the
even integers between $0$ and $2e-1$.

Let $j\in\mathbb{N}$ and let us write
\begin{equation*}
  j=2e[\frac{j}{2e}] +l.
\end{equation*}
Then $l\in\mathbb{N}$ and $0\leq l\leq2e-1$. Moreover, one easily checks that
\begin{equation}\label{eq39}
  \varpi'^{j}\mathcal{O'}\cap(\mathfrak{s}/\nu)=
  \varpi^{[\frac{j}{2e}]}
  \big(\oplus_{i=0}^{e-1}\varpi^{\epsilon_{j,i}}\mathcal{O}\varpi'^{2i}\oplus
  \oplus_{i\in I_{odd}}\varpi^{\epsilon'_{j,i}}\mathcal{O}\varpi'^{i}\bigg),
\end{equation}
where $I_{odd}=I\cap(2\mathbb{N}+1)$ and
\begin{equation*}
  \epsilon_{j,i}=\begin{cases}
  1\mbox{ if }i\leq\frac{l-1}{2},\\
  0\mbox{ otherwise, }
  \end{cases}
  \mbox{ and }   \epsilon'_{j,i}=\begin{cases}
  1\mbox{ if } i\leq l-1,\\
  0\mbox{ otherwise. }
  \end{cases}
\end{equation*}

Let us write $j+1=2e[\frac{j+1}{2e}] +l'$. If $l=2e-1$, then
$[\frac{j+1}{2e}]=[\frac{j}{2e}]+1$, $l'=0$, $\epsilon_{j,i}=1$
and $\epsilon_{j+1,i}=0$, $0\leq i\leq e-1$ while $\epsilon'_{j,i}=1$,
$0\leq i\leq 2(e-1)$ and $\epsilon'_{j+1,i}=0$, $0\leq i\leq 2e-1$. In
this case, assertion (ii) is satisfied if and only if $2e-1\notin
I_{odd}$ or, equivalently, $j\in 2e\mathbb{N}+I'$.

Suppose now that $l\neq2e-1$. Then one has
$[\frac{j+1}{2e}]=[\frac{j}{2e}]$, $l'=l+1$ and assertion (ii) is
equivalent to $[\frac{l+1}{2}]=[\frac{l'+1}{2}]$ and
$I_{odd}\cap[0,l-1]=I_{odd}\cap[0,l'-1]$.

If $j$ is even, so is $l$ and $[\frac{l'+1}{2}]=[\frac{l+1}{2}]+1$. So
condition (ii) is not satisfied while $j\notin2e\mathbb{N}+I'$.

If $j$ is odd, so is $l$ and $[\frac{l'+1}{2}]=[\frac{l+1}{2}]$. Then
assertion (ii) is satisfied if and only if
$I_{odd}\cap[0,l-1]=I_{odd}\cap[0,l]$, that is if and only if $l\notin
I$ or, equivalently, $j\in 2e\mathbb{N}+I'$. This proves the
equivalence of (i) and (ii).

It follows from Corollary \ref{co5.1} that, for $j$ a positive
integer, $\vert S_{j}/S_{j+1}\vert
=\vert\mathfrak{s}\cap\varpi'^{j}\mathcal{O}'\nu/
\mathfrak{s}\cap\varpi'^{2j}\mathcal{O}'\nu\vert
\leq\vert\mathcal{O}'/\varpi'\mathcal{O}'\vert=q$. But equation
\eqref{eq39} implies that
$\vert\mathfrak{s}\cap\varpi'^{j}\mathcal{O}'\nu/
\mathfrak{s}\cap\varpi'^{2j}\mathcal{O}'\nu\vert$ is a power of
$q$. Assertion b) follows.
\end{dem}

\subsection{}\label{5.8}
In this paragraph we determine the unitary characters of $S$ which
appear in the restriction of the metaplectic representation and their
multiplicity therein.

As in paragraph \ref{5.4}, we define the conductor $c(\chi)$ of a
unitary character $\chi$ of $S$ to be the smallest integer $j$ such
that $S_{j}$ is contained in the kernel of $\chi$.

\begin{lem}\label{lem5.12}
  a) We suppose that $v''(u)=0$. If $\chi$ is a unitary character of $T$ with
  odd conductor $c(\chi)\geq3$, one has $c(\chi_{\vert S})=c(\chi)$

  b) We suppose that $v''(u)=1$. If $\chi$ is a unitary character of $T$ with
  even conductor, one has $c(\chi_{\vert S})=c(\chi)$.
\end{lem}
\begin{dem}
%We give the proof in the case a) and leave the other case to the
  %reader.
We may suppose that $\chi$ is not trivial. Let us write
$c(\chi)=2j+1-v''(u)$ with $j\geq1$ an integer.  It follows from
Proposition \ref{pr2.1} that there exists $b\in\mathcal{O}'^{\times}$
such that
\begin{equation*}
  \chi(\rho(X))=\psi(-\frac{1}{4}\sideset{}{_{k''/k}}\tr\varpi'^{\mu-c(\chi)}
  buX)\mbox{, }X\in\varpi'^{j}\mathcal{O}'\nu.
\end{equation*}
As $b\in\mathcal{O}'^{\times}$, we may write
$b=b_{0}+\sum_{l=1}^{2e-1}b_{l}\varpi'^{l}$ with
$b_{l}\in\mathcal{O}$, $1\leq l\leq 2e-1$ and
$b_{0}\in\mathcal{O}^{\times}$.

Now it follows from Corollary \ref{co5.1} that $\rho$ induces the
isomorphism $\rho_{c(\chi)-1}$ from
$\mathfrak{s}\cap\varpi'^{c(\chi)-1}\mathcal{O}'\nu/
\mathfrak{s}\cap\varpi'^{2(c(\chi)-1)}\mathcal{O}'\nu$ onto
$S_{c(\chi)-1}/S_{2(c(\chi)-1)}$. Thus, if
$X\in\mathfrak{s}\cap\varpi'^{c(\chi)-1}\mathcal{O}'\nu$, there exists $x\in
S_{c(\chi)-1}$ such that $x\in\rho(X)T_{2(c(\chi)-1)}$. Then one has
$\chi(x)=\chi(\rho(X))$.
 
Now, according to the preceding lemma, we have $S_{c(\chi)}\subsetneq
S_{c(\chi)-1}$.  So we see that it suffices to prove that there exists
$X\in\mathfrak{s}\cap\varpi'^{c(\chi)-1}\mathcal{O}'\nu$ such that
$\chi(\rho(X))\neq1$.

It follows from Proposition
\ref{pr4.6} that for any $\eta_{0}\in\mathcal{O}^{\times}$ there
exists $X\in\mathfrak{s}$ such that
$X\in\eta_{0}\varpi'^{c(\chi)-1}\nu+\varpi'^{c(\chi)}\mathcal{O}'\nu$. Then,
$X\in\mathfrak{s}\cap\varpi'^{c(\chi)-1}\mathcal{O}'\nu$ and one has
\begin{eqnarray*}
  \chi(\rho(X))&=&
  \psi(-\frac{1}{4}\sideset{}{_{k''/k}}\tr\varpi'^{\mu-1}b_{0}\eta_{0}u\nu)=
  \psi(\frac{1}{2}\sideset{}{_{k'/k}}
  \tr\varpi'^{2e\lambda_{\psi}-\delta-1}b_{0}\eta_{0}d)\\
  &=&\psi(\frac{1}{2}\varpi^{\lambda_{\psi}}\sideset{}{_{k'/k}}\tr\varpi'^{-\delta-1}
  b_{0}\eta_{0}d),
\end{eqnarray*}
which is $\neq1$ for a specific choice of $\eta_{0}\in\mathcal{O}^{\times}$.
\end{dem}

Recall that $I=I_{\underline{d}}$ with
$I'=I'_{\underline{d}}:=\{0,\ldots,2e-1\}\smallsetminus
I_{\underline{d}}$. As $S$ is admissible, according to Proposition
\ref{pr4.6}, on has $I'\subset 2\mathbb{N}+1-v''(u)$.

Then we put
\begin{equation}\label{eq39.5}
  \alpha(j)=\vert\{1,\ldots,j-2\}\cap(2e\mathbb{N}+I')\vert\mbox{,
  }j\in\mathbb{N}.
\end{equation}

\begin{pr}\label{pr5.3}
a) We suppose that $v''(u)=0$. A non trivial character $\chi$ of $S$
appears in the restriction to $S$ of the metaplectic representation if
and only if its conductor is odd. Then it appears with multiplicity
\begin{equation}\label{eq40}
m(\chi)=\vert\mu_{q+1}/\mu_{q+1}\cap S\vert
q^{\alpha(c(\chi))}.
\end{equation}

The trivial character of $S$ appears in the restriction to $S$ of the
metaplectic representation if and only if $\mu_{q+1}\cap
S\neq\mu_{q+1}$. In this case, one has
\begin{equation}\label{eq41}
m(1)=\vert\mu_{q+1}/\mu_{q+1}\cap S\vert-1.
\end{equation}

b) We suppose that $v''(u)=1$. A character $\chi$ of $S$ appears in the
restriction to $S$ of the metaplectic representation if and only if
its conductor is even. If it is non-trivial, it appears with
multiplicity
\begin{equation}\label{eq42}
m(\chi)=\vert\mu_{q+1}/\mu_{q+1}\cap S\vert
q^{\alpha(c(\chi))}.
\end{equation}

The trivial character appears with multiplicity $1$.
\end{pr}
\begin{dem}
  We give the proof in the case a) and leave the other case to the reader.

  First, as $k'$ is tamely ramified over $k$, we have $\delta=2e-1$ and
  thus $\mu$ is odd. It follows from Theorem \ref{theo2.1} that the
  characters of $T$ appearing in the restriction of the metaplectic
  representation are the ones with odd conductor, each one appearing
  with multiplicity one.

  Then, according to Lemma \ref{lem5.12}, the characters of $S/S_{1}$
  appearing in the restriction of the metaplectic representation are
  the restriction to $S/S_{1}$ of the characters of
  $T/T_{1}$ appearing therein. Moreover, if $\eta$ is a character of $S/S_{1}$,
  $m(\eta)$ is the number of characters $\chi$ of $T/T_{1}$ which are
 non-trivial and such that $\chi_{\vert S}=\eta$. One
  immediately sees that
\begin{equation*}
  m(\eta)=\begin{cases}
  \vert\mu_{q+1}/\mu_{q+1}\cap S\vert&\mbox{ if }\eta\neq1,\\
  \vert\mu_{q+1}/\mu_{q+1}\cap S\vert-1&\mbox{ otherwise. }
  \end{cases}
\end{equation*}

Now, still according to Lemma \ref{lem5.12} and Theorem \ref{theo2.1},
the characters of $S$ with conductor at least $2$ appearing in the
restriction of the metaplectic representation are the ones with odd
conductor. Moreover, if $\eta$ is such a character of $S$, its
multiplicity in the restriction of the metaplectic representation is
the number of characters $\chi$ of $T$ such that $c(\chi)=c(\eta)$ and
$\chi_{\vert S}=\eta$.

Let us write $c(\eta)=2j+1$ with $j\geq1$. Then we clearly have
$m(\eta)=\vert T/T_{2j+1}\vert/\vert S/S_{2j+1}\vert$. But $\vert
T/T_{2j+1}\vert=\vert\mu_{q+1}\vert\times\vert
T_{1}/T_{2j+1}\vert=\vert\mu_{q+1}\vert q^{2j}$ while, according to
Lemma \ref{lem5.11}, $\vert S/S_{2j+1}\vert=\vert\mu_{q+1}\cap
S\vert\times\vert S_{1}/S_{2j+1}\vert= \vert\mu_{q+1}\cap S\vert
q^{2j-\vert\{1,\cdots,2j-1\}\cap(2e\mathbb{N}+I')\vert}$. Thus the result.
\end{dem}

\subsection{}\label{5.9}
In the last paragraphs of this section, we still consider the
following cases
\begin{description}
\item[Case A] $k''$ is ramified over $k'$ and $k'$ unramified over $k$,
\item[Case B] $k''$ is unramified over $k'$ and $k'$ tamely and
 totally ramified over $k$. 
\end{description}

In these cases, we relate the multiplicity of a character $\chi$ of
$S$ occurring in the restriction to $S$ of the metaplectic
representation to the volume of the symplectic reduction of the
inverse image under the momentum map of a linear form on
$\mathfrak{s}$ associated to $\chi$.

\begin{lem}\label{lem5.13}
The momentum map $\phi:k''\rightarrow\mathfrak{s}^{*}$ is a submersion
at each point of $k''^{\times}$. Moreover, if $w\in k''^{\times}$,
$\mathfrak{s}\cdot w$ is the kernel of the restriction of $\beta$ to
$\ker\diff\phi(w)$.
\end{lem}
\begin{dem}
Recall that the momentum map is given by 
\begin{equation}\label{eq43}
  \langle\phi(w),X\rangle=-\frac{1}{2}\beta(X\cdot w,w)
  \mbox{, }w\in k''\mbox{, }X\in\mathfrak{s}.
\end{equation}
Let $w\in k''^{\times}$. For $h\in k''$ and $X\in\mathfrak{s}$, we
have
\begin{equation*}
  \langle\diff\phi(w).h,X\rangle=-\beta(X\cdot w,h).
\end{equation*}
It follows that $\ker\diff\phi(w)$ is the orthogonal of
$\mathfrak{s}\cdot w$ with respect to the symplectic form $\beta$. As $w\in
k''^{\times}$, the map $X\mapsto X\cdot w$ is an injection from
$\mathfrak{s}$ into $k''$. Thus $\dim\ker\diff\phi(w)=\dim
k''-\dim\mathfrak{s}$ showing that $\phi$ is a submersion at $w$.

As a consequence of \eqref{eq43} we have $\phi(t\cdot w)=\phi(w)$, $t\in T$,
$w\in W$. Then the last assertion follows from the fact that
$\mathfrak{s}\cdot w$ is the orthogonal of $\ker\diff\phi(w)$ with respect
to $\beta$.
\end{dem}

\begin{co}\label{co5.2}
For every $g\in\im\phi$, $\phi^{-1}(g)$ is a compact smooth
algebraic submanifold in $k''$.
\end{co}
\begin{dem}
As $\phi$ is a proper map and $\phi^{-1}(0)=\{0\}$, this follows from
the previous lemma.
\end{dem}

Let $g\in\im\phi$. As we have seen, the natural
action of $T$ and thus the one of $S$ in $k''$ leaves 
$\phi^{-1}(g)$ invariant. We will see that the set of $S$-orbits
$S\backslash\phi^{-1}(g)$ is a compact manifold and that the
symplectic form $\beta$ induces naturally on it a structure of
symplectic manifold, which is called the symplectic reduction of
$\phi^{-1}(g)$.

%{\color{red}
Let $M$ be the $\mathbb{Z}[\Gamma']$-submodule of
$\mathbf{X}(T)=\mathbb{Z}[\Gamma']$ whose elements are the characters
which are trivial on $S$ and $N$ be the direct sum of the submodules
$M_{d}$ which are not contained in $M$. Let $\overline{N}$ be the set
of elements $\gamma$ of $\mathbf{X}(T)$ such that there exists a
non-zero integer $a\in\mathbb{Z}$ such that $a\gamma\in N$. As
$\mathbb{Z}$-modules, $\overline{N}$ and $N$ have the same rank which
is $\dim\mathfrak{s}$.

%Moreover one has $\tilde{N}\cap M=0$. It
%follows that the image of $M$ (resp. $\tilde{N}$) into
%$\mathbb{Z}[\Gamma']/\tilde{N}$ (resp. $\mathbb{Z}[\Gamma']/M$) by the
%natural projection is a free $\mathbb{Z}$-submodule of same rank.

Then $S'=\{x\in T\vert \gamma(x)=1\mbox{, }\gamma\in\overline{N}\}$ is
a subtorus of $T$ defined over $k$ and its Lie algebra
$\mathfrak{s}'=\cap_{\gamma\in\overline{N}}\ker \gamma$ has dimension
$\codim\mathfrak{s}$. Moreover, as the $\mathbb{Z}$-module
$M\oplus\overline{N}$ as same rank than $\mathbf{X}(T)$,
$\mathfrak{s}'\cap\mathfrak{s}={0}$. Thus $\mathfrak{s}'/\nu$ is a
$\Gamma'$-invariant supplementary subspace of $\mathfrak{s}/\nu$ in
$k'$. Then, it follows from Lemma \ref{lem5.4} in case A and from
\eqref{eq23} applied to $\mathfrak{s}'$ and \eqref{eq24} in case B
that $\mathfrak{s}'/\nu=(\mathfrak{s}/\nu)^{\perp}$.

Recall that $S'_{1}=S'\cap T_{1}$ and $\varepsilon_{S'}=2$ if
$S'\cap\{\pm1\}=\{1\}$ and $1$ otherwise in case A, while
$S'=\mu_{q+1}\cap S'\times S'_{1}$ in case B.

\begin{lem}\label{lem5.14}
The map $(s,s')\mapsto ss'$ is an isomorphism of groups from
$S_{1}\times S'_{1}$ onto $T_{1}$ in cases A) and B), and also from
$S\times S'$ onto $T$ in case A.

As a consequence we have $\varepsilon_{S}\varepsilon_{S'}=2$ in case A.
\end{lem}
\begin{dem}
Let us denote by $\alpha$ the map $(s,s')\mapsto ss'$ which is a
morphism of algebraic groups over $k$, namely $S\times S'$ and $T$. We
also consider $\hat{\alpha}$ the dual morphism from
$\mathbf{X}(T)=\mathbb{Z}[\Gamma']$ into
$\mathbf{X}(S)\times\mathbf{X}(S')=
\mathbb{Z}[\Gamma']/M\times\mathbb{Z}[\Gamma']/\overline{N}$. It is
clear that for $\gamma\in\mathbb{Z}[\Gamma']$,
$\hat{\alpha}(\gamma)=(p(\gamma),p'(\gamma))$ where $p$ (resp. $p'$)
stands for the natural projection from $\mathbb{Z}[\Gamma']$ onto
$\mathbb{Z}[\Gamma']/M$
(resp. $\mathbb{Z}[\Gamma']/\overline{N}$). Then we see that
$\hat{\alpha}(\overline{N})=p(\overline{N})\times\{0\}$ and that
$\hat{\alpha}(M)=\{0\}\times p'(M)$ and consequently that
$p(\overline{N})\times p'(M)\subset\im\hat{\alpha}$. It follows from
this and from Corollary \ref{coap3} that the quotient
$(\mathbb{Z}[\Gamma']/M\times\mathbb{Z}[\Gamma']/\overline{N})/\im\hat{\alpha}$
is a torsion $\mathbb{Z}$-module in which the order of any non-zero
element is a divisor of $2f$ in case A or $2e$ in case B. From that
follows that there exists a basis over $\mathbb{Z}$
$(e_{1},\ldots,e_{\vert\Gamma'\vert})$ of
$\mathbb{Z}[\Gamma']/M\times\mathbb{Z}[\Gamma']/\overline{N}$ and
strictly positive integers $a_{1},\ldots,a_{\vert\Gamma'\vert}$ all
coprime with $p$ such that
$(a_{1}e_{1},\ldots,a_{\vert\Gamma'\vert}e_{\vert\Gamma'\vert})$ is a
basis over $\mathbb{Z}$ of $\im\hat{\alpha}$.

Besides, up to sign the determinant of the matrix of $\hat{\alpha}$ in
a pair of basis of $\mathbb{Z}[\Gamma']$ and of
$\mathbb{Z}[\Gamma']/M\times\mathbb{Z}[\Gamma']/\overline{N}$ does not
depend upon the choice of these basis. We denote it by
$\nu(\alpha)$. Then we have $\nu(\alpha)=\prod_{1\leq
  l\leq\vert\Gamma'\vert}a_{l}$, so that $\nu(\alpha)$ is coprime with
$p$. Thus it follows from \cite[Proposition 2.2.2, p.117 and Theorem
  2.3.1, p. 119]{ono-1961} that $\alpha$ induces an isomorphism from
$S_{1}\times S'_{1}$ onto $T_{1}$, and that, in case A, $\ker\alpha$
and $T/S.S'$ have the same cardinal. The lemma follows easily.
\end{dem}

\subsection{}\label{5.10}

In this paragraph, we suppose we are in case A. Then, we describe the
symplectic reduction of the inverse image under the momentum map of a
linear form on $\mathfrak{s}$ contained in $\im\phi$.

Let $g$ be a non-zero linear form on $\mathfrak{s}$ which is contained
into the image of the momentum map $\phi$. Then there exists
$a\in\mathcal{O}''^{\times}$ and $j\in\mathbb{Z}$ such that, putting
$w_{0}=\varpi''^{\mu-j}a$, one has
\begin{eqnarray}
  \langle g,\eta\varpi''\rangle&=&
  \frac{(-1)^{\mu-j}}{2}\sideset{}{_{k'/k}}\tr\varpi^{\mu-j}
  \sideset{}{_{k''/k'}}\Norm(a\varpi'')\eta\label{eq44}\\
  &=&\langle\phi(w_{0}),\eta\varpi''\rangle\mbox{,
  }\eta\in\mathfrak{s}/\varpi''.\notag
\end{eqnarray}
By definition, for $\alpha\in\Theta(\mathsf{p}(\Norm_{k''/k'}(a)))$
there exists $\gamma_{\alpha}\in\mu_{q'-1}$, unique up multiplication
by $\pm1$, such that
$\mathsf{p}(\gamma_{\alpha}^{2})=\alpha+\mathsf{p}(\Norm_{k''/k'}(a_{0}))$. Then
it follows from \eqref{eq33} that there exists
$b_{\alpha}\in\mathcal{O}'^{\times}$ such that
$b_{\alpha}^{2}=\Norm_{k''/k'}(b_{\alpha})\in
\Norm_{k''/k'}(a)+(\mathsf{p}^{-1}(\alpha)\cap(\mathfrak{s'}/\varpi''))$
(recall that $\mathfrak{s}'/\varpi''=(\mathfrak{s}/\varpi'')^{\perp}$).
%Then one has
%$\mathsf{p}^{-1}(\alpha)\cap(\mathfrak{s}/\varpi'')^{\perp}=
%\Norm_{k''/k'}(b_{\alpha})-\Norm_{k''/k'}(a)+\varpi\mathcal{O}'\cap
%(\mathfrak{s}/\varpi'')^{\perp}$.
Clearly, we may choose $b_{\alpha}$ such that
$b_{\alpha}\in\gamma_{\alpha}(1+\varpi\mathcal{O}')$ and, for $\alpha=0$,
$\Norm_{k''/k'}(b_{0})=\Norm_{k''/k'}(a)$.

Then, denoting by $\sqrt{}$ the inverse of the bijection $z\mapsto
z^{2}$ on $1+\varpi\mathcal{O}'$, we put
\begin{equation}\label{eq45}
  X_{\alpha}=
  \varpi''^{\mu-j}\{b_{\alpha}t\sqrt{1+b_{\alpha}^{-2}\varpi\eta}\,\vert\,
  t\in T\mbox{, }
  \eta\in\mathcal{O}'\cap(\mathfrak{s}'/\varpi'')\}.
\end{equation}
It is clear that $X_{\alpha}$ depends only upon $\alpha$.

\begin{lem}\label{lem5.15}
Let $g\in\im\phi\smallsetminus\{0\}$, $a\in\mathcal{O}''^{\times}$ and
$j\in\mathbb{Z}$ satisfying \eqref{eq44}.

Then we have
\begin{equation*}
\phi^{-1}(g)=\sqcup_{\alpha\in\Theta(\mathsf{p}(\Norm_{k''/k'}(a)))}X_{\alpha}.
\end{equation*}
and the $X_{\alpha}$, $\alpha\in\Theta(\mathsf{p}(\Norm_{k''/k'}(a)))$, are
analytic submanifolds of $k''$ and open subsets in $\phi^{-1}(g)$.

Let $\alpha\in\Theta(\mathsf{p}(\Norm_{k''/k'}(a)))$. The map
\begin{equation}\label{eq46}
  \rho_{\alpha}:(t,\eta)\mapsto
  \varpi''^{\mu-j}b_{\alpha}t\sqrt{1+b_{\alpha}^{-2}\varpi\eta}
\end{equation}
induces a diffeomorphism from
$T\times\mathcal{O}'\cap(\mathfrak{s}'/\varpi'')$ onto
$X_{\alpha}$.
\end{lem}

\begin{dem}
Let $w\in\phi^{-1}(g)$. We may write $w=\varpi''^{\mu-l}b$ with
$b\in\mathcal{O}''^{\times}$ and $l\in\mathbb{Z}$. We easily see that
condition $\phi(w)=g$ is equivalent to
\begin{equation*}
  \varpi^{\mu-j}\sideset{}{_{k''/k'}}\Norm(a\varpi'')-
  \varpi^{\mu-l}\sideset{}{_{k''/k'}}\Norm(b\varpi'')
  \in(\mathfrak{s}/\varpi'')^{\perp}.
\end{equation*}
Suppose that $l\neq j$. Without loss of generality, we may suppose
that $l<j$. We then have
\begin{equation*}
  \sideset{}{_{k''/k'}}\Norm(a)-\sideset{}{_{k''/k'}}\Norm(b)\varpi^{j-l}\in
  (\mathfrak{s}/\varpi'')^{\perp}.
\end{equation*}
But one has
\begin{equation*}
\sideset{}{_{k''/k'}}\Norm(a)-\sideset{}{_{k''/k'}}\Norm(b)\varpi^{j-l}\in
\sideset{}{_{k''/k'}}\Norm(a)(1+\varpi\mathcal{O}')\subset
\sideset{}{_{k''/k'}}\Norm(\mathcal{O}''^{\times}),
\end{equation*}
in contradiction with the hypothesis that $S$ is admissible.

Thus, we have $l=j$ and then $w\in\phi^{-1}(g)$ if and only if
$\Norm_{k''/k'}(b)\in\Norm_{k''/k'}(a)+
(\mathcal{O}'\cap(\mathfrak{s}/\varpi'')^{\perp})$. Using
\eqref{eq33}, we easily see that this last condition means exactly
that there exists $\alpha\in\Theta(\mathsf{p}(\Norm_{k''/k'}(a)))$
such that $w\in X_{\alpha}$. Then it is clear that $\phi^{-1}(g)$ is
the disjoint union of the $X_{\alpha}$,
$\alpha\in\Theta(\mathsf{p}(\Norm_{k''/k'}(a)))$.

For $\gamma\in\mu_{q'-1}$ let $U_{\gamma}=\{w\in k''^{\times}\vert
\varpi''^{-v''(w)}w\in\gamma(1+\varpi''\mathcal{O''})\}$. Then the
$U_{\gamma}$, $\gamma\in\mu_{q'-1}$, are disjoint open subsets of $k''$
and we have $k''^{\times}=\cup_{\gamma\in\mu_{q'-1}}U_{\gamma}$. Moreover,
for any $\alpha\in\Theta(\mathsf{p}(\Norm_{k''/k'}(a)))$, we have
$X_{\alpha}=\phi^{-1}(g)\cap U_{\gamma_{\alpha}}$, showing that
$X_{\alpha}$ is an analytic submanifold of $k''$ and open in $\phi^{-1}(g)$.

Let $\alpha\in\Theta(\mathsf{p}(\Norm_{k''/k'}(a)))$. One checks
easily that $\rho_{\alpha}$ is an analytic map and also a bijection
from $T\times\mathcal{O}'$ onto $X_{\alpha}$. %For
%$\eta\in\mathcal{O}'$, let us put

Moreover, a straightforward computation gives for $t\in T$,
$\eta\in\mathcal{O}'\cap(\mathfrak{s}'/\varpi'')$, $X\in\mathfrak{t}$
and $\xi\in\mathfrak{s}'/\varpi''$
\begin{equation}\label{eq48}
  \diff\rho_{\alpha}(t,\eta)(tX,\xi)
  =t(\rho_{\alpha}(1,\eta)X+
  \frac{1}{2}\varpi^{\mu-j}\rho_{\alpha}(1,\eta)^{-1}\varpi\xi).
\end{equation}
As $\mathfrak{t}\subset k'\varpi''$, it follows that $
\diff\rho_{\alpha}(t,\eta)$ is injective, while it is clear that the
analytic manifolds $T\times\mathcal{O}'\cap(\mathfrak{s}'/\varpi'')$
and $X_{\alpha}$ have same dimension over $k$. Thus $\rho_{\alpha}$ is
a diffeomorphism from
$T\times\mathcal{O}'\cap(\mathfrak{s}'/\varpi'')$ onto
$X_{\alpha}$.
\end{dem}

Let us put
$Y_{\alpha}=\rho_{\alpha}(S'\times\mathcal{O}'\cap(\mathfrak{s}'/\varpi''))$. It
follows from the preceding lemma that $Y_{\alpha}$ is an analytic
submanifold of $X_{\alpha}$. Moreover, the restriction
$\tilde{\rho}_{\alpha}$ of $\rho_{\alpha}$ to
$S'\times\mathcal{O}'\cap(\mathfrak{s}'/\varpi'')$ is a diffeomorphism
onto $Y_{\alpha}$.

\begin{lem}\label{lem5.15.5}
The manifold $Y_{\alpha}$ is transverse to the $S$-orbits into
$X_{\alpha}$. The natural projection from $X_{\alpha}$ onto
$S\backslash X_{\alpha}$ induces an homeomorphism from $Y_{\alpha}$
onto $S\backslash X_{\alpha}$. Via this isomorphism, $S\backslash
X_{\alpha}$ becomes an analytic manifold we identify to
$Y_{\alpha}$. 

The symplectic form $\beta$ induces a closed analytic two-differential
form on $X_{\alpha}$ whose kernel at each point is the tangent space
of the $S$-orbit through this point. It also induces a symplectic form
$\overline{\beta}_{\alpha}$ on $S\backslash X_{\alpha}=Y_{\alpha}$.

The pullback $\tilde{\beta}$  of $\overline{\beta}_{\alpha}$ under
$\tilde{\rho}_{\alpha}$ does not depend upon $\alpha$.  For $s\in S'$ and
$\eta\in\mathcal{O}'\cap(\mathfrak{s}'/\varpi'')$,
$X,X'\in\mathfrak{s'}$ and $\xi,\xi'\in\mathfrak{s'}/\varpi''$ we have
\begin{equation}\label{eq47}
  \tilde{\beta}_{(s,\eta)}((sX,\xi),(sX',\xi'))
  \!=\!\frac{(-1)^{\mu-j}}{2}
  \varpi^{\mu+2-j}\sideset{}{_{k'/k}}\tr(\frac{X}{\varpi''}\xi'\!-\!
  \frac{X'}{\varpi''}\xi).
\end{equation}

The volume form $\wedge^{\dim(\mathfrak{s}')}\tilde{\beta}$ is
translation invariant on the direct product of groups
$S'\times\mathcal{O}'\cap(\mathfrak{s}'/\varpi'')$.
\end{lem}

\begin{dem}
The fact that $Y_{\alpha}$ is transverse to the $S$-orbits into
$X_{\alpha}$ is a consequence of Lemma \ref{lem5.14}. From this
follows that $S\backslash X_{\alpha}$ is an analytic manifold
naturally isomorphic to $Y_{\alpha}$ and that $\tilde{\rho}_{\alpha}$
is a diffeomorphism from
$S'\times\mathcal{O}'\cap(\mathfrak{s}'/\varpi'')$ onto
$Y_{\alpha}$. Lemma \ref{lem5.13} implies that $\beta$ induces an
analytic two-differential form on $X_{\alpha}$ which like $\beta$ is
closed and whose kernel at each point is the tangent space of the
$S$-orbit through this point, and also a symplectic structure on
$S\backslash X_{\alpha}=Y_{\alpha}$.

Formula \eqref{eq47} follows by a straightforward computation from
\eqref{eq48}. The last assertion is then clear.
\end{dem}

\begin{co}\label{co5.3}
The space of $S$-orbits $S\backslash\phi^{-1}(g)$ is an analytic
manifold which is the disjoint union of $S\backslash X_{\alpha}$,
$\alpha\in\Theta(\mathsf{p}(\Norm_{k''/k'}(a)))$, which are themselves
analytic manifolds. Moreover, the form $\beta$ induces a symplectic
two form $\overline{\beta}$ on $S\backslash\phi^{-1}(g)$ whose
restriction to each $S\backslash X_{\alpha}$ is $\overline{\beta}_{\alpha}$.
\end{co}

\subsection{}\label{5.11}

In this paragraph, we suppose we are in case B. Then, we describe the
symplectic reduction of the inverse image under the momentum map of a
linear form on $\mathfrak{s}$  contained in $\im\phi$.

Let $g$ be a non-zero linear form on $\mathfrak{s}$ which is contained
into the image of $\phi$. Then there exists
$a\in\mathcal{O}''^{\times}$ and $j\in\mathbb{Z}$ having same parity
as $\mu$ such that, putting $w_{0}=\varpi'^{\frac{\mu-j}{2}}a$, one has
\begin{eqnarray}\label{eqa}
\langle g,\eta u\rangle&=&
-\frac{1}{2}\sideset{}{_{k'/k}}\tr(\varpi'^{\mu-j}\sideset{}{_{k''/k'}}
\Norm(au)\eta)\\
&=&\langle\phi(w_{0}),\eta u\rangle
\mbox{,
}\eta\in\mathfrak{s}/u.\notag
\end{eqnarray}

As the integer $\mu-j$ is even, it follows from Proposition \ref{pr4.6} that
\begin{equation*}
\mathcal{O}'\cap\varpi'^{j-\mu-2v''(u)}(\mathfrak{s}'/u)
=\mathcal{O}'\cap\varpi'^{j-\mu-2v''(u)}(\mathfrak{s}/u)^{\perp}
\subset\varpi'\mathcal{O}'.
\end{equation*}
Thus, the map
$\rho_{w_{0}}:T\times
\mathcal{O}'\cap\varpi'^{j-\mu-2v''(u)}(\mathfrak{s}'/u)\rightarrow W$ such that
\begin{equation}\label{eqb}
\rho_{w_{0}}(t,\eta)=\varpi'^{\frac{\mu-j}{2}}at\sqrt{1+\sideset{}{_{k''/k'}}
\Norm(a\nu)^{-1}\eta}
\end{equation}
is well defined (as above, $\sqrt{}$ denotes the inverse of the
bijection $z\mapsto z^{2}$ on $1+\varpi'\mathcal{O}'$).

\begin{lem}\label{lema}
Let $g\in\im\phi\smallsetminus\{0\}$, $a\in\mathcal{O}''^{\times}$ and
$j\in\mathbb{Z}$ satisfying \eqref{eqa}.

Then, $\rho_{w_{0}}$ is a diffeomorphism from
$T\times\mathcal{O}'\cap\varpi'^{j-\mu-2v''(u)}(\mathfrak{s}'/u)$ onto
$\phi^{-1}(g)$.
\end{lem}
\begin{dem}
Let $w\in W\smallsetminus\{0\}$. Then there exists $l\in\mathbb{Z}$ and
$b\in\mathcal{O}''^{\times}$ such that $w=\varpi'^{l}b$. The condition
$w\in\phi^{-1}(g)$ implies that
\begin{equation*}
  \varpi'^{2l}\sideset{}{_{k''/k'}}\Norm(bu)-
  \varpi'^{\mu-j}\sideset{}{_{k''/k'}}\Norm(au)\in(\mathfrak{s}/u)^{\perp}.
\end{equation*}
If $2l\neq\mu-j$, arguing like in the first step of the proof of Lemma
\ref{lem5.15} we see that $\sideset{}{_{k''/k'}}\Norm(k''^{\times})
\cap(\mathfrak{s}/u)^{\perp}\neq\emptyset$, in contradiction with the
admissibility of $S$.

Thus, we have $w=\varpi'^{\frac{\mu-j}{2}}b$ and the condition
$w\in\phi^{-1}(g)$ is then equivalent to
$\varpi'^{\mu-j}(\Norm_{k''/k'}(bu)-\Norm_{k''/k'}(au))\in(\mathfrak{s}/u)^{\perp}$
or, as $\nu^{2}\in\mathcal{O}^{\times}$, to
$\Norm_{k''/k'}(b)-\Norm_{k''/k'}(a)\in
\mathcal{O}'\cap\varpi'^{j-\mu-2v''(u)}(\mathfrak{s}/u)^{\perp}
=\mathcal{O}'\cap\varpi'^{j-\mu-2v''(u)}(\mathfrak{s}'/u)$. Then a
straightforward computation shows that $\rho_{w_{0}}$ is a bijection
from $T\times\mathcal{O}'\cap\varpi'^{j-\mu-2v''(u)}(\mathfrak{s}'/u)$
onto $\phi^{-1}(g)$.

By a direct computation we get for $t\in T$,
$\eta\in\mathcal{O}'\cap\varpi'^{j-\mu-2v''(u)}(\mathfrak{s}'/u)$,
$X\in\mathfrak{t}$ and
$\xi\in\varpi'^{j-\mu-2v''(u)}(\mathfrak{s}'/u)$ 
\begin{equation}\label{eqc}
  \diff\rho_{w_{0}}(t,\eta)(tX,\xi)=
  t(\rho_{w_{0}}(1,\eta)X-
  \frac{1}{2d}\varpi'^{\mu-j}(\rho_{w_{0}}(1,\eta)^{\tau})^{-1}\xi).
\end{equation}
where $\tau$ is the conjugation of $k''$ over $k'$. It follows readily
that $\rho_{w_{0}}$ is an immersion and a diffeomorphism onto 
$\phi^{-1}(g)$.
\end{dem}

It follows from the preceding lemma that $\tilde{\rho}_{w_{0}}$, the
restriction to $\mu_{q+1}\times S'_{1}\times
\mathcal{O}'\cap\varpi'^{j-\mu-2v''(u)}(\mathfrak{s}'/u)$ of
$\rho_{w_{0}}$, is a diffeomorphism onto a compact submanifold $X$ of
$\phi^{-1}(g)$. Moreover, $X_{1}=\rho_{w_{0}}(S'_{1}\times
\mathcal{O}'\cap\varpi'^{j-\mu-2v''(u)}(\mathfrak{s}'/u))$ is a
submanifold of $\phi^{-1}(g)$, and the map $(t,x)\mapsto tx$ is a
diffeomorphism from $\mu_{q+1}\times X_{1}$ onto $X$.

\begin{lem}\label{lem5.16.5}
The manifold $X$ is transverse to the $S_{1}$-orbits therein. The
restriction of the symplectic form $\beta$ to $X$ is also symplectic.
 
The quotient space $S\backslash\phi^{-1}(g)$ of $S$-orbits into
$\phi^{-1}(g)$ is an analytic manifold isomorphic to
$\mu_{q+1}/\mu_{q+1}\cap S\times X_{1}$. 

The restriction of $\beta$ to $\phi^{-1}(g)$ is a closed analytic
$T$-invariant two-differen\-tial form whose kernel at each point is
the tangent space of the $S$-orbit through this point.  It also
induces a $\mu_{q+1}/\mu_{q+1}\cap S$-invariant symplectic form
$\overline{\beta}$ on $S\backslash\phi^{-1}(g)=\mu_{q+1}/\mu_{q+1}\cap
S\times X_{1}$ which satisfies $\overline{\beta}_{\vert
  X_{1}}=\beta_{\vert X_{1}}$.  If $\tilde{\beta}$ denotes the
pullback of $\overline{\beta}$ under the restriction of
$\tilde{\rho}_{w_{0}}$ to
$S'_{1}\times\mathcal{O}'\cap\varpi'^{j-\mu-2v''(u)}(\mathfrak{s}'/u)$,
for $s\in S'_{1}$,
$\eta\in\mathcal{O}'\cap\varpi'^{j-\mu-2v''(u)}(\mathfrak{s}'/u)$,
$X,X'\in\mathfrak{s}'$ and
$\xi,\xi'\in\varpi'^{j-\mu-2v''(u)}(\mathfrak{s}'/u)$, one has
\begin{equation}\label{eqd}
  \tilde{\beta}_{(s,\eta)}((sX,\xi),(sX',\xi'))=
  -\frac{1}{2}\sideset{}{_{k'/k}}\tr
  \varpi'^{\mu+v''(u)-j}(\frac{X}{\nu}\xi'-\frac{X'}{\nu}\xi).
\end{equation}
The volume form $\wedge^{\dim(\mathfrak{s}')}\tilde{\beta}$ is
translation invariant on the direct product of groups
$S'_{1}\times\mathcal{O}'\cap\varpi'^{j-\mu-2v''(u)}(\mathfrak{s}'/u)$.
\end{lem}
\begin{dem}
The fact that $X$ is transverse to the $S$-orbits into
$\phi^{-1}(g)$ follows from Lemma \ref{lem5.14}. The fact that the
space of $S$-orbits into $\phi^{-1}(g)$ and $\mu_{q+1}/\mu_{q+1}\cap
S\times X_{1}$ are isomorphic analytic manifolds is clear.

It follows from Lemma \ref{lem5.13} that $\beta$ induces an analytic
symplectic form $\overline{\beta}$ on $S\backslash\phi^{-1}(g)$. Then
it is clear that $\overline{\beta}_{\vert X_{1}}=\beta_{\vert X_{1}}$.

Formula \eqref{eqd} follows by a direct computation from \eqref{eqc}.
The last assertion is then clear.
\end{dem}

\subsection{}\label{5.12}

In this paragraph, we keep the notations introduced in paragraphs
\ref{5.10} and \ref{5.11}. According to \cite[paragraph
  7.4]{igusa-2000} and \cite[pages 9-10]{popa-2011}, once chosen a
Haar measure on the field $k$, the volume form
$\wedge^{2\dim(\mathfrak{s}')}\overline{\beta}$ naturally induces a
measure denoted by $\diff\mu_{\beta}$ on the symplectic manifold
$S\backslash\phi^{-1}(g)$. We will see that there is a canonical
choice for the Haar measure on $k$ with respect to the chosen
character $\psi$ of $k$ and then, we will associate to the symplectic
form $\tilde{\beta}$ a Haar measure $\diff_{\beta}(s,\nu)$ on
$S'\times\mathcal{O}'\cap(\mathfrak{s}'/\varpi'')$
(resp. $S'_{1}\times\mathcal{O}'\cap\varpi'^{j-\mu-2v''(u)}(\mathfrak{s}'/u)$)
in case A (resp. in case B). Then, in case A, $\diff\mu_{\beta}$ will
be the measure on $S\backslash\phi^{-1}(g)$ whose restriction to each
$S\backslash X_{\alpha}$,
$\alpha\in\Theta(\mathsf{p}(\Norm_{k''/k'}(a)))$, is the pushforward
under $\tilde{\rho}_{\alpha}$ of $\diff_{\beta}(s,\nu)$. In case B it
will be the measure on $S\backslash\phi^{-1}(g)$ invariant under the
action of $\mu_{q+1}$ and whose restriction to $X_{1}$ is the
pushforward under $\tilde{\rho}_{w_{0}}$ of $\diff_{\beta}(s,\nu)$.

First, let $G$ be the group of $k$-points of an algebraic group
defined over $k$ and $\mathfrak{g}$ its Lie algebra. Let $\diff g$ be
a left Haar measure on $G$ and $\diff X$ be a Haar measure on the
additive group $\mathfrak{g}$. Let also $U\subset\mathfrak{g}$ an open
subset such that the exponential map $\exp$ is defined on $U$,
$\exp(U)$ is open in $G$ and $\exp$ induces a diffeomorphism from $U$
onto $\exp(U)$. Then the change of variable formula applied to $\exp$
shows that, if $U$ is sufficiently small, there exists a number $c>0$,
independent of the choice of $U$, such that
\begin{equation}\label{eq49}
\int_{\exp(U)}F(g)\diff(g)=c\int_{U}F\circ\exp(X)\diff X,
\end{equation}
for every continuous function $F$ with compact support on $\exp(U)$.
We say that the Haar measures $\diff g$ and $\diff X$ are tangent if
$c=1$.

Recall that in case B, $\varpi'=\varpi''$.

\begin{lem}\label{lem5.16}
Let $\diff s$ (resp. $\diff X$) be the Haar measure on $S'$
(resp. $\mathfrak{s}'$) normalized such that $\int_{S'_{1}}\diff s=1$
(resp. $\int_{\mathfrak{s}'\cap\mathcal{O}'\varpi''}\diff
X=1$). Then these Haar measures are tangent.
\end{lem}
\begin{dem}
It follows easily from Corollary \ref{co2.1} that for $m\in\mathbb{N}$
big enough $\exp$ induces a diffeomorphism from
$\varpi''^{m}\mathcal{O}'\cap\mathfrak{s}'$ onto $S'_{m}=S'\cap
T_{m}$. Let $c>0$ and $m$ big enough be such that equation
\eqref{eq49} is satisfied for $U=\varpi''^{m}\mathcal{O}'\cap\mathfrak{s}'$.

But, Lemma \ref{lem5.3} in case A and Corollary \ref{co5.1} in case B also
apply to $S'$. It follows that $\vert S'_{1}/S'_{m}\vert=
\vert\mathcal\varpi''{O}'\cap\mathfrak{s}'/
\mathcal\varpi''^{m}{O}'\cap\mathfrak{s}'\vert$,
$m\geq1$, and thus, taking $m$ big enough,
\begin{eqnarray*}
  1=\int_{S'_{1}}\diff s&=&\vert S'_{1}/S'_{m}\vert\int_{S'_{m}}\diff
  s\\
  &=&c\vert\mathcal\varpi''{O}'\cap\mathfrak{s}'/
  \mathcal\varpi''^{m}{O}'\cap\mathfrak{s}'\vert
  \int_{\varpi''^{m}\mathcal{O}'\cap\mathfrak{s}'}\diff
  X=c.
\end{eqnarray*}
Thus the lemma.
\end{dem}

Here we explain how to choose a Haar measure $\diff x$ on $k$ which is
canonical with respect to the choice of  $\psi$. Let
$\mathcal{C}^{\infty}_{c}$ be the space of complex valued locally
constant functions with compact support on $k$. We define the Fourier
transforms $\mathcal{F}$ and $\overline{\mathcal{F}}$ of
$F\in\mathcal{C}^{\infty}_{c}$ as
\begin{eqnarray*}
  \mathcal{F}F(y)&=&\int_{k}F(x)\psi(xy)\diff x,\\
  \overline{\mathcal{F}}F(y)&=&\int_{k}F(x)\psi(-xy)\diff x.
\end{eqnarray*}
Then the canonical Haar measure $\diff x$ on $k$ is the one
normalized in such a way that we have
\begin{equation*}
\overline{\mathcal{F}}\circ\mathcal{F}F=F\mbox{,
}F\in\mathcal{C}^{\infty}_{c}.
\end{equation*}
We easily check that the Haar measure $\diff x$ is the one satisfying
\begin{equation}\label{eq50}
\int_{\mathcal{O}}\diff x=q^{\frac{\lambda_{\psi}}{2}}.
\end{equation}

Now let $V$ be a vector-space over $k$ equipped with a symplectic form
$\beta$.  Let $(e_{1},\ldots,e_{2r})$ be a symplectic basis of $V$,
that is which satisfies $\beta(e_{l},e_{m})=\beta(e_{r+l},e_{r+m})=0$
and $\beta(e_{l},e_{r+m})=\delta_{l,m}$, $1\leq l,m\leq r$. Let
$\varrho$ be the linear isomorphism from $k^{2r}$ onto $V$ such that
\begin{equation*}
  \varrho(x_{1},\ldots,x_{2r})=
  \sum_{1\leq i\leq2r}x_{i}e_{i}.
\end{equation*}
Then we equip $V$ with the Haar measure which is the pushforward by
$\varrho$ of the product measure $\diff x_{1}\cdots\diff x_{2r}$ on
$k^{2r}$. This measure does not depend upon the choice of the
symplectic basis and, as each measure $\diff x_{i}$, $1\leq i\leq2r$,
is canonically associated to the chosen character $\psi$ of $k$. We
denote it by $\diff_{\beta} v$.

According to \cite[pages 9-10]{popa-2011}, the volume form
$\wedge^{\dim(\mathfrak{s}')}\tilde{\beta}$ induces on
$S'\times\mathcal{O}'\cap(\mathfrak{s}'/\varpi'')$
(resp. $S'\times\mathcal{O}'\cap\varpi'^{j-\mu-2v''(u)}(\mathfrak{s}'/u)$)
a Haar measure $\diff_{\beta}(s',\eta)$ which is tangent to the
measure $\diff_{\tilde{\beta}_{(1,0)}}(X,\eta)$ on
$\mathfrak{s}'\times \mathfrak{s}'/\varpi''$ (resp
$\mathfrak{s}'\times\varpi'^{j-\mu-2v''(u)}(\mathfrak{s}'/u)$).

Recall that $\diff s$ is the Haar measure on $S'$ such that
$\int_{S'_{1}}\diff s=1$. In case A, let $\diff x$
(resp. $\diff\eta$) be the Haar measure on
$\mathcal{O}'\cap(\mathfrak{s}'/\varpi'')$
(resp. $\mathfrak{s}'/\varpi''$) such that
$\int_{\mathcal{O}'\cap(\mathfrak{s}'/\varpi'')}\diff x=1$
(resp. $\int_{\mathcal{O}'\cap(\mathfrak{s}'/\varpi'')}\diff\nu=1$). In
case B, let $\diff x$ (resp. $\diff\eta$) be the Haar measure on
$\mathcal{O}'\cap\varpi'^{j-\mu-2v''(u)}(\mathfrak{s}'/u)$
(resp. $\varpi'^{j-\mu-2v''(u)}(\mathfrak{s}'/u)$) such that
$\int_{\mathcal{O}'\cap\varpi'^{j-\mu-2v''(u)}(\mathfrak{s}'/u)}\diff x=1$
(resp. $\int_{\mathcal{O}'\cap\varpi'^{j-\mu-2v''(u)}(\mathfrak{s}'/u)}\diff\eta=1$).

Then the Haar measure $\diff s\diff x$ on
$S'\times\mathcal{O}'\cap(\mathfrak{s}'/\varpi'')$
(resp. $S'\times\mathcal{O}'\cap\varpi'^{j-\mu-2v''(u)}(\mathfrak{s}'/u)$)
is tangent to $\diff X\diff\eta$.

\begin{lem}\label{lem5.17}
In case A, we have
\begin{equation*}
\diff_{\beta}(s,\eta)=q^{\dim(\mathfrak{s}')(j-1)}\diff s\diff x.
\end{equation*}

In case B, we have
\begin{equation*}
  \diff_{\beta}(s,\eta)=
  q^{[\frac{j-1}{2e}]\vert I'\vert
    +\vert\{i\in I'\vert 1\leq i\leq l-1\}\vert}
    \diff s\diff x,
\end{equation*}
where $l=j-1-2e[\frac{j-1}{2e}]$.
\end{lem}
\begin{dem}
  Case A).

Let us put $d'=\dim(\mathfrak{s}')$ and let
$e_{1},\ldots,e_{d'}\in\mathfrak{s}'$ such that
$(e_{1}/\varpi'',$ $\ldots,e_{d'}/\varpi'')$ is a basis of
$\mathcal{O}'\cap(\mathfrak{s}'/\varpi'')$ over $\mathcal{O}$. Then
let $e_{d'+1},\ldots,$ $e_{2d'}\in k'\cap(\mathfrak{s}'/\varpi'')$ be
such that $\langle
e_{l}/\varpi'',e_{d'+m}\rangle_{k'/k}=\delta_{l,m}$, $1\leq l,m\leq
d'$. This is possible since the restriction of $\langle\,
,\,\rangle_{k'/k}$ to
$\mathfrak{s}'/\varpi''=(\mathfrak{s}/\varpi'')^{\perp}$ is
non-degenerate. Then we claim that $(e_{d'+1},\ldots,e_{2d'})$, which
is a basis of $k'\cap(\mathfrak{s}'/\varpi'')$ over $k$, is also a
basis of $\mathcal{O}'\cap(\mathfrak{s}'/\varpi'')$ over
$\mathcal{O}$.

On one hand, if $1\leq l\leq d'$, we readily see that $\langle
e_{d'+l},x\rangle_{k'/k}\in\mathcal{O}$ for every $x\in\mathcal{O}'=
\mathcal{O}'\cap(\mathfrak{s}/\varpi'')
\oplus\mathcal{O}'\cap(\mathfrak{s}'/\varpi'')$ which, as $k'$ is
unramified over $k$, implies that $e_{d'+l}\in\mathcal{O}'$. On the
other hand, let $x\in\mathcal{O}'\cap(\mathfrak{s}'/\varpi'')$ and
write $x=\sum_{1\leq l\leq d'} x_{l}e_{d'+l}$ with $x_{1},\dots,
x_{d'}\in k$. Then, for $1\leq l\leq d'$, we have $x_{l}=\langle
x,e_{l}/\varpi''\rangle_{k'/k}\in\mathcal{O}$ showing that
$\mathcal{O}'\cap(\mathfrak{s}'/\varpi'')=\sum_{1\leq l\leq
  d'}\mathcal{O}e_{d'+l}$. Thus the claim.

Let us put $\tilde{e}_{l}=(e_{l},0)$ and
$\tilde{e}_{d'+l}=(0,e_{d'+l})$, $1\leq l\leq d'$. Then it follows
from \eqref{eq47} that, for $1\leq l,m\leq d'$, we have
\begin{equation*}
  \tilde{\beta}_{1,0}(\tilde{e}_{l},\tilde{e}_{d'+m})=
  \frac{(-1)^{\mu-j}}{2}\varpi^{\mu+2-j}\delta_{l,m}.
\end{equation*}

As $k'$ and $k'\varpi''$ are Lagrangian subspaces of $k''$ with respect to
$\beta$, we see that
\begin{equation*}
(\tilde{e}_{1},\ldots,\tilde{e}_{d'},
2(-1)^{\mu-j}\varpi^{j-\mu-2}\tilde{e}_{d'+1},
  \ldots,2(-1)^{\mu-j}\varpi^{j-\mu-2}\tilde{e}_{2d'})
\end{equation*}
is a symplectic basis of $\mathfrak{s}'\times \mathfrak{s}'/\varpi''$.

Let $\varrho:k^{2d'}\rightarrow\mathfrak{s}'\times
\mathfrak{s}'/\varpi''$ be the $k$-linear isomorphism such that
\begin{equation*}
\varrho(x_{1},\ldots,x_{2d'})=
(\sum_{1\leq l\leq d'}x_{l}\tilde{e}_{l},
2(-1)^{\mu-j}\varpi^{j-\mu-2}\sum_{1\leq l\leq d'}x_{d'+l}\tilde{e}_{d'+l}).
\end{equation*}
Then by definition of $\diff_{\tilde{\beta}_{1,0}}(X,\nu)$, we have
\begin{eqnarray*}
  \int_{\mathcal{O}'\varpi''\cap\mathfrak{s}'\times\mathcal{O}'\cap(\mathfrak{s}'/\varpi'')}
  \diff_{\tilde{\beta}_{1,0}}(X,\nu)&=&
  \int_{\varrho^{-1}(\mathcal{O}'\varpi''\cap\mathfrak{s}'\times
    \mathcal{O}'\cap(\mathfrak{s}'/\varpi''))}
  \diff x_{1}\cdots\diff x_{2d'}\\
 &=&\left(\int_{\mathcal{O}}\diff x\times
 \int_{\varpi^{\mu+2-j}\mathcal{O}}\diff x\right)^{d'}\\
 &=&q^{-(\mu+2-j)d'}\int_{\mathcal{O}^{2d'}}\diff x_{1}\cdots\diff x_{2d'}\\
 &=&q^{d'(j-1)}.
\end{eqnarray*}
This shows that $\diff_{\tilde{\beta}_{1,0}}(X,\nu)=q^{d'(j-1)}\diff X\diff\eta$. 
Thus the result in case A.

Case B)

Let us put $e_{i}=(\frac{1}{e\varpi}\varpi'^{i}\nu,0)$ %, $i\in
%I'\backslash\{0\}$, $e_{0}=(\frac{1}{e}\nu,0)$ if $0\in I'$,
and
$f_{i}=(0,\varpi'^{j-\mu-v''(u)+2e-i})$, $i\in I'$.  Then $(e_{i},i\in
I',f_{i}, i\in I')$ is a symplectic basis of $\mathfrak{s}'\times
\varpi'^{j-\mu-2v''(u)}(\mathfrak{s}'/u)$ with respect to
$\tilde{\beta}_{(1,0)}$.

Moreover, one has
\begin{equation*}
  \varpi'\mathcal{O}'\cap\mathfrak{s}'
  =\oplus_{i\in I'\cap\{0\}}\varpi\mathcal{O}\varpi'^{i}\nu
  \oplus
  \oplus_{i\in I'\smallsetminus\{0\}}\mathcal{O}\varpi'^{i}\nu.
\end{equation*}
Let us put $l=j-\mu-v''(u)-2e[\frac{j-\mu-v''(u)}{2e}]$. As
$\delta=2e-1$, we have $\mu=2e(\lambda_{\psi}-1)+1-v''(u)$ and thus
$l=j-1-2e[\frac{j-1}{2e}]$ while
$[\frac{j-\mu-v''(u)}{2e}]=[\frac{j-1}{2e}]+1-\lambda_{\psi}$. Then one has
\begin{equation*}
  \begin{split}
  \mathcal{O}'\cap\varpi'&^{j-\mu-2v''(u)}(\mathfrak{s}'/u)=
  \mathcal{O}'\cap\varpi'^{l}(\mathfrak{s}'/\nu)=
  \varpi^{-1}\varpi'^{l}(\varpi'^{2e-l}\mathcal{O}'\cap(\mathfrak{s}'/\nu))\\
  &\,\,\,\,\,\,\,\,\,\,\,\,\,\,\,=
  \varpi^{-1}\varpi'^{l}(\oplus_{i\in I', i\leq
    2e-l-1}\varpi\mathcal{O}\varpi'^{i}
   \oplus\oplus_{i\in I',
    2e-l\leq i\leq2e-1}\mathcal{O}\varpi'^{i})\\
  &\,\,\,\,\,\,\,\,\,\,\,\,\,\,\,=\oplus_{i\in I', i\leq
    2e-l-1}\mathcal{O}\varpi'^{l+i}
   \oplus\oplus_{i\in I',
     2e-l\leq i\leq2e-1}\varpi^{-1}\mathcal{O}\varpi'^{l+i},
\end{split}   
\end{equation*}
while $f_{i}=(0,\varpi^{[\frac{j-1}{2e}]+1-\lambda_{\psi}}\varpi'^{l+2e-i})$, $i\in I'$.

Let $\varrho:k^{I'}\times k^{I'}\rightarrow\mathfrak{s}'\times
\varpi'^{j-\mu-2v''(u)}\mathfrak{s}'/u$ be the $k$-linear isomorphism such that
\begin{equation*}
\varrho((x_{i},i\in I'),(y_{i},i\in I' ))=
\sum_{i\in I'}x_{i}e_{i}+
\sum_{i\in I'}y_{i}f_{i}.
\end{equation*}
Then by definition of $\diff_{\tilde{\beta}_{1,0}}(X,\nu)$, we have
\begin{multline*}
  \int_{\varpi'\mathcal{O}'\cap\mathfrak{s}'\times
    \mathcal{O}'\cap\varpi'^{j-\mu-2v''(u)}(\mathfrak{s}'/u)}
  \diff_{\tilde{\beta}_{1,0}}(X,\nu)=\\
  \int_{\varrho^{-1}(\varpi'\mathcal{O}'\cap\mathfrak{s}'\times
    \mathcal{O}'\cap\varpi'^{j-\mu-2v''(u)}(\mathfrak{s}'/u))}
  \prod_{i\in I'}\diff x_{i}\prod_{i\in I'}\diff y_{i}=\\
  \left(\int_{\varpi\mathcal{O}}\diff x\right)^{\vert I'\vert}\times
  \left(\int_{\varpi^{-[\frac{j-\mu-v''(u)}{2e}]-1}\mathcal{O}}\diff x
  \right)^{\vert\{i\in I'\vert 1\leq i\leq l-1\}\vert}\\
  \times
  \left(\int_{\varpi^{-[\frac{j-\mu-v''(u)}{2e}]}\mathcal{O}}\diff x
  \right)^{\vert I'\vert-\vert\{i\in I'\vert 1\leq i\leq l-1\}\vert}\\
  =q^{(\frac{\lambda_{\psi}}{2}-1)\vert I'\vert}\times
  q^{(\frac{\lambda_{\psi}}{2}+[\frac{j-\mu-v''(u)}{2e}]+1)
    \vert\{i\in I'\vert 1\leq i\leq l-1\}\vert}\\
  \times
 q^{(\frac{\lambda_{\psi}}{2}+[\frac{j-\mu-v''(u)}{2e}])(\vert I'\vert-
   \vert\{i\in I'\vert 1\leq i\leq l-1\}\vert)} \\
 =q^{[\frac{j-1}{2e}]\vert I'\vert
    +\vert\{i\in I'\vert 1\leq i\leq l-1\}\vert}.
\end{multline*}
The lemma follows.
\end{dem}

\begin{co}\label{co5.4}
In  case A, the volume of $S\backslash\phi^{-1}(g)$ is
\begin{equation*}
  \vol(S\backslash\phi^{-1}(g))=\int_{S\backslash\phi^{-1}(g)}\diff\mu_{\beta}
  =
  \varepsilon_{S}\vert\Theta(\mathsf{p}(\Norm_{k''/k'}(a)))\vert
  q^{\dim(\mathfrak{s}')(j-1)}.
\end{equation*}

In case B, the volume of $S\backslash\phi^{-1}(g)$ is
\begin{equation*}
  \vol(S\backslash\phi^{-1}(g))=\int_{S\backslash\phi^{-1}(g)}\diff\mu_{\beta}=
  \vert\mu_{q+1}/\mu_{q+1}\cap S\vert q^{[\frac{j-1}{2e}]\vert I'\vert
    +\vert\{i\in I'\vert 1\leq i\leq l-1\}\vert},
\end{equation*}
where $l=j-1-2e[\frac{j-1}{2e}]$.
\end{co}
\begin{dem}
  In case A, it follows from the definition of $\diff\mu_{\beta}$ and Lemma
  \ref{lem5.17} that
\begin{eqnarray*}
  \vol(S\backslash\phi^{-1}(g))&=&q^{\dim(\mathfrak{s}')(j-1)}
  \sum_{\alpha\in\Theta(\mathsf{p}(\Norm_{k''/k'}(a)))}
  \int_{S'\times\mathcal{O}'\cap(\mathfrak{s}'/\varpi'')}\diff s\diff x\\
  &=&q^{\dim(\mathfrak{s}')(j-1)}
  \vert S'/S'_{1}\vert\vert\Theta(\mathsf{p}(\Norm_{k''/k'}(a)))\vert.
\end{eqnarray*}
But according to Lemma \ref{lem5.14} one has $\vert
S'/S'_{1}\vert=\varepsilon_{S}$.

In case B, the result follows from the definition of
$\diff\mu_{\beta}$ and Lemma \ref{lem5.17}.
\end{dem}

Now we may state the main results of this section.

\begin{theo}\label{theo5.1}
We suppose that $k''$ is ramified over $k'$, $k'$ is unramified over
$k$ of order $2f$ with $f$ prime to $p$, $\varpi''^{2}=\varpi'$, and
$S$ is a proper and admissible subtorus of $T$. Let $\chi$ be a
non-trivial unitary character of $S$ and $2j$, $j\geq1$, be its
conductor.

Then $\chi$ appears in the restriction to
$S$ of the metaplectic representation if and only if there exists
$g\in\im\phi$ such that
%$a\in\mathcal{O}''^{\times}$ such that, putting $w=a\varpi''^{\mu-j}$
%and $g=\phi(w)\in\mathfrak{s}^{*}$ one has 
\begin{equation}\label{eq51}
\chi(\rho_{[\frac{j}{2}]}(\eta))=\psi(\langle g,\eta\rangle)\mbox{,
}\eta\in\varpi^{[\frac{j}{2}]}\mathcal{O}'\varpi''\cap\mathfrak{s}.
\end{equation}
Then one has $g=\phi(w)$ with $w\in k''$ such that $v''(w)=\mu-j$.

Moreover, the multiplicity $m(\chi)$ of $\chi$ in the restriction to
$S$ of the metaplectic representation is given by
\begin{equation}\label{eq52}
  m(\chi)=\vol(S\backslash\phi^{-1}(g))=
  \varepsilon_{S}\vert\Theta(\mathsf{p}(\Norm_{k''/k'}(a)))\vert
  q^{\dim(\mathfrak{s}')(j-1)}.
\end{equation}
\end{theo}
\begin{dem}
The first assertion of the theorem follows from Proposition
\ref{pr5.2} (i) \eqref{eq30}.

The remaining assertions follow from Proposition \ref{pr5.2},
Corollary \ref{co5.3} and the fact that
$\dim(\mathfrak{s}')=\codim(\mathfrak{s})$.
\end{dem}

\begin{theo}\label{theo5.2}
We suppose that $k''$ is unramified over $k'$ which is totally
and tamely ramified over $k$ of degree $2e$, and $S$ is a proper
admissible subtorus of $T$.

Let $\chi$ be a unitary character of $S$ with conductor $j\geq1$. Then,
$\chi$ appears in the restriction to $S$ of the metaplectic
representation if and only if $\mu-j$ is even.This is equivalent to
the fact that there exists $g$ a linear form on $\mathfrak{s}$
belonging to $\im\phi\smallsetminus\{0\}$ such that
\begin{equation}\label{eq58}
  \chi(\rho_{[\frac{j+1}{2}]}(\eta))=\psi(\langle g,\eta\rangle)\mbox{, }
  \eta\in\varpi'^{[\frac{j+1}{2}]}\mathcal{O}'\nu\cap\mathfrak{s}.
\end{equation}
Then one has $f=\phi(w)$ with $w\in k''$ such that $v''(w)=\frac{\mu-j}{2}$.

Moreover, the multiplicity $m(\chi)$ of $\chi$ in the restriction to
$S$ of the metaplectic representation is given by
\begin{equation}\label{eq59}
  m(\chi)=\vol(S\backslash\phi^{-1}(g))=\vert\mu_{q+1}/\mu_{q+1}\cap S\vert
  q^{\alpha(j)}.
\end{equation}
\end{theo}
\begin{dem}
It follows from Lemma \ref{lem5.12} and Proposition \ref{pr5.3} that
$\chi$ appears in the restriction to $S$ of the metaplectic
representation if and only if $j$ has the same parity than $\mu$ and
that there exists $a\in\mathcal{O}''^{\times}$ such that putting
$w=\varpi'^{\frac{\mu-j}{2}}a$, equation \eqref{eq58} is satisfied
with $g=\phi(w)$.

Recall that for $j\in\mathbb{N}$, the number $\alpha(j)$ has been defined in
formula \eqref{eq39.5}.
Then, according to Corollary \ref{co5.4} and Proposition \ref{pr5.3}
to prove the theorem we have to show that
\begin{equation*}
\vert\{1,\ldots,j-2\}\cap(2e\mathbb{N}+I')\vert=
      [\frac{j-1}{2e}]\vert I'\vert+\vert\{i\in I'\vert1\leq i\leq
      l-1\}\vert,
\end{equation*}
where $l=j-1-2e[\frac{j-1}{2e}]$.

If $m=j-2-2e[\frac{j-2}{2e}]$, one has
\begin{equation}\label{eq61}
  \vert\{1,\ldots,j-2\}\cap(2e\mathbb{N}+I')\vert=
  [\frac{j-2}{2e}]\vert I'\vert+\vert\{i\in I'\vert 1\leq i\leq m\}\vert.
\end{equation}

If $j-1$ is not divisible by $2e$, one has
$[\frac{j-2}{2e}]=[\frac{j-1}{2e}]$ and $m=l-1$. The result follows
from \eqref{eq61}.

If $j-1$ is a multiple of $2e$, then 
$[\frac{j-2}{2e}]=[\frac{j-1}{2e}]-1$, $m=2e-1$ and $l=0$. It follows that
$[\frac{j-2}{2e}]\vert I'\vert+\vert\{i\in I'\vert 1\leq i\leq m\}\vert=
[\frac{j-1}{2e}]\vert I'\vert$, thus the result.
\end{dem}
%\subsection{}\label{5.10}

\begin{appendices}
\section{Appendix}\label{apA}

For the content of this section, we are indebted to François
Courtès$^{\dag}$.

Let $\Gamma$ be a finite cyclic group of order $f$ and $\mathbb{Z}[\Gamma]$
be its group algebra over $\mathbb{Z}$. We choose $\sigma\in\Gamma$ a
generator.

We consider the group algebra $\mathbb{Z}[\Gamma]$ as a module over
itself, acting by multiplication on the left. In the following we will
describe the submodules $M$ of $\mathbb{Z}[\Gamma]$ such that the
quotient $\mathbb{Z}[\Gamma]/M$ is a torsion-free $\mathbb{Z}$-module
and are minimal for this property.

Given a divisor $d$ of $f$ we denote by $\Phi_{d}(X)$ the $d$-th
cyclotomic polynomial and we put
$P_{f,d}(X)=\frac{X^{f}-1}{\Phi_{d}(X)}$,
$Q_{d}(X)=\frac{X^{d}-1}{\Phi_{d}(X)}$ and
$H_{d}(X)=\frac{X^{f}-1}{X^{d}-1}$. These are polynomials with integer
coefficients.

Let $p:\mathbb{Z}/f\mathbb{Z}\rightarrow\mathbb{Z}/d\mathbb{Z}$ be the
natural projection, namely $p : x \pmod{f}\mapsto x\pmod{d}$. For
$l\in\mathbb{Z}/d\mathbb{Z}$, we consider the element of
$\mathbb{Z}[\Gamma]$
\begin{equation*}
Y_{d,l}=\sum_{m\in p^{-1}(l)}\sigma^{m}.
\end{equation*}
Then we have $Y_{d,0}=H_{d}(\sigma)$, $Y_{d,l}=Y_{d,0}\sigma^{l}$,
$l\in\mathbb{Z}/d\mathbb{Z}$ and these elements are linearly
independent over $\mathbb{Z}$.

\begin{defi}\label{defA1}
We denote by $M_{d}$ the
$\mathbb{Z}$-submodule of $\mathbb{Z}[\Gamma]$ whose elements are of
the form $\sum_{l=0}^{d-1}a_{l}Y_{d,l}$ such that
$\sum_{l=0}^{d-1}a_{l}X^{l}$ is a polynomial in $\mathbb{Z}[X]$
divisible by $Q_{d}(X)$
\end{defi}
We see at once that the set
$\{P_{f,d}(\sigma)\sigma^{l}\vert 0\leq l\leq\varphi(d)-1\}$, where
$\varphi$ is the Euler phi or totient function, is a basis of $M_{d}$ over
$\mathbb{Z}$ and $\rank_{\mathbb{Z}}(M_{d})=\varphi(d)$.

\begin{theo}\label{theoap}
  Keep the above notations.

  (i) For $d$ a divisor of $f$, $M_{d}$ is a
  $\mathbb{Z}[\Gamma]$-submodule of $\mathbb{Z}[\Gamma]$ such that
  $\mathbb{Z}[\Gamma]/M_{d}$ is torsion free over $\mathbb{Z}$ and
  minimal for this property. Moreover, the eigenvalues of $\sigma$ in
  $M_{d,\mathbb{C}}=M_{d}\otimes_{\mathbb{Z}}\mathbb{C}$ are the
  $d$-th primitive roots of unity.

  (ii) If $M$ is a $\mathbb{Z}[\Gamma]$-submodule of
  $\mathbb{Z}[\Gamma]$ such that $\mathbb{Z}[\Gamma]/M$ is torsion
  free over $\mathbb{Z}$ and minimal for this property, there exists a
  divisor $d$ of $f$ such that $M=M_{d}$.
\end{theo}
\begin{dem}
(i) We consider
  $\tilde{M}_{d}=\oplus_{l\in\mathbb{Z}/d\mathbb{Z}}\mathbb{Z}Y_{d,l}$. Clearly
  $\tilde{M}_{d}$ is a $\mathbb{Z}[\Gamma]$-submodule of
  $\mathbb{Z}[\Gamma]$ on which $\sigma^{d}$ acts trivially.

  Let $\zeta$ be a $d$-th primitive root of unity in
  $\mathbb{C}$. Then $\mathbb{Q}(\zeta)$ is a Galois extension of
  order $\varphi(d)$ of $\mathbb{Q}$ whose Galois group $G$ is isomorphic
  to $(\mathbb{Z}/d\mathbb{Z})^{\times}$ the group of invertible
  elements in $\mathbb{Z}/d\mathbb{Z}$ : the element
  $u\in(\mathbb{Z}/d\mathbb{Z})^{\times}$ correspond to the element
  $\tau_{u}$ of $G$ such that $\tau_{u}(\zeta)=\zeta^{u}$.

  Into $\tilde{M}_{d}\otimes_{\mathbb{Z}}\mathbb{Q}(\zeta)\subset
  \mathbb{Q}(\zeta)[\Gamma]$ we consider the elements
  \begin{equation}\label{eqap1}
Z_{d,u}=\sum_{l\in\mathbb{Z}/d\mathbb{Z}}\zeta^{-ul}Y_{d,l}\mbox{,
}u\in\mathbb{Z}/d\mathbb{Z}.
\end{equation}
They form a basis of eigenvectors of the $\mathbb{Q}(\zeta)$
vector-space
$\tilde{M}_{d}\otimes_{\mathbb{Z}}\mathbb{Q}(\zeta)$ for the action
of $\sigma$ : $\sigma.Z_{d,u}=\zeta^{u}Z_{d,u}$,
$u\in\mathbb{Z}/d\mathbb{Z}$.

The group $G$ acts naturally on $\mathbb{Q}(\zeta)[\Gamma]$, and one has
\begin{equation*}
  \tau_{v}(Z_{d,u})=Z_{d,uv}\mbox{, }u\in\mathbb{Z}/d\mathbb{Z}\mbox{, }
  v\in(\mathbb{Z}/d\mathbb{Z})^{\times}.
\end{equation*}

Let
$W_{d}\subset\tilde{M}_{d}\otimes_{\mathbb{Z}}\mathbb{Q}(\zeta)$ be
the  subspace generated over $\mathbb{Q}(\zeta)$ by the $Z_{d,u}$,
$u\in(\mathbb{Z}/d\mathbb{Z})^{\times}$. Then $W_{d}$ is invariant
under $G$ and $\sigma$ and we will see that
$W_{d}\cap\mathbb{Z}[\Gamma]=M_{d}$.

We have $W_{d}\cap\mathbb{Q}[\Gamma]=W_{d}^{G}:=\{Z\in
W_{d}\vert\tau_{u}Z=Z\mbox{, }u\in(\mathbb{Z}/d\mathbb{Z})^{\times}\}$.

Let $Z\in W_{d}$ and write
$Z=\sum_{u\in(\mathbb{Z}/d\mathbb{Z})^{\times}}a_{u}Z_{d,u}$ with
$a_{u}\in\mathbb{Q}(\zeta)$. Then one has
$\tau_{v}(Z)=\sum_{u\in(\mathbb{Z}/d\mathbb{Z})^{\times}}\tau_{v}(a_{u})Z_{d,uv}$,
$v\in(\mathbb{Z}/d\mathbb{Z})^{\times}$. It follows that
\begin{equation*}
  W_{d}\cap\mathbb{Q}[\Gamma]=
  \bigg\{\sum_{u\in(\mathbb{Z}/d\mathbb{Z})^{\times}}\tau_{u}(a)Z_{d,u}
  \vert a\in\mathbb{Q}(\zeta)\bigg\}.
\end{equation*}

As the $Z_{d,u}$, $u\in(\mathbb{Z}/d\mathbb{Z})^{\times}$, are linearly
independent, the map
$a\mapsto\sum_{u\in(\mathbb{Z}/d\mathbb{Z})^{\times}}\tau_{u}(a)Z_{d,u}$
is a bijection from $\mathbb{Q}(\zeta)$ onto
$W_{d}\cap\mathbb{Q}[\Gamma]$. Thus
$\dim_{\mathbb{Q}}W_{d}\cap\mathbb{Q}[\Gamma]=\varphi(d)$.

Now let
$Z=\sum_{u\in(\mathbb{Z}/d\mathbb{Z})^{\times}}\tau_{u}(a)Z_{d,u}$
where $a\in\mathbb{Q}(\zeta)$ belonging to
$W_{d}\cap\mathbb{Q}[\Gamma]$ and write $a=\sum_{0\leq
  v\leq\varphi(d)-1}a_{v}\zeta^{v}$ with $a_{v}\in\mathbb{Q}$. Then a
straightforward computation shows that $Z=\sum_{0\leq
  v\leq\varphi(d)-1}a_{v}\tilde{Z}_{d,v}$ with
\begin{eqnarray*}
  \tilde{Z}_{d,v}&=&\sum_{u\in(\mathbb{Z}/d\mathbb{Z})^{\times}}\zeta^{uv}Z_{d,u}\\
  &=&\sum_{l\in\mathbb{Z}/d\mathbb{Z}}\big(\sum_{u\in(\mathbb{Z}/d\mathbb{Z})^{\times}}
  \zeta^{u(v-l)}\big)Y_{d,l}.
\end{eqnarray*}
For $0\leq v\leq\varphi(d)-1$, let $Q_{d,v}$ be the polynomial
\begin{equation*}
Q_{d,v}(X)=\sum_{l\in\mathbb{Z}/d\mathbb{Z}}\big(\sum_{u\in(\mathbb{Z}/d\mathbb{Z})^{\times}}
  \zeta^{u(v-l)}\big)X^{l}.
\end{equation*}
It follows easily from the fundamental theorem of symmetric
polynomials and from the fact that the $\zeta^{u}$,
$u\in(\mathbb{Z}/d\mathbb{Z})^{\times}$ are the roots of the
cyclotomic polynomial $\Phi_{d}\in\mathbb{Z}[X]$, that
$Q_{d,v}\in\mathbb{Z}[X]$, $0\leq v\leq\varphi(d)-1$.

Now, let $\xi$ be a root of $Q_{d}$, that is $\xi=\zeta^{s}$ with $s$
not prime to $d$ or equivalently
$s\in\mathbb{Z}/d\mathbb{Z}\smallsetminus(\mathbb{Z}/d\mathbb{Z})^{\times}$,
and let $0\leq v\leq\varphi(d)-1$. Then one has
\begin{eqnarray*}
  Q_{d,v}(\xi)&=&\sum_{l\in\mathbb{Z}/d\mathbb{Z}}\big
  (\sum_{u\in(\mathbb{Z}/d\mathbb{Z})^{\times}}\zeta^{uv+(s-u)l}\big)\\
  &=&\sum_{u\in(\mathbb{Z}/d\mathbb{Z})^{\times}}\zeta^{uv}
  \big(\sum_{l\in\mathbb{Z}/d\mathbb{Z}}\zeta^{(s-u)l}\big)\\
  &=&0,
\end{eqnarray*}
because $u-s\neq0$ for $u\in(\mathbb{Z}/d\mathbb{Z})^{\times}$ and
$s\in\mathbb{Z}/d\mathbb{Z}\smallsetminus(\mathbb{Z}/d\mathbb{Z})^{\times}$.

It follows that $Q_{d,v}$ is divisible by $Q_{d}$. Then it is clear that
$W_{d}^{G}=M_{d}\otimes_{\mathbb{Z}}\mathbb{Q}$.

In order to prove that $W_{d}\cap\mathbb{Z}[\Gamma]=M_{d}$, we must
show that for any $Z=\sum_{0\leq
  l\leq\varphi(d)-1}a_{l}Q_{d}(\sigma)Y_{d,0}\sigma^{l}\in
M_{d}\otimes_{\mathbb{Z}}\mathbb{Q}\cap\mathbb{Z}[\Gamma]$ all
coefficients $a_{l}$ are in $\mathbb{Z}$. If this is not the case, let
$0\leq l_{0}\leq\varphi(d)-1$ be the biggest integer such that
$a_{l_{0}}\notin\mathbb{Z}$. Then $Z_{0}=\sum_{0\leq l\leq
  l_{0}}a_{l}Q_{d}(\sigma)Y_{d,0}\sigma^{l}\in
M_{d}\otimes_{\mathbb{Z}}\mathbb{Q}\cap\mathbb{Z}[\Gamma]$ and one has
$Z_{0}=a_{l_{0}}\sigma^{f-\varphi(d)+l_{0}}+\sum_{0\leq
  l<f-\varphi(d)+l_{0}}b_{l}\sigma^{l}$. It follows that
$a_{l_{0}}\in\mathbb{Z}$, a contradiction.

Now let $Z\in\mathbb{Z}[\Gamma]$ and $a\in\mathbb{Z}\smallsetminus\{0\}$
such that $aZ\in M_{d}$. Then $Z\in
M_{d}\otimes_{\mathbb{Z}}\mathbb{Q}\cap\mathbb{Z}[\Gamma]$, and by the
above $Z\in M_{d}$. This shows that $\mathbb{Z}[\Gamma]/M_{d}$ is
torsion free over $\mathbb{Z}$. The fact that $M_{d}$ is minimal for
this property will follow from (ii) : if $M_{d}$ was not minimal, it
would strictly contain one such submodule which according to (ii) would
be a certain $M_{d'}$. But, considering the eigenvalues of $\sigma$ in
$M_{d,\mathbb{C}}$ and $M_{d',\mathbb{C}}$, one must have $d=d'$ and
$M_{d'}=M_{d}$, a contradiction.

(ii) Let $M$ be a $\mathbb{Z}[\Gamma]$-submodule of
$\mathbb{Z}[\Gamma]$ such that $\mathbb{Z}[\Gamma]/M$ is torsion free
over $\mathbb{Z}$ and minimal for this property. Then the vector-space
$M_{\mathbb{C}}=M\otimes_{\mathbb{Z}}\mathbb{C}$ %admits a basis of
%eigenvectors for $\sigma$ with eigenvalues roots of unity of order $f$. 
contains an eigenvector for $\sigma$ with eigenvalue $\zeta$ a root of
unity of order $f$. Let $d$ be the divisor of $f$ such that $\zeta$ is
a $d$-th primitive root of unity. Then clearly one has
$M_{d,\mathbb{C}}\subset M_{\mathbb{C}}$. As $\mathbb{Z}[\Gamma]/M$ is
torsion free over $\mathbb{Z}$, one has
$M=M_{\mathbb{C}}\cap\mathbb{Z}[\Gamma]$ and thus $M_{d}\subset
M$. And by minimality of $M$, $M=M_{d}$.
\end{dem}

\begin{co}\label{coap1}
The $M_{d}$ for $d$ a divisor of $f$ are $\mathbb{Z}$-submodules of
$\mathbb{Z}[\Gamma]$ which are pairwise orthogonal for the non-degenerate
bilinear form $\langle\, ,\,\rangle$ such that
$\langle\sigma,\tau\rangle=\delta_{\sigma,\tau}$,
$\sigma,\tau\in\Gamma$. Moreover, their direct summand is a
$\mathbb{Z}$-submodule of $\mathbb{Z}[\Gamma]$ of same rank.
\end{co}
\begin{dem}
Let $\zeta$ be a $f$-th primitive root of $1$ in $\mathbb{C}$. Then
the action of $\sigma$ diagonalizes in $\mathbb{Q}(\zeta)[\Gamma]$
with eigenvalues $\zeta^{u}$, $u\in\mathbb{Z}/f\mathbb{Z}$ each of
multiplicity one. Moreover, the $\mathbb{Q}(\zeta)$-hermitian form $h$
on $\mathbb{Q}(\zeta)[\Gamma]$ such that
$h(\sigma,\tau)=\delta_{\sigma,\tau}$, $\sigma,\tau\in\Gamma$ is
clearly $\Gamma$-invariant and its restriction to $\mathbb{Z}[\Gamma]$
is the bilinear form $\langle\, ,\,\rangle$.  It then follows easily
that $\mathbb{Q}(\zeta)[\Gamma]$ is the orthogonal direct sum of the
$\mathbb{Q}(\zeta)\otimes_{\mathbb{Z}}M_{d}$ for $d$ a divisor of
$f$. Thus the result.
\end{dem}
\begin{defi}\label{defA2}
If $\underline{d}$ is a set of divisors of $f$, we put
$M_{\underline{d}}=\oplus_{d\in\underline{d}}M_{d}$ and
$\overline{M}_{\underline{d}}=\{\chi\in\mathbb{Z}[\Gamma']\vert\exists
a\in\mathbb{Z}\smallsetminus\{0\}\mbox{ such that } a\chi\in
M_{\underline{d}}\}$.
\end{defi}
We remark that $M_{d}=M_{\{d\}}=\overline{M}_{\{d\}}$.

\begin{co}\label{coap2}
  The $\mathbb{Z}[\Gamma]$-submodules $M$ of $\mathbb{Z}[\Gamma]$ such
  that $\mathbb{Z}[\Gamma]/M$ is torsion free are the
  $\overline{M}_{\underline{d}}$ for $\underline{d}$ a set of divisors
  of $f$.
\end{co}
\begin{dem}
By construction any $\overline{M}_{\underline{d}}$ is a
$\mathbb{Z}[\Gamma]$-submodule of $\mathbb{Z}[\Gamma]$ and
$\mathbb{Z}[\Gamma]/\overline{M}_{\underline{d}}$ is torsion free.

Conversely let $M$ be $\mathbb{Z}[\Gamma]$-submodule of
$\mathbb{Z}[\Gamma]$ such that $\mathbb{Z}[\Gamma]/M$ is torsion
free. Let $\zeta$ be a $f$-th primitive root of $1$ in
$\mathbb{C}$. One checks easily that
$\mathbb{Q}(\zeta)\otimes_{\mathbb{Z}}M$ is the direct sum of the
$\mathbb{Q}(\zeta)\otimes_{\mathbb{Z}}M_{d}$ it contains. The corollary
follows easily.
\end{dem}

\begin{co}\label{coap3}
  The order of any non-zero element in the torsion module
  $\mathbb{Z}[\Gamma]/(\oplus_{d\vert f}M_{d})$ is a divisor of $f$.
\end{co}
\begin{dem}
We use the notations introduced in the proofs of Theorem \ref{theoap}
and Corollary \ref{coap1} and consider the $\Gamma$-modules
$\mathbb{Z}[\zeta][\Gamma]$ and
$M_{d}[\zeta]=\mathbb{Z}[\zeta]\otimes_{\mathbb{Z}}M_{d}$. It is clear
that $M_{d}=M_{d}[\zeta]\cap\mathbb{Z}[\Gamma]$.  Recall
  the hermitian form $h$ introduced in the proof of Corollary
  \ref{coap1}. Then, the $Z_{d,u}$,
  $u\in(\mathbb{Z}/d\mathbb{Z})^{\times}$ are linearly independent
  elements of $M_{d}[\zeta]$ and one checks easily that they satisfy
  $h(Z_{d,u},Z_{d,v})=f\delta_{u,v}$. Let us denote by $M'_{d}[\zeta]$
  the $\mathbb{Z}[\zeta]$-submodule of $M_{d}[\zeta]$ generated by
  these vectors. Then, the $M'_{d}[\zeta]$ are clearly pairwise
  orthogonal with respect to $h$.

Now let $x\in\mathbb{Z}[\Gamma]$. Then one has $x=\sum_{d\vert
  f}\sum_{u\in(\mathbb{Z}/d\mathbb{Z})^{\times}}a_{d,u}Z_{d,u}$ with
$a_{d,u}\in\mathbb{Q}[\zeta]$. Taking scalar products, one has
$\langle x,Z_{d,u}\rangle=fa_{d,u}\in\mathbb{Z}[\zeta]$, $d\vert f$,
$u\in(\mathbb{Z}/d\mathbb{Z})^{\times}$. It follows that
$fx\in
(\oplus_{d\vert f}M'_{d}[\zeta])\cap
\mathbb{Z}[\Gamma]\subset\oplus_{d\vert f}M_{d}$.
Thus the result.
\end{dem}

\begin{co}\label{coap4}
If $\underline{d}$ is a set of divisors of $f$, the order of any
non-zero element in the torsion module $\mathbb{Z}$-module
$\overline{M}_{\underline{d}}/M_{\underline{d}}$ is a divisor of $f$.
\end{co}
\begin{dem}
It suffices to prove that
$\overline{M}_{\underline{d}}\cap(\oplus_{d\vert
  f}M_{d})=M_{\underline{d}}$ : it will follow that
$\overline{M}_{\underline{d}}/M_{\underline{d}}$ is a submodule of
$\mathbb{Z}[\Gamma]/(\oplus_{d\vert f}M_{d})$ and the corollary will
be a consequence of the previous one.

Thus, let $\chi\in \overline{M}_{\underline{d}}\cap(\oplus_{d\vert
  f}M_{d})$. Write $\chi=\sum_{d\vert f}\chi_{d}$
with $\chi_{d}\in M_{d}$. By hypothesis, there exists
$a\in\mathbb{Z}\smallsetminus\{0\}$ such that $a\chi=a\sum_{d\vert
  f}\chi_{d}\in\sum_{d\in\underline{d}}M_{d}$ which implies
$\chi_{d}=0$ for $d\notin\underline{d}$. The result follows.
\end{dem}

\section{Appendix}\label{apB}

In this section, we prove that there exist anisotropic tori in $Sp(W)$
which are not contained in an elliptic maximal torus.

Let $\ell$ and $\ell'$ two Lagrangians subspaces of $W$ such that
$W=\ell\oplus\ell'$. As we have seen in paragraph \ref{1.4}, the
symplectic form $\beta$ induces a duality between $\ell$ and $\ell'$
and the map
\begin{equation}\label{eqapb1}
a\mapsto\begin{pmatrix} a&0\\0&^{t}a^{-1}\end{pmatrix}
\end{equation}
is an embedding of $GL(\ell)$ into $Sp(W)$.
\begin{lem}\label{lemapb1}
Le $T\subset GL(\ell)$ be an anisotropic torus. We suppose that the
action of $T$ in $\ell$ is irreducible over $k$. Then, if the image of
$T$ by the embedding \eqref{eqapb1} is contained in an elliptic
maximal torus of $Sp(W)$ the $T$ modules $\ell$ and $\ell'$ are
isomorphic.
\end{lem}
\begin{dem}
Let $V\subset W$ be a non-zero $T$-invariant symplectic subspace in
which $T$ acts irreducibly. Let $p:W\rightarrow\ell$ be the linear
projection parallel to $\ell'$. Then $p$ intertwines the
actions of $T$ on $V$ and $\ell$, and $\ker p_{\vert
  V}=V\cap\ell'=0$. As $\ell$ is an irreducible $T$-module, $p_{\vert
  V}$ is an isomorphism of $T$-modules from $V$ onto $\ell$. Exchanging
the role of $\ell$ and $\ell'$ we see that the $T$-modules $V$ and
$\ell'$ are also isomorphic, thus the lemma.
\end{dem}

\begin{pr}\label{prapb1}
Le $k'$ be a finite extension of $k$ and let $T\subset k'^{\times}$ be an
anisotropic subtorus. We equip $k'$ with the structure of $T$-module
induced by the multiplication and suppose that $k'$ is an irreducible
$T$-module over $k$. We equip the dual space $k'^{*}$ of $k'$ over $k$
with the contragredient structure of $T$-module. Then the following
assertions are equivalent :

(i) the $T$-modules $k'$ and $k'^{*}$ are isomorphic over $k$

(ii) there exists $\tau$, an automorphism of order two of $k'$ over
$k$, such that $T\subset\{x\in k'^{\times}\vert x\tau(x)=1\}$.
\end{pr}
\begin{dem}
Let us denote by $\gamma$ the action of $T$ into $k'$ and by
$\gamma^{*}$ the contragredient action : $\gamma(t)x=tx$, $t\in T$,
$x\in k'$ and $\gamma^{*}(t)=\,^{t}\gamma(t)^{-1}$, $t\in T$.

The non-degenerate $k$-bilinear form $(x,y)=\tr_{k'/k}xy$, $x,y\in k'$
induces a natural $k$-isomorphism between $k'$ and $k'^{*}$. We
identify $k'$ and $k'^{*}$ via this isomorphism. Then, we see that
$\gamma^{*}(t)x=t^{-1}x$, $t\in T$, $x\in k'$.

Suppose that (i) is satisfied and let $A:k'\rightarrow k'$ be a $k$-linear
isomorphism intertwining $\gamma$ and $\gamma^{*}$. Let us put
\begin{equation*}
\tau(x)=(A1)^{-1}Ax\mbox{, }x\in k'.
\end{equation*}
Then, we have\begin{equation}\label{eqapb2}
\tau(tx)=(A1)^{-1}Atx=(A1)^{-1}t^{-1}Ax=t^{-1}\tau(x)\mbox{, }t\in
T\mbox{, } x\in k'.
\end{equation}
Applying this formula with $x=1$, we get
\begin{equation}\label{eqapb3}
t\tau(t)=1\mbox{, }t\in T,
\end{equation}
while applying again $\tau$ to \eqref{eqapb2} we get
$\tau^{2}(tx)=\tau(t^{-1}\tau(x))=t\tau^{2}(x)$, $t\in T$, $x\in k'$.
As $k'$ is an irreducible $T$-module, $k'$ is generated by $T$ as a
$k$-linear subspace. It follows that the centralizer of $T$ in the
algebra of $k$-linear endomorphisms of $k'$ is also the centralizer of
$k'$ and is thus equal to $k'$.  Then, as $\tau^{2}(1)=1$, we see that
$\tau^{2}=Id_{k'}$, that is $\tau$ is a $k$-linear endomorphism of
order two.

In order to prove that (i) $\Rightarrow$ (ii), it remains to show that
$\tau(xy)=\tau(x)\tau(y)$, $x,y\in k'$. But, this follows easily from the fact that $k'$ is generated by $T$ over $k$.

%Let $x,y\in k'$. As $k'$ is
%generated by $T$ over $k$, we may write $x=\sum_{1\leq i\leq
%  m}\lambda_{i}s_{i}$ and $y=\sum_{1\leq j\leq n}\mu_{j}t_{j}$ where
%$s_{i},t_{j}\in T$ and $\lambda_{i},\mu_{j}\in k$. Then one has
%\begin{eqnarray*}
%\tau(x)\tau(y)
%&=&(\sum_{1\leq i\leq m}\lambda_{i}s_{i}^{-1})(\sum_{1\leq j\leq n}\mu_{j}t_{j}%^{-1})\\
%&=&\sum_{1\leq i\leq m,1\leq j\leq n}\lambda_{i}\mu_{j}\tau(s_{i}t_{j})\\
%&=&\tau(\sum_{1\leq i\leq m,1\leq j\leq n}\lambda_{i}\mu_{j}s_{i}t_{j})\\
%&=&\tau(xy)
%\end{eqnarray*}
%Thus the desired implication.

The proof of the converse implication is straightforward. 
\end{dem}

\begin{co}\label{coapb1}
There exist anisotropic tori in $Sp(W)$ which are not contained in
any elliptic maximal torus in $Sp(W)$.
\end{co}
\begin{dem}
Let $k'$ be an extension of $k$ of degree at least $3$ and take
$W=k'\oplus k'$ with $\beta$ such that
$\beta((x,y),(x',y'))=\tr_{k'/k}(xy'-x'y)$, $x,x',y,y'\in k'$. Then
the two copies of $k'$ are Lagrangian subspaces in $W$. Now let
$T=\{x\in k'\vert \Norm_{k'/k}(x)=1\}$. It is clear that $k'$ is an
irreducible $T$-module for the multiplication. We consider $T$ as a
subgroup of $GL(k')$ and embed it in $Sp(W)$ via the morphism
\eqref{eqapb1}. It follows from Lemma \ref{lemapb1} and Proposition
\ref{prapb1} that $T$ is not contained in any elliptic maximal torus
in $Sp(W)$.

One also sees that, if the degree of $k'$ over $k$ is odd, no non-trivial
anisotropic subtorus of $k'^{\times}$ is contained in an elliptic
maximal torus.
\end{dem}

\end{appendices}

\bibliographystyle{amsplain} \bibliography{biblio}

\end{document}